\def\DATE{{Feb.\ 11, 2008}}
\def\DATE{\relax}
\magnification=1100
\baselineskip=14truept
\voffset=.75in  
\hoffset=1truein
\hsize=4.5truein
\vsize=7.75truein
\parindent=.166666in
\pretolerance=500 \tolerance=1000 \brokenpenalty=5000

\footline={\hfill{\rm\the\pageno}\hfill\llap{\sevenrm\DATE}}

\def\note#1{%
  \hfuzz=50pt%
  \vadjust{%
    \setbox1=\vtop{%
      \hsize 3cm\parindent=0pt\eightpoints\baselineskip=9pt%
      \rightskip=4mm plus 4mm\raggedright#1%
      }%
    \hbox{\kern-4cm\smash{\box1}\hfil\par}%
    }%
  \hfuzz=0pt
  }
\def\note#1{\relax}

\def\anote#1#2#3{\smash{\kern#1in{\raise#2in\hbox{#3}}}%
  \nointerlineskip}     

\newcount\equanumber
\equanumber=0
\newcount\sectionnumber
\sectionnumber=0
\newcount\subsectionnumber
\subsectionnumber=0
\newcount\snumber  
\snumber=0

\def\section#1{%
  \subsectionnumber=0%
  \snumber=0%
  \equanumber=0%
  \advance\sectionnumber by 1%
  \noindent{\bf \the\sectionnumber .~#1.~}%
}%
\def\subsection#1{%
  \advance\subsectionnumber by 1%
  \snumber=0%
  \equanumber=0%
  \noindent{\bf \the\sectionnumber .\the\subsectionnumber .~#1.~}%
}%
\def\prevs{\the\sectionnumber .\the\subsectionnumber .\the\snumber }
\long\def\Definition#1{\bigskip\noindent{\bf Definition.}\quad{\it#1}}
\long\def\Claim#1{\noindent{\bf Claim.}\quad{\it#1}}
\long\def\Corollary#1{%
  \global\advance\snumber by 1%
  \bigskip
  \noindent{\bf Corollary~\prevs .}%
  \quad{\it#1}%
}%
\long\def\Lemma#1{%
  \global\advance\snumber by 1%
  \bigskip
  \noindent{\bf Lemma~\prevs .}%
  \quad{\it#1}%
}%
\def\Proof{\noindent{\bf Proof.~}}
\long\def\Proposition#1{%
  \advance\snumber by 1%
  \bigskip
  \noindent{\bf Proposition~\prevs .}%
  \quad{\it#1}%
}%
\long\def\Remark#1{%
  \bigskip
  \noindent{\bf Remark.~}#1%
}%
\long\def\Theorem#1{%
  \advance\snumber by 1%
  \bigskip
  \noindent{\bf Theorem~\prevs .}%
  \quad{\it#1}%
}%
\long\def\Statement#1{%
  \advance\snumber by 1%
  \bigskip
  \noindent{\bf Statement~\prevs .}%
  \quad{\it#1}%
}%
\def\ifundefined#1{\expandafter\ifx\csname#1\endcsname\relax}
\def\labeldef#1{\global\expandafter\edef\csname#1\endcsname{\prevs}}
\def\labelref#1{\expandafter\csname#1\endcsname}
\def\label#1{\ifundefined{#1}\labeldef{#1}\note{$<$#1$>$}\else\labelref{#1}\fi}

\def\preveq{(\the\sectionnumber .\the\subsectionnumber .\the\equanumber)}
\def\neq{\global\advance\equanumber by 1\eqno{\preveq}}

\def\ifundefined#1{\expandafter\ifx\csname#1\endcsname\relax}

\def\equadef#1{\global\advance\equanumber by 1%
  \global\expandafter\edef\csname#1\endcsname{\preveq}%
  \preveq}

\def\equaref#1{\expandafter\csname#1\endcsname}
 
\def\equa#1{%
  \ifundefined{#1}%
    \equadef{#1}%
  \else\equaref{#1}\fi}

\font\eightrm=cmr8%
\font\sixrm=cmr6%

\font\eightsl=cmsl8%

\font\eightbf=cmb8%

\font\eighti=cmmi8%
\font\sixi=cmmi6%

\font\eightsy=cmsy8%
\font\sixsy=cmsy6%

\font\eightex=cmex8%
\font\sixex=cmex6%
\font\fiveex=cmex5%

\font\eightit=cmti8%

\font\eighttt=cmtt8%

\font\tenbb=msbm10%
\font\eightbb=msbm8%
\font\sevenbb=msbm7%
\font\sixbb=msbm6%
\font\fivebb=msbm5%
\newfam\bbfam  \textfont\bbfam=\tenbb  \scriptfont\bbfam=\sevenbb  \scriptscriptfont\bbfam=\fivebb%

\font\tenbbm=bbm10

\font\tencmssi=cmssi10%
\font\sevencmssi=cmssi7%
\font\fivecmssi=cmssi5%
\newfam\ssfam  \textfont\ssfam=\tencmssi  \scriptfont\ssfam=\sevencmssi  \scriptscriptfont\ssfam=\fivecmssi%

\font\tenfrak=cmfrak10%
\font\eightfrak=cmfrak8%
\font\sevenfrak=cmfrak7%
\font\sixfrak=cmfrak6%
\font\fivefrak=cmfrak5%
\newfam\frakfam  \textfont\frakfam=\tenfrak  \scriptfont\frakfam=\sevenfrak  \scriptscriptfont\frakfam=\fivefrak%
\def\frak{\fam\frakfam\tenfrak}%

\font\tenmsam=msam10%
\font\eightmsam=msam8%
\font\sevenmsam=msam7%
\font\sixmsam=msam6%
\font\fivemsam=msam5%

\def\bb{\fam\bbfam\tenbb}%

\def\hexdigit#1{\ifnum#1<10 \number#1\else%
  \ifnum#1=10 A\else\ifnum#1=11 B\else\ifnum#1=12 C\else%
  \ifnum#1=13 D\else\ifnum#1=14 E\else\ifnum#1=15 F\fi%
  \fi\fi\fi\fi\fi\fi}
\newfam\msamfam  \textfont\msamfam=\tenmsam  \scriptfont\msamfam=\sevenmsam  \scriptscriptfont\msamfam=\fivemsam%
\def\msam{\msamfam\tenmsam}%
\mathchardef\leq"3\hexdigit\msamfam 36%
\mathchardef\geq"3\hexdigit\msamfam 3E%

\def\eightpoints{%
\def\rm{\fam0\eightrm}%
\textfont0=\eightrm   \scriptfont0=\sixrm   \scriptscriptfont0=\fiverm%
\textfont1=\eighti    \scriptfont1=\sixi    \scriptscriptfont1=\fivei%
\textfont2=\eightsy   \scriptfont2=\sixsy   \scriptscriptfont2=\fivesy%
\textfont3=\eightex   \scriptfont3=\sixex   \scriptscriptfont3=\fiveex%
\textfont\itfam=\eightit  \def\it{\fam\itfam\eightit}%
\textfont\slfam=\eightsl  \def\sl{\fam\slfam\eightsl}%
\textfont\ttfam=\eighttt  \def\tt{\fam\ttfam\eighttt}%
\textfont\bffam=\eightbf  \def\bf{\fam\bffam\eightbf}%

\textfont\frakfam=\eightfrak  \scriptfont\frakfam=\sixfrak \scriptscriptfont\frakfam=\fivefrak  \def\frak{\fam\frakfam\eightfrak}%
\textfont\bbfam=\eightbb      \scriptfont\bbfam=\sixbb     \scriptscriptfont\bbfam=\fivebb      \def\bb{\fam\bbfam\eightbb}%
\textfont\msamfam=\eightmsam  \scriptfont\msamfam=\sixmsam \scriptscriptfont\msamfam=\fivemsam  \def\msam{\msamfam\eightmsam}

\rm%
}

\def\poorBold#1{\setbox1=\hbox{#1}\wd1=0pt\copy1\hskip.25pt\box1\hskip .25pt#1}

\mathchardef\lsim"3\hexdigit\msamfam 2E%
\mathchardef\gsim"3\hexdigit\msamfam 26%

\def\bcup{{\textstyle\bigcup}}
\def\bcap{\,{\textstyle\bigcap}\,}
\def\complement{{\rm c}}
\def\d{\,{\rm d}}
\def\D{{\rm D}}
\def\ds{\displaystyle}
\long\def\DoNotPrint#1{\relax}
 
\def\finetune#1{#1}
\def\fixedref#1{#1\note{fixedref$\{$#1$\}$}}
\def\g#1{g_{[0,#1)}}
\def\Id{\hbox{\rm Id}}
\def\limeps{\lim_{\epsilon\to 0}}
\def\liminft{\liminf_{t\to\infty}}
\def\limn{\lim_{n\to\infty}}
\def\limsupt{\limsup_{t\to\infty}}
\def\limt{\lim_{t\to\infty}}
\def\oF{\overline F{}}
\def\oH{\overline H{}}
\def\qed{{\vrule height .9ex width .8ex depth -.1ex}}
\def\ss{\scriptstyle}

\def\calL{{\cal L}}
\def\calM{{\cal M}}
\def\calN{{\cal N}}

\def\II{\hbox{\tenbbm 1}}
\def\NN{{\bb N}\kern .5pt}
\def\RR{{\bb R}}
\def\ZZ{{\bb Z}}

\def\fS{{\frak S}}

\pageno=1

\centerline{\bf VERAVERBEKE'S THEOREM AT LARGE}
\vskip 8pt
\centerline{\bf\poorBold{$\widetilde{\hbox to 1cm{\hfill}}$}}
{\eightpoints
\centerline{\bf ON THE MAXIMUM OF SOME PROCESSES}
\centerline{\bf  WITH NEGATIVE DRIFT}
\centerline{\bf AND HEAVY TAIL INNOVATIONS}
}
\bigskip
 
\centerline{Ph.\ Barbe$^{(1)}$ and W.P.\ McCormick$^{(2)}$}
\centerline{${}^{(1)}$CNRS, France, ${}^{(2)}$University of Georgia}
 
{\narrower
\baselineskip=9pt\parindent=0pt\eightpoints

\bigskip

{\bf Abstract.}  Veraverbeke's (1977) theorem relates the tail of the
distribution of the supremum of a
random walk with negative drift to the tail of the distribution of its
increments, or equivalently, the probability that a centered random walk with 
heavy-tail increments hits a moving linear boundary. We study similar problems
for more general processes. In particular, we derive an analogue of
Veraverbeke's theorem for fractional integrated ARMA models without 
prehistoric influence, when 
the innovations have regularly varying tails. Furthermore, we prove some limit
theorems for the trajectory of the process, conditionally on a large maximum.
Those results are obtained by using a general scheme of proof which we
present in some detail and should be of value in other related problems.

\bigskip


\bigskip
 

}

\vfill\eject
{\eightpoints
\obeylines\baselineskip=9pt\setbox1=\hbox{1. }
\newdimen\dd
\dd=\wd1
\def\quad{\hskip\dd}
{\bf \noindent Contents}

\medskip

{\sectionnumber=0
\def\section#1{\advance\sectionnumber by 1\subsectionnumber=1\line{\the\sectionnumber. #1\hfill}}
\def\subsection#1{\line{\quad\the\sectionnumber.\the\subsectionnumber. #1\hfill}\advance\subsectionnumber by 1}

\section{Introduction}
\section{$(g,F)$-processes, their maximum and sample paths}
\subsection{$(g,F)$-processes and their maximum}%
\subsection{Removing the positivity of the coefficients}%
\subsection{Typical trajectories leading to a large maximum}%
\subsection{Examples}
\subsection{Note on the quantiles of the maximum of $(g,F)$-processes}%
\subsection{Concluding remarks}%
\section{Veraverbeke's theorem at large}%
\subsection{The single large jump heuristic}%
\subsection{From the heuristic to a theorem}%
\subsection{Analysis of the paths leading to a large maximum}%
\section{A large deviation inequality and a Karamata type theorem}
\subsection{A large deviation inequality}%
\subsection{A Karamata type theorem}%
\section{Some asymptotic analysis related to analytic functions}%
\subsection{Preliminaries}%
\subsection{The functions $\Psi_n$ and their inverses}%
\subsection{Approximation of $\Psi_n^{-1}$ when $(g_n)_{n\geq 0}$ tends to $0$}%
\subsection{Asymptotic analysis of $\Psi_n^{-1}$ when $(g_n)_{n\geq 0}$ tends to infinity}%
\subsection{Approximation of $\Psi_n^{-1}$ when $(g_n)_{n\geq 0}$ has a positive and finite limit}
\section{Proof of the results of section 2}%
\subsection{Proof of Theorem 2.1.1 -- upper bound}%
\subsection{Proof of Theorem 2.1.1 -- lower bound}%
\subsection{Proof of Theorem 2.1.2 -- upper bound}%
\subsection{Proof of Theorem 2.1.2 -- lower bound}%
\subsection{Proof of Theorem 2.1.3}%
\subsection{Proof of Theorem 2.2.1}%
\subsection{Proof of Theorem 2.3.1}%

\noindent References
}
}

\vfill\eject

\def\preveq{(\the\sectionnumber .\the\equanumber)}
\def\prevs{\the\sectionnumber .\the\snumber }

\section{Introduction}
Veraverbeke's (1977) theorem relates the tail behavior of the maximum
of some random walks with negative drift to the tail behavior of their
increments. The purpose of this paper is to show that such relation
holds for a much larger class of discrete time stochastic processes
which encompass some nonstationary FARIMA ones.

Before going further, let us recall Veraverbeke's theorem. Let $X_i$,
$i\geq 1$, be a sequence of independent and identically distributed
random variables, having a negative mean $\mu$. Define the random walk $S_n$
by $S_0=0$ and $S_n=S_{n-1}+X_n$ for all positive $n$. Since the
increments $X_i$ have negative mean, the maximum of the walk,
$M=\max_{n\geq 0} S_n$, is almost surely finite.

Let $F$ be the common distribution function of the $X_i$, and let $\oF=1-F$
be its tail. This tail is regularly varying (see e.g.\ Bingham, Goldie,
Teugels, 1989) if there exists a nonnegative $\alpha$ such that
$$
  \lim_{t\to\infty} \oF(\lambda t)/\oF(t) = \lambda^{-\alpha}
$$
for any positive $\lambda$. The number $-\alpha$ is called the index
of regular variation.

Veraverbeke's theorem asserts that if the distribution of the increments
of the random walk has negative mean and regularly varying tail with index
$-\alpha$ less than $-1$, then
$$
  P\{\, M>t\,\} \sim {1\over -\mu} \int_t^\infty \oF(u)\d u \, ,
  \eqno{\equa{Vint}}
$$
as $t$ tends to infinity; or, equivalently, using Karamata's theorem,
$$
  P\{\, M>t\,\} \sim {1\over -\mu}\, {t\,\oF(t)\over \alpha-1}
  \eqno{\equa{VK}}
$$
as $t$ tends to infinity.

The original question which motivated us to write this paper is the following
simple one: if the increment of the `random walk' are themselves a random
walk, should we replace $t\oF(t)$ in \VK\ by $t^2\oF(t)$? The answer to 
this question is given in subsection 2.4 of this paper.

Veraverbeke's theorem has been extended in several directions. Clearly, one can
seek to prove that it holds for a larger class of
distributions for the increments. In that line of investigation, Veraverbeke's
(1977) original result asserts in fact that \VK\ holds whenever $F$ is 
subexponential. Later, Korshunov (1997) obtained necessary and sufficient 
condition for \Vint\ to hold, building upon Borovkov (1971), Pakes (1975),
Veraverbeke (1977), Embrechts and Veraverbeke (1982). We also mention that 
it is not necessary that the increments of the random
walk have negative mean for the maximum of the process to be almost surely
finite. An analysis of the tail distribution for the global maximum of a
random walk with heavy tail incremenents without mean and with a left tail
dominance can be found in Borovkov (2003). 

In a different direction,
Mikosch and Samorodnitsky (2000) replace the independent increments of
the random walk by an infinite order moving average process; they consider
$X_n=\mu+\sum_{j\in\ZZ}\phi_{n-j}\epsilon_j$, where $\mu$ is negative, the 
$\epsilon_j$ are centered, independent, equidistributed with a common
distribution having regularly varying tail. They assume further a natural tail
balance condition so that the distribution of $X_n$ is itself regularly 
varying even if some of the $\phi_j$ are negative. They prove that 
if $\sum_{j\in\ZZ} |j\phi_j|$ is finite, then
an analogue of Veraverbeke's theorem holds, that is, there exists a constant
$c$ such that
$$
  P\{\, M> t\,\}\sim c\, {t\oF(t)\over\alpha-1}
$$
as $t$ tends to infinity. Interestingly, the constant $c$ is explicit
--- but its value is irrelevant here --- and is, in general, different
than the factor $1/(-\mu)$ involved in \Vint\ or \VK. 
However, Mikosch
and Samorodnitsky's (2000) result shows that the decay in Veraverbeke's
theorem remains unchanged in their more general setting. They also
note that their assumptions exclude some fractional integrated ARMA
models. We remark that such models are considered in section 2.4 of
the current paper.

Yet in another direction, Konstantinides and Mikosch (2005) make a study
of the tail behavior of the global maximum of partial sums of dependent
heavy tailed summands with negative drift. More specifically, they
consider the stationary solution $Y_n$ of a stochastic recurrence
equation $Y_n=A_nY_{n-1}+B_n$ and consider the random walk $S_n=Y_1+\cdots +
Y_n$. Considering a negative real number $\mu$, they provide an 
asymptotic equivalent for the tail distribution
of the maximum of the process $(S_n-ES_n+\mu n)_{n\geq 1}$. One of their 
findings is that, with their assumptions on the coefficients $A_n$ and $B_n$,
the order of decay of the classical Veraverbeke result is preserved but the
constant in the asymptotic expression changes. They interpret this change
of constant as a measure of clustering of extremes for the stationary
$(Y_n)_{n\geq 1}$ sequence.

As indicated, a purpose of this paper is to describe the tail behavior 
of the maximum 
of some processes which generalize in a natural way the random walk model,
which are nonstationary and exhibit long range dependence --- for instance,
in the sense that the series of the correlations between the process 
at a fixed time and time $n$ is not summable in $n$. Of Veraverbeke's
original result, only that there is a relation between the tail of the
maximum and the tail of the innovation will be preserved; neither the
constant $1/(-\mu)$ nor, in contrast to the models studied by Mikosch and 
Samorodnitsky (2000) or Konstantinides and Mikosch (2005), the order of 
decay, $t\oF(t)$, will be preserved in general.

To understand the motivation for the class of processes which we are
going to introduce and to which we will extend Veraverbeke's result, 
as well as to frame the contribution of this paper in a larger context, 
it is necessary to recall two results on random walks and fractional
ARIMA processes; such presentation requires defining some notation related
to the latter. Toward this end, for any nonpositive integer $n$, we set
$X_n$ to be $0$. The backward shift $B$ acts on the sequence $(X_n)_{n\in\ZZ}$
by $BX_n=X_{n-1}$. As usual, this operator can be raised to a nonnegative 
power, with
$B^0$ being the identity and $B^n$ being defined inductively as $BB^{n-1}$. 
For any positive real number $d$ and any polynomials $\Phi$ and $\Theta$
with both $\Phi (1)$ and $\Theta (1)$ nonzero, the nonstationary 
$\hbox{FARIMA}(\Phi,d,\Theta)$
process $(Y_n)_{n\geq 0}$ with innovations $(X_n)_{n\geq 1}$ is defined 
by the formula
$$
  Y_n=(1-B)^{-d}\Phi(B)^{-1}\Theta(B) X_n \, ;
$$
the actual meaning of this expression is obtained by expanding the function
$$
  g(x)=(1-x)^{-d} \Phi(x)^{-1}\Theta(x)
  \eqno{\equa{gFARIMA}}
$$ 
as a Taylor series $\sum_{i\geq 0} g_i x^i$ and setting $Y_n=g(B)X_n$, that
is,
$$
  Y_n=\sum_{0\leq i<n} g_iX_{n-i} \, , \qquad n\geq 0 \, .
$$
Note that $Y_0$ vanishes. One sees that the random walk is obtained for
$d=1$ and the polynomials $\Phi$ and $\Theta$ being constant equal to $1$. 
The classical ARMA processes are obtained for $d=0$.

The two results alluded to --- which we will not use but put the present paper
in a broader perspective --- are that as $n$ tends to infinity, a random 
walk up to time $n$,
suitably rescaled and under the proper moment conditions on the increments
(the exact assumptions are irrelevant to this discussion)

\medskip

\noindent (i) converges in distribution to a Wiener process (see e.g.\ 
Billingsley, 1968);

\smallskip

\noindent (ii) obeys a large deviation principle with rate function involving
a derivative (Varadhan, 1966; see e.g.\ Dembo and Zeitouni, 1992).

\medskip

Extending those two results with possibly stronger assumptions on the 
increments, some FARIMA processes, properly rescaled

\medskip

\noindent (i) converge in distribution to a fractional Brownian motion (Akonom 
and Gouri\'eroux, 1987; see also Wang, Lin and Gulati, 2002)

\smallskip

\noindent (ii) obey a large deviation principle with rate function involving a 
fractional derivative (Barbe and Broniatowski, 1998). 

\medskip

Thus, at a broad level, underlying these results is the idea that some 
statements valid for the partial sum processes may be extended to
some FARIMA processes, replacing integrals or derivative by their
fractional analogue (see Oldham and Spanier, 2006, for fractional calculus).
This suggests that for some FARIMA processes, an analogue of Veraverbeke's
theorem might be true, replacing the integrated tail by a fractional
integrated tail.

\bigskip

There are further general motivations for results of this paper, which
are related to the disparate reasons for studying the maximum of random
walks with negative drifts and FARIMA processes. One area where the 
interest in this type
of result is clear is that of insurance risk. For the classical model
of claims arriving according to a homogeneous Poisson process and constant
premium rate, the surplus claim process viewed at lattice time points
forms such a random walk and the assumption of profitability ensures that
the increments have negative mean. The distribution of the global maximum
of the walk describes the ruin probability over an infinite horizon. This
was one of the motivation of Embrechts and Veraverbeke (1982). For dependent
heavy tail claims, asymptotic bounds on ruin probability have been given
by Nyrhinen (2005). A good reference for ruin probability calculations
under a wide variety of model assumptions can be found in Asmussen (2000).
From this perspective, our results allow calculation of ruin probability
when the claim process is a nonstationary FARIMA one.

In queuing theory, for a GI/G/1 queue with traffic intensity less than
$1$, the stationary distribution of the waiting time is given by the
distribution of the global maximum of a random walk with negative
drift and if the service time distribution has a heavy tail, we are
precisely in the situation governed by Veraverbeke's result; see Pakes
(1975) and Asmussen (1987). For extension of the theory to a dependent
setting, see Asmussen, Schmidli, Schmidt (1999). For related 
information in the case of
queuing networks, we refer to Baccelli and Foss (2004) and Baccelli,
Foss and Lelarge (2005). Again, our result could be converted into
statements on waiting time distribution for some queue.

From a modeling perspective, the processes which we will study extend
the FARIMA ones. FARIMA processes possess the desirable property that
both short-\ and long-memory components of a time series can be
accounted for. For example, in hydrology, Montanari, Rosso and Taqqu
(1997) use a FARIMA(1,$d$,1) process to model deseasonalized daily
flows into a lake. A value for the parameter $d$ in the range$0<d<1/2$
corresponds to a long-memory process. By way of illustration, a value
of $d=0.26$ was obtained for the lake inflow data studied by
Montanari, Rosso and Taqqu (1997), indicating that long memory
models are of value in this type of application. Further discussion
of applications of those models can be found in Samorodnitsky and
Taqqu (1992). For the estimation and theoretical properties of FARIMA
processes with heavy tails, we refer to Kokoszka and Taqqu
(1995). Resnick (2007) is also a source of information concerning
long-range dependence and heavy-tailed modeling.

FARIMA processes have also been of much use in econometric and time
series analysis, in part because their occurrence in aggregation of
light-tailed time series --- see Granger (1980) and the clear
exposition in Beran (1994) --- and also in connection with the problem
of testing for unit root (Akonom and Gouri\'eroux, 1987; Phillips, 1987; 
Tanaka, 1999).

\bigskip


\def\preveq{(\the\sectionnumber .\the\subsectionnumber .\the\equanumber)}
\def\prevs{\the\sectionnumber .\the\subsectionnumber .\the\snumber }

\section{\poorBold{$(g,F)$}-processes, their maximum and sample paths}\
This section contains our main concrete results. A more abstract formulation
is presented in section \fixedref{3}. In the first subsection
we define the $(g,F)$-processes, which generalize in a natural way the
FARIMA processes, and we state our tail equivalent of the distribution
of their supremum. This is done under some positivity assumption which
we remove in the second subsection. In the third subsection we analyze
the likely paths for such processes to reach a high level.
In the fourth subsection we discuss two examples.

\bigskip

\subsection{\poorBold{$(g,F)$}-processes and their maximum}%
Let $g$ be a real analytic function on the segment $(-1,1)$. 
Its Taylor series expansion
$$
  g(x)=\sum_{i\geq 0} g_i x^i
$$
allows one to define the nonstationary process $S_n=g(B)X_n$. We call such
a process a $(g,F)$-process, $F$ being the common distribution of the $X_n$
with $n$ positive, and with the convention that $X_n$ is $0$ if $n$ is 
nonpositive. Thus, $S_n=\sum_{0\leq i<n} g_i X_{n-i}$. We see that
if all the $g_i$ are equal to $1$, that is if $g(x)=1/(1-x)$, then $S_n$ is 
a random walk. In this
section we give an analogue of Veraverbeke's theorem for the maximum of
some $(g,F)$-processes.

Defining for any nonnegative integers $k$ and $n$ with
$k$ less than $n$,
$$
  g_{[k,n)} = \sum_{k\leq i<n} g_i \, ,
$$
and keeping the notation $\mu$ for the mean of the $X_j$,
the mean of $S_n$ is $\mu \g{n}$. For the maximum $M$ of this
process to be finite in a setting which extends that of a random walk, it is
natural to require $\mu$ to be negative and $\lim_{n\to\infty}\g{n}
=+\infty$. In particular, this latter requirement suggests that
$g$ should have a singularity at $1$.  Because we will need to have some 
estimation on the decay of the expectation of $S_n$ toward minus infinity,
because it is sufficient to encompass the FARIMA processes, and because it
yields a nice mathematical theory, we restrict the singularity of $g$
by requiring $g$ to be regularly varying at $1$, meaning the existence of some
$\gamma$ such that
$$
  \limeps {g(1-\lambda \epsilon)\over g(1-\epsilon)} = \lambda^{-\gamma} \, .
$$
For $g$ to have a singularity at $1$, it is then necessary that $\gamma$
is nonnegative. We will only consider positive $\gamma$.

If $(g_n)_{n\geq 0}$ is asymptotically equivalent to a monotone sequence 
as $n$ tends
to infinity, then a straightforward variant of Karamata's Tauberian theorem
for power series (Bingham, Goldie and Teugels, 1989, Corollary 1.7.3) shows
that regular variation of $g$ with index $\gamma$ is equivalent to that
of the sequence $(g_n)_{n\geq 0}$ with index $\gamma-1$; furthermore, in this
case, writing $\Gamma(\cdot)$ for the gamma function,
$$
  g_n\sim {g(1-1/n)\over n\Gamma(\gamma)}
  \eqno{\equa{gnEquiv}}
$$
as $n$ tends to infinity. Thus, when applying our results, either
the function $g$ or its coefficients can be given. 

Note that if the sequence $(g_n)_{n\geq 0}$ is regularly varying with index
$\gamma-1$, then, whenever $\gamma$ is positive and different than $1$,
the sequence $(g_n)_{n\geq0}$ is asymptotically equivalent to a
monotone sequence;
indeed, this follows from Bojanic and Seneta's theorem (Bingham, Goldie
and Teugels, 1989, Theorem 1.5.3). 
An alternative point of view, replacing
any assumption
on the coefficients by assumptions solely on the function $g$, is given by
Braaksma and Stark's (1997) complex variable analogue of Karamata's power
series theorem.

\bigskip

\noindent{\bf Note.} In the remainder of this section, whenever we use
an analytic function $g$, we assume that its Taylor 
coefficients at $0$ are nonnegative. This assumption will be dropped in 
section \fixedref{2.2}.

\bigskip

Since we are only interested in situations where the drift pushes 
the $(g,F)$-process toward minus infinity, it is natural to assume also that
$\limeps g(1-\epsilon)=+\infty$. Thus, if $g$ is regularly varying, there 
exists a regularly varying function $U$, of index $1/\gamma$, unique up 
to an asymptotic equivalence, such that
$$
  g\Bigl( 1-{1\over U(t)}\Bigr) \sim t
$$
as $t$ tends to infinity. This function $U$ appears in our results.

Our first result gives an asymptotic equivalent of the tail of the
distribution of $M$ when the sequence $(g_n)_{n\geq 0}$ converges to $0$.
In this case, the number
$$
  g^*=\sup_{i\geq 0} g_i
$$
is well defined, and in fact the supremum is even a maximum.

Recall the beta integral,
$$
  B(p,q)
  =\int_0^1 u^{p-1}(1-u)^{q-1}\d u
  =\int_0^\infty u^{p-1}(1+u)^{-p-q}\d u \, .
$$

\Theorem{\label{Thg}%
  Let $g$ be a real analytic function on $(-1,1)$ whose Taylor coefficients
  $(g_n)_{n\geq 0}$ are nonnegative, regularly varying of index $\gamma-1$, 
  not summable, and tend
  to $0$ at infinity. If $\gamma=1$, assume further that $(g_n)_{n\geq 0}$ is
  asymptotically equivalent to a monotone sequence.
  
  Let $F$ be a distribution function with negative mean and whose tail is 
  regularly varying with index $-\alpha$. If $1/\alpha<\gamma<1$
  then the maximum $M$ of the 
  corresponding $(g,F)$-process satisfies
  $$
    P\{\, M>t\,\} 
    \sim 
    {{g^*}^\alpha\over\gamma} B\Bigl({1\over\gamma},\alpha-{1\over\gamma}\Bigr)
    \Bigl( {\Gamma(1+\gamma)\over -\mu}\Bigr)^{1/\gamma} (U\oF) (t) \, ,
    \eqno{\equa{Vg}}
  $$
  as $t$ tends to infinity.
}

\bigskip

In the case of Theorem \Thg, the distribution tail of $M$ decays at
rate $U\oF$ which, because of \gnEquiv\ and the convergence of 
$(g_n)_{n\geq 0}$ to $0$, is slower than that in Veraverbeke's theorem.
Note that $U\oF$ is regularly varying of index $(1/\gamma)-\alpha$, which
can assume any value between $1-\alpha$ and $0$ by a proper choice
of $\gamma$.

\bigskip

We next consider the case where the sequence $(g_n)_{n\geq 0}$ tends
to infinity, which forces $\gamma$ to be at least $1$. This case is
more involved in particular when $\gamma$ is $1$. Indeed, in this latter case
our proof involves some more refined asymptotic analysis which requires 
more precise assumptions; a general result in this situation remains elusive.
We limit ourself to models neighboring the classical random walk and find
sufficient conditions for preserving the result obtained when $\gamma$ 
exceeds $1$. The analysis requires some extra notions related to the theory
of regularly varying functions, notions which we now introduce. 

Recall that a slowly
varying function $\ell$ has a Karamata representation (Bingham, Goldie and
Teugels, 1989, Theorem 1.3.1) which asserts the existence of a function
$\varepsilon(\,\cdot\,)$ with limit $0$ at infinity and a 
function $a(\,\cdot\,)$ with a finite positive limit at infinity, such that
$$
  \ell(x)=a(x)\exp \Bigl( \int_1^x{\varepsilon (u)\over u} \d u\Bigr) 
$$
ultimately. When $\gamma$ is equal to $1$, it is natural to assume that 
$$
  g_k\sim \ell(k)
  \eqno{\equa{HypGOneA}}
$$
where $\ell$ is slowly varying and tends to infinity at infinity. It then
follows from Corollary 1.3.5 in Bingham, Goldie and Teugels (1989) that
we can take $\varepsilon(\,\cdot\,)$ to be nonnegative. In fact, by the same
argument used to prove their Corollary 1.3.5, the 
function $\varepsilon(\,\cdot\,)$
can be taken positive if $\ell$ is increasing. Note that in the Karamata
representation of $\ell$, if the function $\varepsilon(\,\cdot\,)$ is equal to
$1/\log^p$ for some $p$ greater than $1$, then $\ell$ has 
a finite limit at infinity.
Thus, to fix the ideas, under some extra conditions which we will not assume, 
we could force $\varepsilon(\,\cdot\,)$ to tend to $0$ at a rate slower 
than $1/\log^p$ for any $p$ greater than $1$. The point of this remark is to
suggest that, actually, it is natural to assume 
that $\varepsilon(\,\cdot\,)$ 
is itself slowly varying. We will assume in fact a weak form of 
super-slow variation (see Bingham, Goldie and Teugels, 1989, \S 3.12.2)
for $\varepsilon(\,\cdot\,)$, namely that
$$
  \limt {\varepsilon\bigl (\lambda t\varepsilon(t)\bigr)\over\varepsilon(t)} 
  = 1
  \eqno{\equa{HypGOneB}}
$$
uniformly in $\lambda$ in any compact subset of the positive half-line. 
Following the terminology given in Bingham, Goldie and Teugels (1989, \S 2.11)
we will also assume further that
$$
  \log \varepsilon(e^t) \hbox{ is self-neglecting,}
  \eqno{\equa{HypGOneC}}
$$
meaning that the function $\varphi(t)=\log\varepsilon(e^t)$ satisfies
locally uniformly, that is, uniformly for $\lambda$ in any fixed compact set,
$$
  \limt{\varphi\bigl( t+\lambda\varphi(t)\bigr)\over\varphi(t)}=1 \, .
$$

These unfortunately technical looking conditions still allow a wide
array of interesting examples. For instance, if $g_k=\log^p k$ for some
positive $p$, then $\ell(x)=\log^p x$ and $\varepsilon (x)=p/\log x$ satistfies
both \HypGOneB\ and \HypGOneC. Similarly, $\ell(x)=\exp(\log^p x)$ with $p$
less than $1$ yields $\varepsilon (x)=p\log^{p-1}x$, and that latter function
satisfies \HypGOneB\ and \HypGOneC.

Note that the function \gFARIMA\ is regularly varying at $1$, with
index $d$.

For $\gamma$ positive we need to introduce the functions
$$
  \rho_\gamma(u)
  =\min_{y>0} {\Gamma(1+\gamma)+(u+y)^\gamma\over\gamma y^{\gamma-1}} \, ,
$$
whose argument $u$ is nonnegative. In general, it does not seem possible
to obtain an explicit form of the minimum. However, one
can easily see that $\rho_1(u)=1+u$ and $\rho_2(u)=u+\sqrt{2+u^2}$.
We will see that $\rho_\gamma^{-\alpha}$ is integrable on the positive 
half-line whenever $\alpha$ is greater than $1$ --- see 
Lemma \fixedref{5.4.2}.

The following result is an analogue of Veraverbeke's theorem for 
$(g,F)$-processes when the sequence $(g_n)_{n\geq 0}$ diverges to infinity.
Note that by \gnEquiv\ this condition necessitates that $\gamma$ is at least
$1$. We write $\Id$ for the identity function on the real line.

\Theorem{\label{ThG}
  Let $g$ be a real analytic function on $(-1,1)$ whose Taylor coefficients
  $(g_n)_{n\geq 0}$ tend to infinity and are regularly varying of 
  index $\gamma-1$. If $\gamma=1$ assume further that $(g_n)_{n\geq 0}$ is 
  asymptotically equivalent to a monotone sequence and that
  \HypGOneA, \HypGOneB\ and \HypGOneC\ hold.

  Let $F$ be a distribution function with negative mean and whose 
  tail is regularly varying with index $-\alpha$ less than $-1$.
  The maximum $M$ of the corresponding $(g,F)$-process satisfies
  $$
    P\{\, M>t\,\} 
    \sim 
    (-\mu)^{\alpha (1-1/\gamma)-1/\gamma} (\Id\oF)\circ U(t) \int_0^\infty
    \rho_\gamma^{-\alpha}(u)\d u
    \eqno{\equa{VG}}
  $$
  as $t$ tends to infinity.
}

\bigskip

When Theorem \ThG\ applies, the tail probability of $M$ decays 
like $(\Id\oF)\circ U$. Under the assumptions of Theorem \ThG, $U$ tends to 
infinity at a rate slower than that of the identity. Then, the rate of
decay of $(\Id\oF)\circ U$ is slower than that involved in Veraverbeke's
(1977) theorem. By a suitable choice of $U$, hence of $g$, the index of 
regular variation of $(\Id\oF)\circ U$ can assume any value between $1-\alpha$
and $0$.

\bigskip

The last result of this subsection somewhat fills the main gap left 
by Theorems \Thg\ and \ThG,
assuming now that the sequence $(g_n)_{n\geq 0}$ converges to a positive
and finite limit. Recall that $g^*$ is defined as the supremum of the sequence
$(g_n)_{n\geq 0}$.

\Theorem{\label{ThOne}%
  Let $g$ be a real analytic function on $(-1,1)$ whose Taylor coefficients
  $(g_n)_{n\geq 0}$ converge to a finite and positive limit $g_\infty$.
  
  Let $F$ be a distribution function with negative mean and whose tail is 
  regularly varying with index $-\alpha$ less than $-1$.  Then, the 
  maximum $M$ of the associated $(g,F)$-process satisfies
  $$
    P\{\, M>t\,\} \sim {{g^*}^\alpha\over-\mu g_\infty} 
    {t\oF(t)\over\alpha-1}
  $$
  as $t$ tends to infinity.
}

\bigskip

Under the assumptions of Theorem \ThOne, it is easy to see that 
$U(t)\sim t/g_\infty$ at infinity. Thus,
the result of Theorem \ThOne\ is formally that of Theorem
\Thg\ when $\gamma$ is $1$ and $U(t)\sim t/g_\infty$. This is no 
coincidence, and the proof of Theorem \ThOne\ builds upon that of Theorem \Thg.

Clearly, if $g(x)=1/(1-x)$, then all the Taylor coefficients $g_i$ are
equal to $1$, and Theorem \ThOne\ implies Veraverbeke's (1977) result.

It is interesting to compare the index
of regular variation of the tail 
probability of $M$ in Theorems \Thg, \ThG\ and \ThOne. We see that
it varies between $1-\alpha$ and $0$, and that within the class of
$(g,F)$-processes considered, the 
classical random walk of Verarverbeke's original result is an extreme
case where the index is equal to $1-\alpha$.

\bigskip

\subsection{Removing the positivity of the coefficients}%
For application to FARIMA processes without too restrictive assumptions
on the polynomials $\Theta$ and $\Phi$, it is necessary to extend the results
of the previous section to the case where some coefficients $g_i$ may be
negative. Our technique allows this extension. It is particularly simple if
only a finite number of $g_i$ are negative. It is still simple if
$\sum_{0\leq i<n}g_i$ is of the same order as $\sum_{0\leq i< n}|g_i|$ as
$n$ tends to infinity. However, the discussion is more delicate 
if $\sum_{0\leq i<n}g_i=o(\sum_{0\leq i< n}|g_i|)$ as $n$ tends to 
infinity --- this is related to the asymptotic behavior of the 
function $\psi_n$ to be defined in section \fixedref{3} and analyzed in some
detail in section \fixedref{4}. Since this latter case does not seem to
have much bearing to applications, we limit ourselves to the situation
where only a finite number of $g_i$ have negative sign. As we will see, this
is sufficient to cover the FARIMA processes.

If $\limn g_n=+\infty$, then it is easy to see from its proof that Theorem
\ThG\ remains valid without any change, even if finitely many coefficients
are negative.

If $\limn g_n=0$, then the lower tail of the distribution may play a role
in the tail behavior of $M$. The usual assumption is the following tail
balance condition. Let $F_*$ be the distribution function of $|X_1|$. The
tail balance condition asserts that
$$
  \limt {\oF(t)\over\oF_*(t)}\in (0,\infty)\qquad\hbox{and}\qquad
  \limt {F(-t)\over \oF_*(t)}\in [\, 0,\infty) \, .
  \eqno{\equa{TailBalance}}
$$
Writing $g_*$ for the smallest negative Taylor coefficient of $g$ if it exists
and for $0$ otherwise, that is
$$
  g_*=0\wedge\inf_{i\geq 0} g_i \, ,
$$
we have the following extension of Theorem \Thg. Recall that, by our 
convention, a regularly varying sequence is ultimately positive; in 
particular this forces $g^*$ to be positive.

\Theorem{\label{VgSigned}%
  Assume that $(g_n)_{n\geq 0}$ is regularly varying, not summable and tends
  to $0$ at infinity. Let $F$ be a distribution function with negative mean,
  satisfying the tail
  balance condition \TailBalance, and whose upper tail
  is regularly varying with index $-\alpha$. If $1/\alpha<\gamma<1$, 
  the maximum $M$ of the
  corresponding $(g,F)$-process satisfies
  $$\displaylines{\qquad
    P\{\, M>t\,\}
    \sim {1\over\gamma} B\Bigl( {1\over\gamma},\alpha-{1\over\gamma}\Bigr)
    \Bigl({\Gamma(1+\gamma)\over-\mu}\Bigr)^{1/\gamma} 
    \hfill\cr\noalign{\vskip 3pt}\hfill
    U(t)\bigl( {g^*}^\alpha \oF(t)+(-g_*)^\alpha F(-t)\bigr)
    \qquad\cr}
  $$
  as $t$ tends to infinity.
  }

\bigskip

Note that with regard to the previous theorem, writing $\widetilde F$ for 
the distribution function $x\mapsto \oF(-{x-})$ --- that is, if $X$ has 
distribution function $F$ then
$\widetilde F$ is the distribution function of ${-X}$ --- a $(g,F)$-process
has the same distribution as a $(-g,\widetilde F)$ one. However, the theorem
does not allow such a substitution for it is assumed that $F$ has negative
mean and $g$ tends to $+\infty$ at $1$.

\bigskip

\subsection{Typical trajectories leading to a large maximum}%
The purpose of this subsection is to describe the most likely
trajectories of $(g,F)$-processes leading to a large value of their
maximum. For simplicity, we consider only the case where all the $g_n$
are positive; the extension to the setting of the previous subsection
poses no real difficulties and does not seem to bring any further
understanding. Our result can also be viewed as an extension of those
of Asmussen and Kl\"uppelberg (1996).  For the classical random walk
with negative drift and heavy-tail increments, they prove that
conditionally on having the maximum of the process larger than $t$, the
process properly normalized, run through the proper time scale and up
to the time at which it reaches its maximum, converges to a straight
line with slope equal to the mean of the increments --- indicating
that the process behaves as expected in this time frame --- whereas
the overshoot at the jump time properly normalized converges to a
Pareto distribution. Further information on the time of jump and
conditional path behavior in this context is given in Asmussen
(2000). This section addresses a similar problem and carries the analysis
further by providing a description of $(g,F)$-processes both before
and after the jump time.

To analyse those trajectories, we write $N_t$ for the first passage 
time of the process $(S_n)_{n\geq 0}$ over the threshold $t$, that is
$$
  N_t=\min\{\, n \, :\, S_n>t\,\} \, ,
$$
with the convention that the minimum of the empty set is $+\infty$. Clearly,
$M$ exceeds $t$ if and only if $N_t$ is finite. The single large jump
heuristic described in detail in section \fixedref{3} suggests that $M$
exceeds $t$ because, most likely, one of the $X_i$, $1\leq i\leq N_t$, is
large. Thus, it is natural to consider the index $J_t$ of occurrence
of the `big jump', that is, the integer between $0$ and $N_t$ such that
$$
  X_{J_t}=\max\{\, X_i \, :\, 1\leq i\leq N_t\,\} \, .
$$
In case of ties, we take $J_t$ to be the smallest such index. Furthermore, 
recalling that
the function $U$ is defined by the asymptotic equivalence $g(1-1/U)\sim\Id$
at infinity, we consider the rescaled process 
$$
  \fS_t(\lambda)=S_{\lfloor \lambda U(t)\rfloor}/t \,
$$
as well as the rescaled random variables
$$
  \tau_t=J_t/U(t) \qquad\hbox{ and } \qquad Y_t=X_{J_t}/U(t) \, .
$$
The rescaled process $\fS_t$ is right continuous with left limit, and therefore
is viewed here in the space $\D[\,0,\infty)$ of all c\`adl\`ag functions
equipped with the Skorohod topology (see Billingsley, 1968; Pollard, 1984).
In order to obtain a pleasing result, we will assume that the tail balance 
condition \TailBalance\ holds. We will restrict ourselves to what
happens under the assumptions of Theorem \ThG. 

\Theorem{\label{TrajG}%
  Let $\gamma$ be greater than $1$ and assume that the hypotheses of Theorem
  \ThG\ hold as well as the tail balance condition \TailBalance. If $\mu=-1$,
  then the conditional distribution of $(\fS_t,\tau_t,Y_t)$ conditional on
  $M>t$ converges weakly$*$
  to that of $\bigl(\fS,\tau,\rho_\gamma(\tau)Y\bigr)$
  where $Y$ and $\tau$ are independent and
  
  \smallskip

  \noindent (i) $Y$ has a Par\'eto distribution on $[\,1,\infty)$ with 
  parameter $\alpha$,

  \smallskip

  \noindent (ii) $\tau$ has density proportional to $\rho_\gamma^{-\alpha}$,

  \smallskip

  \noindent (iii) $\fS(\lambda)={\ds 1\over\ds\Gamma(1+\gamma)}
  \Bigl( -\lambda^\gamma+\II\{\, \lambda\geq\tau\,\} \gamma
  (\lambda-\tau)^{\gamma-1} \rho_\gamma(\tau) Y\Bigr)$.
}

\bigskip

When $\mu$ is an arbitrary negative number, we will prove by a  rescaling
argument that the limiting triple is
$$
  \Bigl(\fS\bigl((-\mu)^{1/\gamma}\,\cdot\,\bigr) \, ,\, 
        (-\mu)^{-1/\gamma}\tau \, , \,
        (-\mu)^{1-1/\gamma}\rho_\gamma\bigl((-\mu)^{-1/\gamma}\tau\bigr)Y
  \Bigr) \, .
$$

It is not difficult to adapt our proof to the case of the random walk
with negative drift, for which $\gamma=1$, and show that the same result
remains true.

\medskip

One can check, starting from the definition of $\rho_\gamma$ and the fact 
that $Y$ is at least $1$ almost surely, that the maximum of the process
$\fS$ is at least $1$ almost surely.
Note that the trajectories of the limiting process $\fS$ are infinitely
differentiable on the positive half-line, except at $\tau$. If $\gamma<2$,
the trajectories are not differentiable at $\tau+$ but are H\"olderian of
index $\gamma-1$. Thus, on the right of the random time $\tau$ they have
a vertical tangent going upward. If $\gamma> 2$, the trajectories 
are differentiable. If $\gamma=2$, the trajectories are not differentiable
at $\tau$ but admit left and right tangents.

The following pictures show typical paths of the limiting process $\fS$
for different values of $\gamma$.

\bigskip

\setbox1=\hbox to 1.8in{\vbox to 1.5in{
  \hsize=1.8in
  \kern 1.8in
  \includegraphics{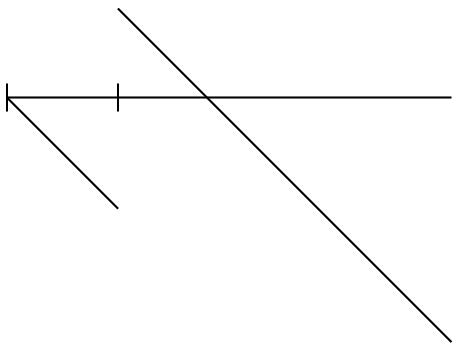}
  \anote{1}{.7}{$\ss\gamma=1$}
  \anote{-.05}{1.5}{$\ss 0$}
  \anote{.4}{1.5}{$\ss\tau$}
  \vfill
  \vss
}\hfill}    
\setbox2=\hbox to 1.8in{\vbox to 1.5in{
  \hsize=1.8in
  \kern 1.8in
  \includegraphics{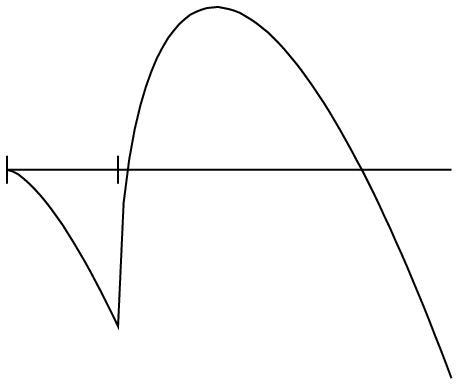}
  \anote{1}{.7}{$\ss\gamma=1.5$}
  \anote{-.05}{1.5}{$\ss 0$}
  \anote{.4}{1.5}{$\ss\tau$}
  \vfill\vss
}\hfill}  

\bigskip

\centerline{\box1\hfill\box2}

\setbox1=\hbox to 1.8in{\vbox to 2in{
  \hsize=1.8in
  \kern 2.2in
  \includegraphics{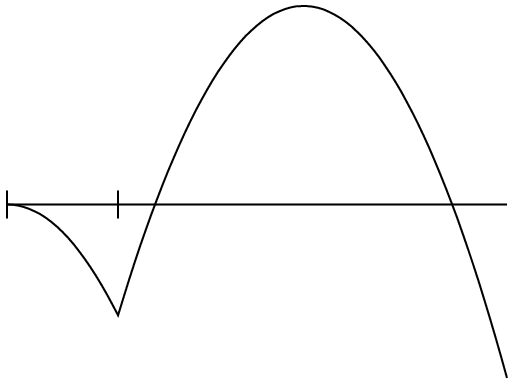}
  \anote{1}{.7}{$\ss\gamma=2$}
  \anote{-.05}{1.5}{$\ss 0$}
  \anote{.4}{1.5}{$\ss\tau$}
  \vfill
  \vss
}\hfill}    
\setbox2=\hbox to 1.9in{\vbox to 2in{
  \hsize=1.8in
  \kern 2.2in
  \includegraphics{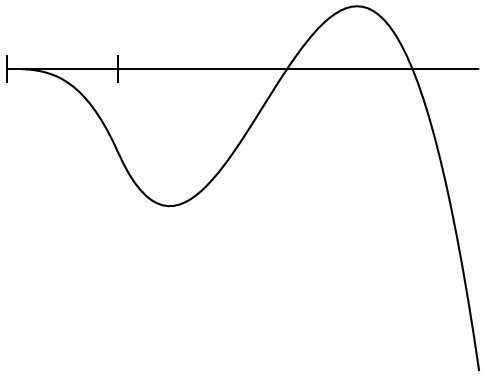}
  \anote{1}{.7}{$\ss\gamma=3$}
  \anote{-.05}{1.5}{$\ss 0$}
  \anote{.4}{1.5}{$\ss\tau$}
  \vfill\vss
}\hfill}  

\bigskip

\centerline{\box1\hfill\box2}

\bigskip

It is noticeable that for $\gamma>2$, the limiting process keeps decreasing
for some time after $\tau$, before increasing to reach its maximum
and finally regains a path asymptotic to that given by the law of large
number. Thus there is a delay not only in reaching the maximum but also
in changing from the negative drift to a positive one which will lead to 
the maximum.

It is not difficult to adapt the proof of Theorem \TrajG\ to study
more precisely what happens near the jump when $\gamma$ is equal to $1$
and the sequence $(g_n)_{n\geq 0}$ tends to infinity. Similarly, one
can see with our proof that when $\gamma$ is less than $1$, the
processes $\fS_t$ do not converge in distribution in the space
$\D[\,0,\infty)$ --- one of the many reasons is that the function $h$ defined
in section \fixedref{3.3.1} is not continuous from $[\,0,1)^2$ to
$\D[\,0,\infty)$. The phenomenological reason is that in order for those 
processes to converge, our proof shows that the limiting function needs to be
$$
  {-1\over\Gamma(1+\gamma)} \bigl( -\lambda^\gamma +\II\{\, \lambda=\tau\,\}
  \gamma\rho_\gamma(\tau)Y\bigr)
$$
which, because of the jump at $\tau$, is not a c\`adl\`ag function. 
It is not particularly difficult, though
somewhat lengthy, to adapt our proof to study the process rescaled
with different time scales, $U(t)$ and $1$, before and after $J_t$, and
show that with these different scalings, both parts converge.

Define the random variables $\calN$ and $\calM$ by
$$
  \calN=\inf\{\, \lambda \,:\, \fS(\lambda)>1\,\} \, ,
$$
and
$$
  \calM
  =\sup_{\lambda>0}\fS(\lambda)
  = {1\over\Gamma(1+\gamma)} \max_{\lambda\geq\tau} -\lambda^\gamma
  +\gamma(\lambda-\tau)^{\gamma-1}\rho_\gamma(\tau)Y
 \, .
$$
By the continuous mapping theorem (Billingsly, 1968, \S 5; Pollard, 1984, 
\S VI.1, example 2), we deduce from Theorem \TrajG\ that
when $\gamma$ is greater than $1$,
the conditional distribution of $\bigl(N_t/U(t),M/t\bigr)$ given $M>t$ 
converge to that of $((-\mu)^{-1/\gamma}\calN,\calM)$.
It follows that the conditional limiting distribution of the
overshoot, $(M-t)/t$, given $M>t$ converges to that of $\calM-1$.
In the same spirit, write $L$ for the first time that the process 
$(S_n)_{n\geq 0}$ attains it maximum. Defining $\calL$ to be the largest 
solution in $[\,\tau,\infty)$ of the equation
$$
  -\calL^{\gamma-1} +(\gamma-1)(\calL-\tau)^{\gamma-2}Y\rho_\gamma(\tau)=0 \, ,
$$
the conditional distribution of 
$L/U(t)$ given $M>t$ converges to the distribution of 
$(-\mu)^{-1/\gamma}\calL$.

The same conclusions hold when $\gamma$ is $1$ and the hypotheses of 
Theorem \ThG\ are satisfied or in the case of a random walk.

\bigskip

\subsection{Examples}
For a FARIMA$(\Phi,d,\Theta)$ process, $g$ is given by \gFARIMA. Let
$c=\Phi(1)/\Theta(1)$. To calculate the associated function $U$, note that
if $u$ tends to infinity, then
$$
  g\Bigl(1-{1\over u}\Bigr)
  =u^d {\Theta\over\Phi}\Bigl( 1-{1\over u}\Bigr)
  \sim {u^d\over c} \, .
$$
We assume that $c$ is positive; if this is not the case, we should replace
$g$ by $-g$ and $X_n$ by $-X_n$ and permute the upper and lower tails in what
follows. Since $c$ is positive, we obtain
$$
  U(t)\sim (ct)^{1/d}
$$
as $t$ tends to infinity. For the FARIMA processes, Akonom and 
Gouri\'eroux (1987) showed directly that
$$
  g_n\sim {n^{d-1}\over c\Gamma(d)} \, ,
$$
as $n$ tends to infinity. Thus, if $d$ is larger than $1$, Theorem \ThG\ yields
$$\eqalign{
  P\{\, M>t\,\}
  &{}\sim (-\mu)^{\alpha(1-1/d)-1/d} (\Id\oF)(c^{1/d}t^{1/d})
   \int_0^\infty\rho_d(v)^{-\alpha} \d v \cr
  &{}\sim (-\mu)^{\alpha (1-1/d)-1/d} c^{(1-\alpha)/d} t^{1/d}\oF(t^{1/d})
   \int_0^\infty \rho_d(v)^{-\alpha}\d v
  }
$$
as $t$ tends to infinity. In contrast, if $d$ is less than $1$ and the tail
balance condition \TailBalance\ holds, then, writing $p$ 
for $\limt \oF(t)/\oF_*(t)$ and $q$ for $\limt F(-t)/\oF_*(t)$, 
Theorem \VgSigned\ yields
$$\displaylines{\quad
  P\{\, M>t\,\} 
  \hfill\cr\hfill
  \sim {1\over d} B\Bigl( {1\over d},\alpha-{1\over d}\Bigr) 
  \Bigl({\Gamma(1+d)\over -\mu}\Bigr)^{1/d} c^{1/d} 
  \bigl(p{g^*}^\alpha+q(-g_*)^{\alpha}\bigr)t^{1/d}\oF_*(t)
  \cr}
$$
as $t$ tends to infinity. There is no explicit expression for $g^*$ and $g_*$
but those can be calculated numerically if needed.

The above asymptotic equivalent sheds further light on the parallel mentioned
in the introduction between FARIMA process and the classical random walk.
Indeed, one may consider the fractional integrated tail of order $1/d$,
$$
  I_{1/d}\oF_*(t)=\int_t^\infty (x-t)^{(1/d)-1}\oF_*(x)\d x \, .
$$
A change of variable $x=(1+\lambda)t$ shows that
$$
  I_{1/d}\oF_*(t)\sim B\Bigl({1\over d},\alpha-{1\over d}\Bigr) t^{1/d}\oF_*(t)
$$
as $t$ tends to infinity. Therefore, when $d$ is less than $1$ 
we have
$$
  P\{\, M>t\,\} 
  \sim \Bigl({\Gamma(1+d)\over -\mu}\Bigr)^{1/d} {c^{1/d}\over d} 
  \bigl(p{g^*}^\alpha+q(-g_*)^{\alpha}\bigr) I_{1/d}\oF_*(t)
$$
as $t$ tends to infinity, continuing the similarity mentioned in the 
introduction between FARIMA processes and the usual random walk, namely
that an integrated tail is replaced with a fractionally integrated tail.

\bigskip

The `random walk' whose increments are themselves a random walk,
corresponds to $g(x)=(1-x)^{-2}$. To apply Theorem \ThG, we need to 
consider $\rho_2(u)=u+\sqrt{2+u^2}$. The change of variable $s=u+\sqrt{2+u^2}$
shows that
$$
  \int_0^\infty \rho_2(u)^{-\alpha}\d u
  = {1\over 2^{\alpha+1}} {3\alpha+1\over\alpha^2-1} \ .
$$
Hence, in this case, we obtain
$$
  P\{\, M>t\,\}
  \sim {(-\mu)^{(\alpha-1)/2}\over 2^{\alpha+1}} {3\alpha+1\over\alpha^2-1}
  \sqrt{t}\, \oF(\sqrt t)
$$
as $t$ tends to infinity.

More generally, a `random walk' whose increments are
a $(g,F)$-process is a $\bigl((1-\Id)^{-1}g,F\bigr)$-process.
Provided the drift is negative, Theorem \ThG\ or \VgSigned\ yield an 
asymptotic estimate on the tail of the distribution of its maximum.

\bigskip


\subsection{Note on the quantiles of the maximum of \poorBold{$(g,F)$}-processes}
Motivated by the last paragraph of the previous section, the purpose of
this section is to describe how the tail of the maximum of a $(g,F)$-process
changes when the increments derive from a process of the same type.
Interestingly, we will see that various quantile functions are asymptotically 
related. Though this is not our main purpose, this connection between 
high-order quantiles is of potential interest in the theory of value
at risk (see e.g.\ Mc Neil, Frey and Embrechts, 2005). Indeed, motivated
by the Basel II regulatory framework, there has been some studies of how
extreme quantiles behave under addition of random variables, with the most
recent research emphasizing cases where the random variables are
dependent (Barbe, Foug\`eres, Genest, 2006; Embrechts, Ne\v slehov\'a, 
W\"utrich, 2008; Embrechts, Lambrigger, W\"utrich, 2008, and references 
therein). The following result gives corresponding results when aggregation
is made according to some moving average scheme and the quantity of interest
is the global maximum of the process.

To state our next results, let $H_+$ be the set of all analytic functions
on $(-1,1)$ whose sequence of Taylor coefficients at $0$ tends to infinity
and are regularly varying of positive index. This set of functions is a 
semi-group both under addition and multiplication. Let $F$ be a fixed 
distribution function with regularly varying tail of index less than $-1$
and with negative mean. 
Let $M_g$ be the maximum of the corresponding $(g,F)$-process. We write
$q_g$ for the function,
$$
  q_g(s)=\inf\bigl\{\, t \, : P\{\, M_g\leq t\,\}\geq 1-1/s\,\bigr\} \, .
$$
This is the quantile function of $M_g$ evaluated at $1-1/s$. We also write
$$
  c_g=(-\mu)^{\alpha (1-1/\gamma)-1/\gamma}
  \int_0^\infty \rho_\gamma^{-\alpha}(u)\d u
$$
for the constant involved in the statement of Theorem \ThG. The quantile
function of $M_g/c_g^{\gamma/(1-\alpha)}$ evaluated at $1-1/s$ 
is 
$$
  \tilde q_g(s)=c_g^{-\gamma/(1-\alpha)} q_g(s) \, . 
$$
We then have the following asymptotic relations showing that the map
$g\in H_+\mapsto \tilde q_g$ is a linear morphism of semigroups in an 
asymptotic sense.

\Proposition{ Let $g$, $g_1$ and $g_2$ be some functions in $H_+$. Then for
  any positive $\lambda$, the following asymptotic equivalences hold at
  infinity.

  \smallskip

  \noindent (i) $\tilde q_{\lambda g}\sim \lambda \tilde q_g$.

  \smallskip

  \noindent (ii) $\tilde q_{g_1+g_2}\sim \tilde q_{g_1}+\tilde q_{g_2}$.

  \smallskip

  \noindent (iii) $\tilde q_{g_1g_2}\sim \tilde q_{g_1}\tilde q_{g_2}$.
}

\bigskip

\Proof Note that with some obvious notations, the relation defining $U_g$, that
is, $g(1-1/U_g)\sim \Id$,
implies $U_g^\leftarrow\sim g(1-1/\Id)$ at infinity. Consequently
$$
  U_{\lambda g}^\leftarrow \sim \lambda U_g^\leftarrow \, ,
  \quad
  U_{g_1+g_2}^\leftarrow\sim U_{g_1}^\leftarrow+U_{g_2}^\leftarrow \, ,
  \quad\hbox{and}\qquad
  U_{g_1g_2}^\leftarrow\sim U_{g_1}^\leftarrow U_{g_2}^\leftarrow \, .
  \eqno{\equa{Quantile}}
$$
Theorem \ThG\ shows that
$$
  c_g (\Id \oF)\circ U_g\circ q_g\sim 1/\Id \, .
$$
Thus,
$$
  q_g\sim U_g^\leftarrow\circ (\Id\oF)^\leftarrow(1/c_g\Id) \, .
$$
In particular, $q_g$ is regularly varying of index $\gamma_g/(\alpha-1)$ and
$$
  q_g\sim c_g^{\gamma_g/(\alpha-1)} U_g^\leftarrow\circ (\Id\oF)^\leftarrow
  (1/\Id) \, .
$$
Hence, $\tilde q_g\sim U_g^\leftarrow\circ (\Id \oF)^\leftarrow (1/\Id)$.
The result then follows from \Quantile.\hfill\qed

\bigskip


\subsection{Concluding remarks}%
The next sections will show the technique used to prove the results of
the current section. As it should be clear at the end of this paper, this
technique can be used to give extensions of Veraverbeke's result in
different directions.

For instance, Veraverbeke's theorem can be interpreted as a statement on 
the probability that a centered random walk crosses a linear moving boundary.
Indeed, let $Z_i$ be the centered random variable $X_i-\mu$. The maximum of
the random walk based on the $X_i$ exceeds $t$ if the random walk based
on the $Z_i$ crosses the boundary $t-\mu i$, which is called a moving
boundary since $t$ translates it upward. In nonlinear renewal theory
(see, e.g., Woodroofe, 1982), calculations of crossing probabilities of 
nonlinear moving boundaries are questions of importance. Our technique 
yields estimates for the crossing of some nonlinear boundaries. However, 
it does not allow one to recover all known results
and in particular our technique does not work for boundaries near the
range of the law of the iterated logarithm --- compare to the remarkably 
general Theorem 1 in Foss, Palmowski and Zachary (2005).

In a similar spirit, but further away from Veraverbeke's theorem, our
technique can be used to evaluate the probability that at least one of the
$X_i$ exceeds $t+b(i)$ where $b$ is an increasing function. Clearly, this
probability can be evaluated directly, and so the purpose of this
remark is only to delineate further the range of usefulness of our technique.
Our technique applies
when $b$ is regularly varying of index greater $1/\alpha$, but not if
$b$ is regularly varying of index $1/\alpha$. Therefore, it does not yield
estimates as sharp as extreme value theory.

On a more positive note, our technique seems useful when some form of 
dependence is present, provided that one has a good representation of the
random variables involved. The $(g,F)$-processes provide a nontrivial example.
Another interesting example is as follows.
Let $(Z_i)_{i\geq 1}$ be a sequence of independent random variables,
equidistributed and centered. Let $p$ be a fixed integer, let $\mu$ be
a negative real number, and let further $X_i=Z_iZ_{i+1}\ldots Z_{i+p}+\mu$. 
Consider the `random walk' $(S_n)_{n\geq 0}$
associated to the $X_i$ and defined by $S_0=0$ and $S_i=S_{i-1}+X_i$ 
if $i\geq 1$. The increments of this random walk are the $p$-dependent
sequence $(X_i)_{i\geq 1}$. Such process has been considered in the
context of large deviations by Choi, Cover and Csiz\'ar (1987) as well as
Bolthausen (1993). Our technique allows one to show that if
$X_i$ has a distribution function $F$ whose tail is regularly varying
of index $\alpha$, then Veraverbeke's result remains valid, namely that
$$
  P\{\, \max_{n\geq 0} S_n>t\,\}\sim {1\over-\mu} {t\oF(t)\over \alpha-1} 
  \eqno{\equa{conclusionA}}
$$
as $t$ tends to infinity. To sketch the proof of this assertion requires
the notation to be developed in the next section, and, perhaps, the next
paragraph can only be understood after reading the remainder of this
paper; however, the shortness of the following sketch seems a compelling
argument in favor of the general framework which we will develop in the
next section.

To prove \conclusionA, we evaluate
$$
  \psi_{i,n}(x)=E(S_i\,|\,X_n=x) 
  = \cases{i\mu & if $i<n$\cr
           (i-1)\mu +x & if $i\geq n$\cr}
$$
Thus, $\psi_n(x)=\max_{i\geq 1} \psi_{i,n}(x)$ is invertible on the preimage
of some interval $(t_0,\infty)$ and $\psi_n^{-1}(t)=t-(n-1)\mu$ on this
preimage. It follows that we can take $U$ and $\chi$ to be the identity
function, and $\rho(x)=1-\mu x$. From this, the result can be guessed
using formula \fixedref{(3.1.2)}. To actually prove the result, steps
3--7 of subsection \fixedref{3.1} are established as follows. We first split
$S_n$ as a sum of the $p$ random walks with independent increments
$$
  S_{n,j}
  =\sum_{\ss 1\leq i\leq n\hfill\atop\ss i\equiv j\, {\rm mod } p\hfill} X_i 
  \, , \qquad 0\leq j<p \, .
$$ 
Then, whenever we need to estimate a probability involving the event $S_n>t$
we note that this event is included in $\cup_{0\leq j<p} \{\, S_{n,j}>t/p\,\}$
and use Bonferroni's inequality. The result then follows by the estimates
of section \fixedref{6.5}.

\bigskip


\def\preveq{(\the\sectionnumber .\the\subsectionnumber .\the\equanumber)}
\def\prevs{\the\sectionnumber .\the\subsectionnumber .\the\snumber }
\section{Veraverbeke's theorem at large}%
The proofs of our theorems are conceptually simple, but this simplicity is
somewhat lost in the many steps needed in its execution. In order to 
make this simplicity more 
obvious and intuitive, as well as in order to make our technique easy to adapt
to different problems, we first describe a general scheme for how to prove
the type of results which we are aiming for. In the second
subsection we prove a theorem which asserts that, indeed, whenever 
this general scheme can be applied, it yields the correct result. It
will be used to prove the results of section 2.

In this section only, we consider a stochastic process $(S_n)_{n\geq 1}$
built through functions $S_n$ mapping a sequence of independent random 
variables $(X_n)_{n\geq 1}$ into the real line. We are seeking some
tail estimate for the maximum of $(S_n)_{n\geq 1}$. Clearly, nothing
useful can be said with that level of generality, but our purpose is to
describe a technique at a conceptual level. We make no claim that this
technique yields the correct result in general --- as a matter of fact,
it is very easy to find counterexamples --- but the remainder of this
paper will show that this description can be most useful.

We will use the following notation.

\Definition{ A function $f$ mapping a neighborhood of infinity to a 
neighborhood of infinity has an asymptotic inverse if it is asymptotically
equivalent to a monotone function and there exists a function $f^\leftarrow$
such that $f\circ f^\leftarrow\sim f^\leftarrow\circ f \sim \Id$ at
infinity.
}

\bigskip

\subsection{The single large jump heuristic}%
Our purpose in this subsection is not to do rigorous mathematics but to give
some useful intuitions. The basic idea underlying the single large jump 
heuristic is that a single large $X_n$ is what is likely to make the maximum
$M$ of the process $(S_n)_{n\geq 1} $ to be large, and that the other $X_i$ 
contribute to the
process in an average way, having in mind some form of law of large numbers.
This leads us to consider the conditional expectations functions
$$
  \psi_{i,n}(x)=E(S_i|X_n=x) \, , \qquad i,n\geq 1 \, .
$$
If $X_n$ is large, we expect the process $(S_k)_{k\geq 0}$ to reach the 
level $\psi_{i,n}(X_n)$ at time $i$. Therefore, defining
$$
  \psi_n(x)=\max_{i\geq 1} \psi_{i,n}(x) \, ,
$$
the maximum of the process, given $X_n$, is expected to reach the level
$\psi_n(X_n)$. It exceeds $t$ if $\psi_n(X_n)$ does, that is, if $\psi_n$
is increasing and invertible, if $X_n>\psi_n^{-1}(t)$. So, we anticipate 
that, as $t$ tends to infinity,
$$
  P\{\, M>t\,\} 
  \sim P\{\, \exists n\geq 1 \, :\, X_n>\psi_n^{-1}(t)\,\}
  \sim\sum_{n\geq 1}\oF\circ \psi_n^{-1}(t) \, .
$$
In what follows, let $r$ be a function such that 
$$
  r(t)\sim \sum_{n\geq 1}\oF\circ \psi_n^{-1}(t) \, 
$$
at infinity. It is of course assumed that this function tends to $0$ at
infinity, that is $M$ is almost surely finite.
Thus, $r$ is a tentative asymptotic equivalent for the probability
that $M$ exceeds $t$. A general scheme to turn this tentative equivalent into
an actual one is as follows. It has two parts, an analytical one and a 
probabilistic one.

\medskip

\noindent {\it Analytical part.} The purpose of this part is to obtain useful
information on the function $\psi_n^{-1}$.

\medskip

\noindent{\bf Step 1.} Find two regularly varying functions $U$ and $\chi$ 
with limit
infinity at infinity and a function $\rho$ continuous on the positive 
half-line, bounded away from $0$ on any
compact subset of the positive half-line, such that $\rho^{-\alpha}$ is
Lebesgue integrable on the nonnegative half-line, and such that the asymptotic 
factorization
$$
  \psi_{\lfloor xU(t)\rfloor}^{-1} (t)
  \sim \chi (t)\rho (x)
  \eqno{{\rm (S1)}}
$$
holds as $t$ tends to infinity, uniformly in $x$ in any compact subset of the 
positive half-line. Implicit in this assertion is 
that for $t$ large enough, $\psi_{\lfloor xU(t)\rfloor}$ is invertible on 
the preimage of $[\, t,\infty)$.

The functions $\chi$ and $\rho$ are not unique. However, (S1) implies
that for another such pair, say $(\chi_1,\rho_1)$, we have, 
for any positive $x$,
$$
  {\chi\over\chi_1}(t)\sim {\rho_1\over\rho}(x)
$$
as $t$ tends to infinity. This forces $\rho/\rho_1$ to be constant, equal
to some $c$ say, and then $\chi_1\sim c \chi$. The constant $c$ is positive
for both $\chi$ and $\chi_1$ are assumed to tend to infinity at infinity. 
It follows that even though the functions
$\chi$ and $\rho$ are not unique, they are asymptotically unique up to
a positive multiplicative constant. Our results do not depend on the
choice of the constant, and in applications we will choose whatever
constant makes the calculation less cumbersome.

\medskip

In what follows, to any positive real number $\epsilon$ less than $1$ we 
associate the set
$$
  I_{\epsilon,t}
  =\{\, n\in\NN \, : \, \epsilon U(t)\leq n\leq U(t)/\epsilon \,\} \, ,
$$
and we write
$$
  N_{\epsilon,t}=\lfloor U(t)/\epsilon\rfloor
$$
for its largest element. Often we will drop the subscripts, writing $I$ and $N$
for $I_{\epsilon,t}$ and $N_{\epsilon,t}$.

\medskip

\noindent{\bf Step 2.} Prove that the asymptotic behavior of $\sum_{n\geq 1}
\oF\circ\psi_n^{-1}(t)$ as $t$ tends to infinity is driven by the terms
for which $n$ is of order $U(t)$, that is,
$$\displaylines{\qquad
    1=\lim_{\epsilon\to 0}\liminft {1\over r(t)} \sum_{n\in I}
    \oF\circ \psi_n^{-1}(t) 
  \hfill\cr\hfill
    {}\leq \lim_{\epsilon\to 0}\limsupt {1\over r(t)}
    \sum_{n\in I} \oF\circ \psi_n^{-1}(t)
    = 1\, . 
  \qquad {\rm (S2)}\cr }
$$
In particular, this step suggests that in relation to step 1, we should also
obtain order of
magnitudes or rather crude bounds for $\psi_n^{-1}$ when $n$ is outside $I$.

Note that the completion of step 2 implies that the function $r$ can be
identified as follows. Since $\oF$ is regularly varying, steps 1 and 2 imply
$$\eqalignno{ 
  \sum_{n\in I} \oF\circ \psi_n^{-1}(t)
  &{}\sim \sum_{n\in I} \oF\Bigl( \chi (t)\rho\Bigl({n\over U(t)}\Bigr)\Bigr)
   \cr
  &{}\sim \int_{\epsilon U(t)}^{U(t)/\epsilon} \oF\circ \chi (t)
   \rho\Bigl({u\over U(t)}\Bigr)^{-\alpha} \d u \cr
  &{}\sim (U\oF\circ\chi) (t) \int_\epsilon^{1/\epsilon} \rho(u)^{-\alpha}
   \d u \, . 
  &\equa{rcalculation}\cr
  }
$$
Thus, completion of steps 1 and 2 implies that
$$
  r(t)= (U\oF\circ \chi) (t)\int_0^\infty \rho(v)^{-\alpha}\d v
  \, .
  \eqno{\equa{rDef}}
$$
In particular $r$ is regularly varying. Given how $r$ was initially defined, 
this suggests that one could easily guess the 
tail behavior of $M$ by a simple examination of $\psi_n^{-1}$, just guessing
what $U$, $\chi$ and $\rho$ are. This has been illustrated in the last example
discussed in section 2.5.

\medskip

\noindent{\it Probabilistic part.} This part consists in proving that if the
process reaches the level $t$, then it is unlikely to occur at a time too
small or too large, and that some form of law of large numbers holds.

\medskip

\noindent{\bf Step 3.} Prove that the process is unlikely to reach the
level $t$ at a time of smaller order than $U(t)$, that is
$$
  \limeps\limsupt {1\over r(t)} 
  P\{\, \exists n\leq \epsilon U(t) \, :\, S_n>t\,\} = 0 \, .
  \eqno{{\rm (S3)}}
$$
The proof of such result is sometimes made easier by the following remark.
For any fixed $\theta$, consider the events
$$
  B_i=\{\, X_i\leq \theta\chi (t)\,\} \, ,
$$
whose notation does not keep track of the dependence on $\theta$ and $t$.
We write $B_i^\complement$ for the complement of $B_i$.
Since
$$
  P\bigl\{\, \bcup_{i\leq \epsilon U(t)} B_i^\complement\,\bigr\} 
  \sim \epsilon U(t)\oF\bigl( \theta\chi(t)\bigr)
  \sim \epsilon\theta^{-\alpha}(U\oF\circ\chi)(t) \, ,
$$
it suffices to prove that
$$
  \limeps\limsupt {1\over r(t)} 
  P\bigl\{\,  \exists n\leq\epsilon U(t)\, :\, S_n>t \, ;\, 
              \bcap_{i\leq\epsilon U(t)} B_i\,\bigr\} = 0 \, .
  \eqno{\equa{SThreeAlt}}
$$
The advantage of this formulation is that on the event 
$\bcap_{i\leq\epsilon U(t)} B_i$ the random variables are bounded, and
many more inequalities exist for bounded random variables than for unbounded
ones.

\medskip

\noindent{\bf Step 4.} Prove that the process is unlikely to reach the level
$t$ at a time of larger order than $U(t)$, that is,
$$
  \limeps\limsupt {1\over r(t)} 
  P\{\, \exists n\geq U(t)/\epsilon\, :\, S_n>t\,\} = 0 \, .
  \eqno{{\rm (S4)}}
$$
As in step 3, it is sometimes useful to replace the unbounded variables $X_i$
by bounded ones. Let $\kappa$ be the index of regular variation of $\chi$ and
assume that $\alpha\kappa\gamma$ is greater than $1$. We have
$$\displaylines{\qquad
    P\bigl\{\, \exists n\geq 0 \, :\, X_n>\theta 
    \bigl( \chi\circ U^\leftarrow(n)+\chi(t)/\eta\bigr)\,\bigr\}
  \hfill\cr\noalign{\vskip 3pt}\hfill
    {}\sim\sum_{n\geq 1} \oF\Bigl(\theta \bigl( \chi\circ U^\leftarrow(n)
                +\chi(t)/\eta\bigr)\Bigr) \, .
  \qquad\cr}
$$
Replacing this series by the corresponding Riemann integral,
$$
  \int_1^\infty \oF\Bigl( \theta\bigl(\chi\circ U^\leftarrow(u)+\chi(t)
  /\eta\bigr)\Bigr) \d u \, ,
$$
making the change of variable $u=\lambda U(t)$ and using the regular variation
of $\oF$, the series is 
asymptotically equivalent to
$$
  \theta^{-\alpha} (U\oF\circ\chi)(t)
  \int_0^\infty (\lambda^{\kappa\gamma}+1/\eta)^{-\alpha} \d\lambda
$$
as $t$ tends to infinity. Since $(\lambda^{\kappa\gamma}+1/\eta)^{-\alpha}$
is at most $\eta^\alpha$ when $\lambda$ is in $[\,0,1\,]$ and at 
most $\lambda^{-\alpha\kappa\gamma}$ when $\lambda$ is at least $1$, 
by dominated convergence,
$$
  \lim_{\eta\to 0}\limsupt {1\over r(t)} P\bigl\{\, \exists n \geq 0\, :\,
    X_n>\theta \bigl( \chi\circ U^\leftarrow(n)+\chi (t)/\eta\bigr)\,\bigr\}
  = 0 \, .
$$
Therefore, we can replace the original problem of this step by that 
of proving 
$$\displaylines{\quad
  \lim_{\eta\to 0}\limsupt {1\over r(t)} 
  P\Bigl\{\, \exists n\geq U(t)/\eta \, :\, S_n>t \, ; \,
  \hfill\cr\hfill
  \bcap_{i\geq 1} \bigl\{\, X_i\leq \theta \bigl( 
  \chi\circ U^\leftarrow(i)+\chi(t)/\eta\bigr)\,\bigr\} \,\Bigr\} =0 \, .
  \quad\equa{StepFourAlt}
  \cr}
$$

\medskip

The next step consists in formalizing the single large jump heuristic, proving
that for $M$ to exceed $t$ then at least one large jump had likely occurred.

\medskip

\noindent{\bf Step 5.} Prove that for the process to exceed $t$ at a time
in $I$, we need at least one variable prior to that time to 
exceed $\theta \chi (t)$, that is, for any positive $\epsilon$, there 
exists some positive $\theta$ such that
$$
  P\Bigl(\, \bcup_{i\in I} \Bigl(\{\, S_i>t \,\}\cap \bcap_{n\leq i} B_n \Bigr)
  \Bigr)
  = o\bigl( r(t)\bigr) \, .
  \eqno{{\rm (S5)}}
$$

Note that if this holds for some $\theta$ then it holds for any smaller one,
because the sets $B_i$ are decreasing in $\theta$.

Interestingly, completion of steps 1, 2 and 5 are enough to show that it is
unlikely that $M$ reaches $t$ because two of the $X_i$'s are large. Indeed,
we have
$$
  P\Bigl(\, \bcup_{\ss i,j\leq N\hfill\atop\ss i\not= j\hfill} 
  B_i^\complement\cap B_j^\complement\,\Bigr)
  \sim {N\choose 2} \oF\bigl( \theta \chi(t)\bigr)^2
  \sim {1\over 2\theta^{2\alpha} \epsilon^2} (U\oF\circ \chi)^2 (t) \, ,
$$
and since $r\asymp U\oF\circ \chi$ tends to $0$, this implies
$$
  \limt {1\over r(t)} 
  P\Bigl\{\, \bcup_{\ss i,j\leq N\hfill\atop\ss i\not= j\hfill} 
   B_i^\complement\cap B_j^\complement \,\Bigr\}
  = 0 \, .
$$

\medskip

We can then move on to the next step, showing that when a single $X_n$
is large, then the process can be approximated by its conditional expectation
given that large random variable. Recall that $N=\lfloor U(t)/\epsilon\rfloor$.
It is convenient for what follows to 
introduce the events that all $X_i$, $1\leq i\leq N$, are at most $\theta
\chi(t)$ except perhaps $X_n$, that is,
$$
  C_n=\bcap_{\ss 1\leq i \leq N\hfill\atop\ss i\not=n\hfill} B_i \, ,
$$
and the events that the $S_i$, $1\leq i\leq n$, are well approximated 
by $\psi_{i,n}(X_n)$, that is,
$$
  D_n=\bcap_{1\leq i\leq N} \{\, |S_i-\psi_{i,n}(X_n)|\leq \delta t\,\}
  \, .
$$

\medskip

\noindent{\bf Step 6.} Prove that if $X_n$ exceeds $\theta \chi (t)$ and all
the other $X_i$, $i\leq N$, are at most $\theta\chi (t)$, then each $S_i$
is about $\psi_{i,n}(X_n)$; more precisely, prove that for any positive 
$\delta$,
$$
  P\bigl(\, \bcup_{n\in I}
  B_n^\complement\cap C_n \cap D_n^\complement\,\bigr)
  = o\bigl( r(t)\bigr) \, .
  \eqno{{\rm (S6)}}
$$

This is in fact a little stronger than what we need, and sometimes a one-sided 
bound, replacing $|S_i-\psi_{i,n}(X_n)|$ by $S_i-\psi_{i,n}(X_n)$, may
suffice.

A naive and yet effective way to prove such law of large numbers is to show
first that 
$$\displaylines{\qquad
  \max_{n\in I} P( B_n^\complement\cap C_n \cap D_n^\complement )
  {}= o\bigl( r(t)/U(t)\bigr) \, ,
  \qquad\cr}
$$
and then use Bonferroni's inequality, upon noting that the cardinality of
$I$ is of order $U(t)$.

\medskip

The combination of all these steps suggests that $M>t$ occurs most likely 
because a single $X_n$ is large, that conditionally on this event, $S_i$ is
about $\psi_{i,n}(X_n)$, and that the maximum of the process will indeed
exceed $t$ if some $\psi_{i,n}(X_n)$ exceeds $t$, that is if $\psi_n(X_n)$
does. In fact, we will prove rigorously in the next subsection that 
completion of steps 1--6, in other words (S1)--(S6), implies the upper bound
$$
  \limsup_{t\to\infty} {1\over r(t)} P\{\, M>t\,\} \leq 1 \, .
$$

To obtain a matching lower bound, additional knowledge seems needed for the
following reason. Let $\delta$ be a positive real number. If $\psi_n(X_n)$
exceeds $(1+2\delta)t$, there exists an integer $i$ such that 
$$
  \psi_{i,n}(X_n)>(1+\delta)t \, .
  \eqno{\equa{HeurA}}
$$
We would like to use step 6 to prove that, perhaps up to intersecting with
a further set,
$$
  S_i\geq\psi_{i,n}(X_n)-\delta t > t
$$
and so $M$ exceeds $t$, suggesting 
that $\sum_{n\in I}\oF\circ\psi_n^{-1}(t)$
is an asymptotic lower bound for the probability that the process reaches the
level $t$ at some time. The problem with this approach is that nothing
guarantees that the $i$ involved in \HeurA\ stays of order $U(t)$, and,
therefore, that step 6 gives the needed law of large numbers on the 
proper range
of $i$. Various assumptions could be made to remove this difficulty. In
some cases we may adapt an argument due to Zachary (2004) while in others
the following may do.

\medskip

\noindent{\bf Step 7.} Let $i(n,x)$ be an integer which maximizes 
$\psi_{i,n}(x)$. Prove that for any positive $\delta$ and $\theta$ less 
than $1$, there exists a positive $\eta$ such that for any $t$ large enough,
$$
  \{\, i(n,x) \, :\, \theta\chi (t) \leq x\leq \chi(t)/\theta \, , \, 
       n\in I_{\delta,t}\,\}
  \subset I_{\eta,t} \, .
  \eqno{{\rm (S7)}}
$$
We will show in the next subsection that if in addition to the previous steps
this last one can be completed then
$$
  \liminf_{t\to\infty} {1\over r(t)} P\{\, M>t\,\} \geq 1\, .
$$
Therefore, once all seven steps have been verified we obtain
$$
  P\{\, M>t\,\}\sim U(t)\oF\circ \chi(t) \int_0^\infty \rho(u)^{-\alpha}\d u
$$
as $t$ tends to infinity.

\bigskip

\noindent{\bf Remarks.} While we defined $\psi_{i,n}$ as a conditional 
expectation, and from there $\psi_n$, we could as well have defined those two
functions in a more axiomatic way, with no connection to conditional 
expectation, as functions which allows us to
carry out steps 1--6 if we are seeking only an upper bound, or 1--7 if we
are seeking an asymptotic equivalent for the tail probability of $M$.

In step 1, we assumed that $\psi_n$ is invertible 
on the preimage of some interval $(t_0,\infty)$. This assumption could be
replaced by a weaker one using asymptotic inverse; however, that requires
some form of uniformity with respect to $n$ in the asymptotic inversion.
While technically possible, such refinement does not seem relevant in
applications.

\bigskip

\subsection{From the heuristic to a theorem}%
The previous subsection sketched a possible path to obtain an asymptotic 
equivalent of the probability that the process reaches a large level at some 
time. In this subsection, we prove
rigorously that this scheme, if it can be completed, indeed yields an 
asymptotic equivalent of the probability that $M$ exceeds $t$. Given
its unsightly assumptions, one may be skeptical that the following
theorem is of any value, but the next sections will demonstrate that
its virtue is to break somewhat complicated problems into bits far
more tractable. In particular, this theorem will be used to prove the
results given in section \fixedref{2}.

\Theorem{\label{UsefulStatement}%
  Referring to the previous subsection, if steps 1--7 have been completed,
  that is, if (S1)--(S7) hold, then
  $P\{\, M>t\,\}\sim r(t)$ as $t$ tends to infinity.
}

\bigskip

Needless to say that the function $r$ in this statement refers to that 
defined in \rDef.

\Remark The proof of Theorem \UsefulStatement\ shows that under (S1)--(S5)
and the one-sided version of (S6) with $D_n$ replaced by
$$
  D_n=\bcap_{i\leq N}\{\, S_i-\psi_{i,n}(X_n)\leq \delta t\,\}
$$
then
$$
  \limsup_{t\to\infty} P\{\, M>t\,\}/r(t) \leq 1 \, .
$$
If in addition the two sided version of (S6) holds as well as (S7) then
$$
  \liminft P\{\, M>t\,\}/r(t)\geq 1\, . 
$$

\bigskip

\Proof We first derive an upper bound  for the probability that $M$ exceeds $t$
under (S1)--(S6), and, with the addition of (S7), a matching lower bound. 

\noindent{\it Upper bound.} We set
$$
  A_i=\{\, S_i>t\,\} \, .
$$
Let $\delta$ be a positive real number. 
Steps 3 and 4 show that we can find a positive $\epsilon$ such that, ultimately
in $t$,
$$
  P\{\, M>t\,\}
  \leq P\bigl( \bcup_{i\in I} A_i\bigr) + \delta r(t) \, .
$$
Using step 5, find $\theta$ such that, for any $t$ large enough,
$$
  P\Bigl( \bcup_{i\in I} \Bigl( A_i \cap \bcap_{n\leq i} B_i \Bigr) 
  \leq \delta r(t) \, .
$$
Then, the equality
$$
  \bcup_{i\in I} A_i
  = \Bigl( \bcup_{i\in I} \Bigl(A_i\cap\bcup_{n\leq i} B_n^\complement\Bigr)
    \Bigr)
  \cup \Bigl(\bcup_{i\in I} \Bigl( A_i\cap \bcap_{n\leq i} B_n\Bigr)\Bigr) 
$$
yields, for any $t$ large enough,
$$
  P\{\, M>t\,\} \leq 
  P\Bigl( \bcup_{i\in I} A_i\cap \bcup_{n\leq N}B_n^\complement\Bigr)
  + 2\delta r(t) \, .
  \eqno{\equa{HToThA}}
$$
Consider the event that all random variables before $N$, except perhaps $X_n$,
are at most $\theta \chi (t)$, that is
$$
  C_n=\bcap_{\ss 1\leq i\leq N\hfill\atop\ss i\not= n\hfill} B_i \, .
$$
The event $B_n^\complement\cap C_n$ expresses that all the $X_i$, $i\leq N$,
but $X_n$ are at most $\theta\chi (t)$. We have the identity
$$\displaylines{\qquad
  A_i\cap \bcup_{n\leq N} B_n^\complement
  \hfill\cr\hfill
  \eqalign{
  {}={}& A_i \cap \bcup_{n\leq N} \bigl( (B_n^\complement\cap C_n)\cup
         (B_n^\complement\cap C_n^\complement)\bigr) \cr
  {}={}& \bcup_{n\leq N} (A_i\cap B_n^\complement\cap C_n)\cup\bcup_{n\leq N}
         (A_i\cap B_n^\complement\cap C_n^\complement) \, . \cr
  }\qquad\cr}
$$
But
$$
  B_n^\complement\cap C_n^\complement
  = B_n^\complement\cap \bcup_{\ss i\leq N\hfill\atop\ss i\not=n\hfill} 
    B_i^\complement
  = \bcup_{\ss i\leq N\hfill \atop \ss i\not= n\hfill} 
  B_n^\complement\cap B_i^\complement
 \, ,
$$
and we saw after step 5 that 
$\bcup_{\ss n,i\leq N\hfill\atop\ss i\not= n\hfill}B_n^\complement\cap 
B_i^\complement$ has
probability $o\bigl( r(t)\bigr)$. Thus, for $t$ large enough, \HToThA\ shows
that
$$
  P\{\, M>t\,\}
  \leq P\Bigl( \bcup_{i\in I}\bcup_{n\leq N} 
  A_i\cap B_n^\complement\cap C_n\Bigr) + 3\delta r(t) \, .
  \eqno{\equa{HToThB}}
$$
We then consider the events
$$
  D_n=\bcap_{i\leq N} \{\, S_i-\psi_{i,n}(X_n)\leq \delta t\,\} \, .
$$
Completion of step 6 --- in fact, the one-sided version would suffice here ---
ensures that
$$
  P\Bigl(  \bcup_{n\in I} B_n^\complement\cap C_n\cap D_n^\complement \Bigr)
  = o\bigl( r(t)\bigr)
  \eqno{\equa{HToThC}}
$$
as $t$ tends to infinity. Thus, \HToThB\ implies that for any $t$ large enough,
$$
  P\{\, M>t\,\}
  \leq P\Bigl( \bcup_{i\in I}\bcup_{n\leq N} A_i\cap B_n^\complement\cap C_n
  \cap D_n\Bigr) + 4\delta r(t) \, .
$$
On $A_i\cap D_n$,
$$
  t<S_i < \delta t +\psi_{i,n}(X_n) \, .
$$
Thus, we proved that, ultimately in $t$,
$$\eqalign{
  P\{\, M>t\,\}
  &{}\leq P\Bigl(\bcup_{i\in I}\bcup_{n\leq N} 
    \bigl\{\, \psi_{i,n}(X_n)> (1-\delta) t\,\bigr\}\Bigr) + 4\delta r(t) \cr
  &{}\leq \sum_{n\geq 1} 
    P\bigl\{\, \max_{i\geq 1} \psi_{i,n}(X_n) > (1-\delta)t\,\bigr\} 
    + 4\delta r(t) \cr
  &{}\leq \sum_{n\geq 1} \oF\circ \psi_n^{-1} \bigl( (1-\delta)t\bigr)
    + 4\delta r(t) \cr
  &{}\leq r\bigl( (1-\delta)t\bigr) + 5\delta r(t) \, . \cr
  }
$$
Since $\delta$ is arbitrary and $r$ is regularly varying, it follows that
$$
  \limsupt {1\over r(t)} P\{\, M>t\,\} \leq 1 \, .
$$

\noindent{\it Lower bound.} Let $\epsilon$ be a positive real number.
Using step 2, let $\delta$ be such that
$$
  (1-\epsilon)r(t)
  \leq \sum_{n\in I_{\delta,t}}\oF\circ\psi_n^{-1}(t)
  \leq (1+\epsilon) r(t)
$$
ultimately. Let $\theta$ be small enough so that
$\theta^\alpha/\delta \leq\epsilon$. Consider the events
$$
  F_n =\{\, \psi_n(X_n)>(1+\delta)t \, ;\, 
            \theta \chi(t)<X_n\leq \chi(t)/\theta \,\} \, ,
$$
and let $B_n$, $C_n$ be the same events as defined previously, with $C_n$ 
defined with reference to the set $I_{\eta,t}$ obtained from step 7. 
Let $N=\lfloor U(t)/\eta\rfloor$ be the largest element of that $I_{\eta,t}$ 
and redefine $D_n$ to be
$$
  D_n=\bcap_{i\leq N} \{\, S_i-\psi_{i,n}(X_n)\geq -\delta t\,\} \, .
$$
Recall the notation $i(n,x)$ introduced in step 7. For $n$ in $I_{\delta,t}$
and on $F_n\cap D_n$,
$$
  \psi_{i(n,X_n),n}(X_n) > (1+\delta)t
$$
and
$$
  S_{i(n,X_n)} 
  \geq \psi_{i(n,X_n),n} (X_n)-\delta t
  > t \, .
$$
Thus, 
$$\eqalignno{
  P\{\, M>t\,\}
  &{}\geq P\bigl\{\, \bcup_{n\in I_{\delta,t}} 
   F_n\cap B_n^\complement\cap C_n\cap D_n\,\bigr\} \cr
  &{}\geq P\bigl\{\, \bcup_{n\in I_{\delta,t}} F_n\,\bigr\}
   - P\bigl\{\, \bcup_{n\in I_{\delta,t}} F_n\cap 
   (B_n^\complement\cap C_n\cap D_n)^\complement\,\bigr\} \, .\cr
  &
  &\equa{HToThD}\cr
  }
$$
We consider the event
$$
  (B_n^\complement\cap C_n\cap D_n)^\complement
  = B_n
    \cup (B_n^\complement\cap C_n^\complement) 
    \cup (B_n^\complement\cap C_n\cap D_n^\complement) \, .
$$
Note that $F_n\cap B_n=\emptyset$ for $\theta\chi(t)<X_n$ in $F_n$ while
$X_n\leq\theta\chi(t)$ on $B_n$. Recall that the event $C_n$ is defined
with reference to the set $I_{\eta,t}$ obtained from step 7, while
\HToThD\ involves the different set $I_{\delta,t}$. Taking $\eta$ to be 
at most $\delta$, which can be done without any loss of generality,
guarantees
$$
  \bcup_{n\in I_{\delta,t}}(B_n^\complement\cap C_n^\complement)
  \subset \bcup_{n\in I_{\eta,t}} (B_n^\complement\cap C_n^\complement)
$$
and, as mentioned after step 5, the event in the right hand side of this
inclusion had probability $o\bigl(r(t)\bigr)$ as $t$ tends to infinity. 
Finally, by step 6, one has
$$
  P\Bigl\{\, \bcup_{n\in I_{\eta,t}}
  (B_n^\complement\cap C_n\cap D_n^\complement)\,\Bigr\}
  = o\bigl( r(t)\bigr)
$$
as $t$ tends to infinity. Therefore, \HToThD\ yields
$$
  P\{\, M>t\,\}
  \geq P\bigl\{\, \bcup_{n\in I_{\delta,t}} F_n\,\bigr\} + o\bigl(r(t)\bigr)
  \, .
  \eqno{\equa{HToThE}}
$$
Given step 1, we can also take $\theta$ small enough so that 
$\psi_n(X_n)>(1+\delta)t$ and $n\in I_{\delta,t}$ guarantees 
$X_n>\theta\chi(t)$. Then \HToThE\ implies
$$
  P\{\, M>t\,\}
  \geq \sum_{n\in I_{\delta,t}} \oF\circ\psi_n^{-1}\bigl((1+\delta)t\bigr)
   - \sharp I_{\delta,t}\oF\bigl( \chi(t)/\theta\bigr) + o\bigl( r(t)\bigr)
  \, .
$$
Thus, given our choice of $\delta$, we obtain for any $t$ large enough
$$\displaylines{\qquad
  P\{\, M>t\,\}
  \geq (1-\epsilon) r\bigl( (1+\delta)t\bigr) 
  \hfill\cr\hfill
  {}- \theta^\alpha 
  \Bigl( {1\over\delta}-\delta\Bigr) U(t)\oF\circ\chi(t) \bigl(1+o(1)\bigr)
  + o\bigl( r(t)\bigr) \, .
  \qquad\cr}
$$
Since $\theta^\alpha/\delta$ is at most $\epsilon$ and $\epsilon$ is arbitrary,
regular variation of $r$ yields
$$
  \liminft {1\over r(t)} P\{\, M>t\,\} \geq 1 \, .
  \eqno{\qed}
$$

\subsection{Analysis of the paths leading to a large maximum}%
The purpose of this subsection is to examine the likely trajectories of the
process which lead to a large maximum, in the same formal framework as in 
the previous subsection. That is, we are seeking for the limiting distribution
of the process $(S_n)_{n\geq 0}$ conditionally on $M$ exceeding $t$, as $t$
tends to infinity. Clearly the process needs to be rescaled
to avoid degeneracy. The right rescaling is suggested by the proof of
Theorem \UsefulStatement\ and that proof also suggests introducing
other random variables of interests. We define $N_t$ to be the first time
that $S_n$ exceeds $t$ and $J_t$ the index of the largest random variable
$X_i$ among $X_1,\ldots, X_{N_t}$, that is
$$
  X_{J_t}=\max_{1\leq i\leq N_t} X_i \, ,
$$
with the convention that $J_t$ is minimal in case of ties. The proof of Theorem
\UsefulStatement\ suggests that $J_t$ is of order $U(t)$ while $X_{J_t}$ is
of order $\chi(t)$. Thus, it is natural to introduce the random variables
$$
  \tau_t=J_t/U(t)\qquad\hbox{ and }\qquad Y_t=X_{J_t}/\chi (t) \, ,
$$
as well as the rescaled process
$$
  \fS_t(\lambda)=S_{\lfloor\lambda U(t)\rfloor}/t \, .
$$
This process belongs to the space $\D[\,0,\infty)$ of all real-valued 
c\`adl\`ag functions endowed with the projective topology inherited 
from the Skorokhod topology on $\D[\,0,1/\epsilon\,]$ for any positive
$\epsilon$ (see e.g.\ Billingsley, 1968, chapter 3; Pollard, 1984, chapter 6). 
Step 6 of subsection
\fixedref{3.1} also suggests that $S_{\lfloor \lambda U(t)\rfloor}$ should
be about $\psi_{\lfloor\lambda U(t)\rfloor,\lfloor\tau_t U(t)\rfloor}
\bigl(\chi(t)Y_t\bigr)$. Therefore, for the process to converge it is natural
to assume that there is a function $h$ on $[\,0,\infty)^3$ such that for any
$\lambda$, $\tau$ and $y$,
$$
  \limt t^{-1}\psi_{\lfloor\lambda U(t)\rfloor,\lfloor\tau U(t)\rfloor}
\bigl(\chi(t)y\bigr)=h(\lambda,\tau,y) \, .
$$
This pointwise convergence is not sufficient to guarantee the convergence
in distribution of the process $\fS_t$ in $\D[\, 0,\infty)$. To 
strengthen it, set
$$
  h_t(\lambda,\tau,y)
  = t^{-1}\psi_{\lfloor\lambda U(t)\rfloor,\lfloor\tau U(t)\rfloor}
  \bigl(\chi(t)y\bigr) \, .
$$
We assume that
\setbox1=\vbox{\hsize=3.6truein\par\noindent
  $(\tau,y)\mapsto h_t(\,\cdot\,,\tau,y)$ and
  $(\tau,y)\mapsto h(\,\cdot\,,\tau,y)$ are measurable, and,
  for Lebesgue almost all $(\tau,y)$ in $[\, 0,\infty)^2$, 
  the functions $h_t(\,\cdot\,,\tau,y)$
  converge to $h(\,\cdot\,,\tau,y)$ in $\D[\, 0,\infty)$ as $t$ tends to 
  infinity.}
$$
  \raise -24pt\box1
  \eqno{\equa{PathA}}
$$

Though this is not important for our purpose,
assumption \PathA\ is not independent of (S1)
and there is a somewhat complicated though explicit relation between the
functions $h$ and $\rho$.

The following result describes the most likely trajectories of the process
leading to a large maximum.

\Theorem{\label{LimitingTraj}%
  Under (S1)--(S7) and \PathA, the conditional distributions of
  $(\fS_t,\tau_t,Y_t)$ given $M>t$ converges weakly$*$ to the distribution
  of $\bigl(\fS,\tau,\rho(\tau)Y\bigr)$ where $\tau$ and $Y$ are independent
  and

  \noindent (i) $Y$ has a Par\'eto distribution on $[\,1,\infty)$ with
  parameter $\alpha$,

  \noindent (ii) $\tau$ has density proportional to $\rho^{-\alpha}$,

  \noindent (iii) $\fS(\lambda)=h\bigl(\lambda,\tau,\rho(\tau)Y\bigr)$.
}

\bigskip

\Proof The proof requires establishing a couple of lemmas, describing
the limiting behavior of $(\tau_t,Y_t)$ given that $M$ exceeds $t$, as 
$t$ tends to infinity.

\Lemma{\label{TauYLimit}%
  Under the assumption of Theorem \LimitingTraj, for any nonnegative $u$ 
  and $y$
  $$
    \limt P\{\, \tau_t\leq u\,;\, Y_t>y\mid M>t \,\}
    = { \int_0^u\bigl( y\vee \rho(v)\bigr)^{-\alpha} \d v
        \over
        \int_0^\infty \rho(v)^{-\alpha}\d v } \, . 
  $$
}

\Proof Let $\epsilon$ and $\eta$ be two positive real numbers. 
Referring to the sets introduced in subsection \fixedref{3.1} 
and \fixedref{3.2}
and with $N=\lfloor U(t)/\epsilon\rfloor$, consider the event
$$
  F=\bcup_{i\in I} \bcup_{1\leq n\leq N} 
  A_i\cap B_n^\complement\cap C_n\cap D_n \, ,
$$
whose dependence on $t$, $\epsilon$ and $\theta$ is not kept track of. Note
that if $\theta$ is small enough and $t$ is large enough and if $F$ occurs,
the proof of the upper bound of Theorem \UsefulStatement\ shows that
$$
  P(\{\, M>t\,\}\setminus F)\leq 4\eta r(t) \, .
$$
Therefore, ultimately,
$$\displaylines{\qquad
  P\{\, \tau_t\leq u \,;\, Y_t>y\, ;\, M>t\,\}
  \hfill\cr\noalign{\vskip 3pt}\hfill
  \leq P\{\, \tau_t\leq u\, ;\, Y_t>y \,;\, M>t\, ;\,  F\,\} + 4\eta r(t) \, .
  \qquad
  \equa{TauYLimitA}\cr}
$$
We can also write $F$ as
$$
  F=\bcup_{1\leq n\leq N}\Bigl(\bcup_{i\in I}
  A_i\cap B_n^\complement\cap C_n\cap D_n\Bigr) \, .
$$
Take $\theta$ smaller than $y$ and $\epsilon$ sufficiently small so that $u$
lies between $\epsilon$ and $1/\epsilon$.
If $\tau_t\leq u$ and $Y_t>y$ and $B_n^\complement$ occur and if $J_t\not= n$,
then the two distinct random variables $X_{J_t}$ and $X_n$ exceed 
$\theta \chi(t)$ and both $J_t$ and $n$ are at most $N$. But we have seen
after (S5) that the probability for two distinct $X_i$ with $1\leq i\leq N$
to exceed $\theta\chi(t)$ is $o\bigl(r(t)\bigr)$ as $t$ tends to infinity.
Therefore, \TauYLimitA\ is ultimately at most $5\eta r(t)$ plus
$$\displaylines{\quad
  P\Bigl\{\, \bcup_{1\leq n\leq N} 
  \Bigl( \{\, \tau_t\leq u\,;\, J_t=n\,;\, Y_t>y\,;\, M>t\,\}
  \hfill\cr\hfill
         {}\cap \bcup_{i\in I} A_i\cap B_n^\complement\cap C_n\cap D_n\,\Bigr)
  \,\Bigr\}\, .
  \quad\cr}
$$
Using Bonferroni's inequality, ultimately, this is at most
$$
  \sum_{1\leq n\leq uU(t)} 
  P\{\, X_n\geq y\chi (t)\,;\, \bcup_{i\in I} A_i\cap D_n\,\} 
  + 5\eta r(t) \, .
$$
Next, on $A_i\cap D_n$,
$$
  t<S_i\leq \delta t + \psi_{i,n}(X_n) \, .
$$
Therefore, taking $\epsilon$ small enough so that $u$ is between $\epsilon$
and $1/\epsilon$, the right hand side in \TauYLimitA\ is at most
$$
  \sum_{1\leq n\leq uU(t)} P\Bigl\{\, X_n> y \chi(t) \,;\, \max_{i\in I}
  \psi_{i,n}(X_n)> (1-\delta)t \,\Bigr\} + 5 \eta r(t) \, .
$$
Since
$$
  \sum_{1\leq n\leq \epsilon U(t)} P\{\, X_n> y\chi(t)\,\}
  \sim \epsilon y^{-\alpha}(U\oF\circ\chi)(t)\, ,
$$
as $t$ tends to infinity, we see that, provided $\epsilon$ is small enough,
\TauYLimitA\ is at most
$$\displaylines{\quad
  \sum_{\epsilon U(t)\leq n\leq uU(t)}
  P\Bigl\{ X_n>y\chi (t)\,;\, \max_{i\in I}\psi_{i,n}(X_n)>(1-\delta)t\,\Bigr\}
  \hfill\cr\hfill
  {}+ 6\eta (U\oF\circ\chi)(t)\, ,
  \quad\cr}
$$
that is, at most
$$\displaylines{\quad
  \sum_{\epsilon U(t)\leq n\leq uU(t)}
  P\Bigl\{\, X_n> y\chi(t)\,;\, \psi_n(X_n)> (1-\delta)t\,\Bigr\}
  \hfill\cr\hfill
  {}+ 6\eta (U\oF\circ\chi)(t) \, .
  \quad\equa{TauYLimitB}\cr
  }
$$
Note that (S1) implies
$$
  \psi_{\lfloor u U(t)\rfloor}^{-1}\bigl( (1-\delta)t\bigr)
  \sim \chi\bigl( (1-\delta)t\bigr) \rho\bigl((1-\delta)^{-1/\gamma}u\bigr) 
  \, .
$$
Recall we set $\kappa$ for the index of regular variation of $\chi$.
The same argument used to derive \rDef, that is, regular variation and
comparison to a Riemann integral, shows that the sum in \TauYLimitB\ is
asymptotically equivalent to
$$
  (U\oF\circ\chi)(t)\int_\epsilon^u \Bigl( y\vee 
  (1-\delta)^\kappa\rho\bigl( v(1-\delta)^{-1/\gamma}\bigr)
  \Bigr)^{-\alpha} \d v\, .
$$
We let $\epsilon$, then $\delta$ and then $\eta$ tend to $0$ to obtain
$$
  \limsupt { P\{\, \tau_t\leq u \,;\, Y_t>y\,;\, M>t\,\} \over
             (U\oF\circ\chi)(t) }
  \leq \int_0^u \bigl( y\vee \rho(v)\bigr)^{-\alpha} \d v \, .
$$

To obtain a matching lower bound, we refer to how we proved the lower bound
of Theorem \UsefulStatement. In particular, keeping the notation of that proof
and remembering that the sets $B_n^\complement\cap C_n$ are disjoint for 
different values of $n$, we see that
$$
  P\{\, \tau_t\leq u\,;\, Y_t>y\,;\, M>t\,\} 
$$
is at least (cf.\ \HToThD)
$$\displaylines{\qquad
  \sum_{\ss n\in I_{\delta,t}\hfill\atop\ss n\leq uU(t)\hfill}
            P\bigl\{\, J_t=n\,;\, X_n> y\chi(t)\,;\, 
            F_n\cap B_n^\complement\cap C_n\cap D_n\,\bigr\} 
  \hfill\cr\hfill
  {}\geq \sum_{\ss n\in I_{\delta,t}\hfill\atop\ss n\leq uU(t)\hfill}
            P\bigl\{\, X_n> y\chi(t)\vee 
            \psi_n^{-1}\bigl( (1+\delta)t\bigr)\,\Bigr\}
  \hfill\cr\hfill
  {}-\sharp I_{\delta,t}\oF\bigl( \chi(t)/\theta\bigr) 
            + o\bigl(r(t)\bigr)
  \qquad\equa{TauYLimitC}\cr
  }
$$
as $t$ tends to infinity. Since (S1) implies that for $n$ in $I_{\delta,t}$
$$
  \psi_n^{-1} \bigl( (1+\delta)t\bigr)
  \sim \chi\bigl( (1+\delta)t\bigr) \rho\bigl( n/U\bigl((1+\delta)t\bigr)\bigr)
$$
and $\chi$ is regularly varying with index $\kappa$, the sum 
in \TauYLimitC\ is asymptotically equivalent to
$$\displaylines{\qquad
  \int_{\delta U(t)}^{uU(t)} \oF\Bigl( \chi\bigl( (1+\delta)t\bigr)
  \rho\Bigl({s\over U\bigl((1+\delta)t\bigr)}\Bigr)\vee \chi(t)y\Bigr) \d s
  \hfill\cr\hfill
  \eqalign{
  {}\sim{}&\int_\delta^u 
           \oF\Bigl( \chi(t)\bigl( (1+\delta)^\kappa
           \rho\bigl((1+\delta)^{-1/\gamma}v\bigr)\vee y\bigl)\Bigr)
           U(t)\d v \cr
  {}\sim{}&(U\oF\circ\chi)(t)\int_\delta^u \bigl( (1+\delta)^\kappa
           \rho\bigl((1+\delta)^{-1/\gamma}v\bigr)\vee
           y\bigl)^{-\alpha} \d v \cr} 
  \qquad\cr}
$$
as $t$ tends to infinity. Since $\delta$ can be made as small as desired,
it follows that
$$
  \liminft { P\{\, \tau_t\leq u \,;\, Y_t>y\,;\, M>t\,\}
             \over (U\oF\circ\chi)(t) }
  \geq \int_0^u \bigl( y\vee\rho(v)\bigr)^{-\alpha}\d v 
$$
and this completes the proof.

\bigskip

The next lemma gives a simple 
representation of random variables having the limiting distribution 
written in Lemma \TauYLimit.

\Lemma{\label{cleanDistribution}
  Let $\tau$ be a random variable having density proportional to 
  $\rho^{-\alpha}$, and let $Y$ be a random variable independent of
  $\tau$, having a Par\'eto distribution of index $\alpha$ on $[\,1,\infty)$.
  Then
  $$
     P\{\,\tau\leq u\,;\, \rho(\tau)Y>y\,\}
     = {\int_0^u \bigl( y\vee\rho(v)\bigr)^{-\alpha}\d v \over
        \int_0^\infty \rho(v)^{-\alpha} \d v } \, .
  $$
}

\Proof We simply write
$$\displaylines{\qquad
    P\{\,\tau\leq u\,;\, \rho(\tau)Y>y\,\}
       \int_0^\infty \rho(v)^{-\alpha} \d v
  \hfill\cr\hfill
  \eqalign{
    {}={}&\int_0^u P\{\, Y>y/\rho(v)\,\} \rho(v)^{-\alpha}\d v \cr
    {}={}&\int_0^u\Bigl( {y\over\rho(v)}\vee 1\Bigr)^{-\alpha} 
          \rho(v)^{-\alpha}\d v \cr
    {}={}&\int_0^u \bigl( y\vee \rho(v)\bigr)^{-\alpha} \d v \, . 
    \quad\qquad\qquad\qed\cr}%
  \cr}
$$

We can now conclude the proof of Theorem \LimitingTraj. Let $\epsilon$ be a 
positive real number and consider the event
$$
  G=\bigl\{\, \exists\lambda \leq 1/\epsilon \, :\, 
  |\fS_t(\lambda)-h_t(\lambda,\tau_t,Y_t)|>\delta\,\bigr\} \, .
$$
Recall that the event $F$ introduced in the proof of the upper bound related
to Lemma \TauYLimit\ depends on a paramater $\theta$ through the events
$B_n$ and $C_n$. As we saw in that proof, for any fixed $\eta$ and for any
$\theta$ small enough
$$
  P(G\,;\, M>t)
  \leq P(G\cap F\,;\, M>t)+4\eta r(t)
$$
ultimately. If $F$ occurs, then $N_t$ is at most $N$ and, also,
$\bcup_{1\leq n\leq N}B_n^\complement\cap C_n \cap D_n$ occurs. Again,
this union is disjoint for the sets $B_n^\complement\cap C_n$ are
disjoint for different values of $n$. If $B_n^\complement\cap C_n$
occurs, it is tempting to conclude that $J_t=n$. Since $J_t$ is at
most $N_t$, this is true provided that $N_t$ is at least $n$.
We now show that assumption (S5) guarantees that this is likely to be the
case. Clearly, from the definition of the event $F$,
$$\displaylines{\qquad
    F\subset
    \bcup_{1\leq n\leq N}\bcup_{i<n}A_i\cap B_n^\complement\cap C_n\cap D_n
  \hfill\cr\hfill
    {}\cup\, \bcup_{1\leq n\leq N} \bcup_{n\leq i\leq N} 
    \Bigl( A_i\cap B_n^\complement\cap C_n\cap D_n\cap 
    \bigl(\bcup_{j<n}A_j\bigr)^\complement\Bigr) \, .\qquad\cr}
$$
Note that if $A_i\cap B_n^\complement\cap C_n\cap D_n$ occurs, then so does
$A_i\cap C_n$. In this case, if $i$ is less than $n$ then the definition of
$C_n$ shows that $A_i\cap\bcap_{1\leq j\leq i} B_j$ occurs as well. But (S5)
combined with (S3) imply 
that $\bcup_{1\leq i\leq N}( A_i\cap \bcap_{1\leq j\leq i}B_j)$ occurs with 
probability at most $2\delta r(t)$ provided $\epsilon$ is small enough and $t$
is large enough. Therefore,
\hfuzz=2pt
$$\displaylines{
  P(G\,;\, M>t)
  \hfill\cr\hfill
  \leq P\Bigl(G\cap \bcup_{1\leq n\leq N}\bcup_{n\leq i\leq N}\bigl(A_i\cap 
  B_n^\complement
  \cap C_n\cap D_n\cap\bigl(\bcup_{j<n}A_j\bigr)^\complement\bigr)\Bigr)
  +6\eta r(t)\, .
  \cr}
$$
\hfuzz=0pt
Considering the events involved in this upper bound,
if $A_i\cap (\bcup_{j<n}A_j)^\complement$
occurs and $i$ is at least $n$, then
$N_t$ is at least $n$ and, as anounced, $J_t=n$; moreover, if $D_n$ occurs, 
then $|S_i-\psi_{i,n}(X_n)|$ is at most $\delta t$ for any $i$ at most
$N$. In that case, uniformly in $\lambda$ in $[\,0,1/\epsilon\,]$,
$$
  \bigl|t^{-1}S_{\lfloor \lambda U(t)\rfloor} 
  - t^{-1}\psi_{\lfloor\lambda U(t)\rfloor,J_t}(X_{J_t})\bigr|
  \leq \delta \, ,
$$
that is,
$$
  |\fS_t(\lambda)-h_t(\lambda,\tau_t,Y_t)|\leq\delta \, ,
$$
and so $G$ does not occur. Therefore,
$$
  \limt P(G\mid M>t )=0\, .
  \eqno{\equa{TrajLimitingA}}
$$
Recall that the Skorohod topology is metric.
Combining Lemmas \TauYLimit, \cleanDistribution\ and Theorem 5.5 
in Billingsley (1968) upon
using \PathA, imply that the conditional distribution of
$h_t(\,\cdot\,,\tau_t, Y_t)$ given $M>t$ converges weakly$*$ to that of
$h\bigl(\,\cdot\,,\tau,\rho(\tau)Y\bigr)$ as $t$ tends to infinity. The
result then follows from \TrajLimitingA\ which asserts that, on
compact sets, $\fS_t-h_t(\,\cdot\,,\tau_t,Y_t)$ converges uniformly to
$0$ in probability under the conditional probability given
$M>t$ as $t$ tends to infinity, and uniform convergence on compact
sets implies convergence in $\D[\, 0,\infty)$ under the Skorohod 
metric.\hfill\qed

\bigskip


\section{A large deviation inequality and a Karamata type theorem}
The folklore attributes to Kolmogorov that behind every limit theorem
there is an inequality. The purpose of this short section is to derive
the large deviation inequality behind some of our theorems as well as
to state a Karamata type theorem which we will be needing.

\bigskip


\subsection{A large deviation inequality}%
The result of this subsection is of a more technical nature. 
It provides a bound on the
moment generating function of a centered random variable truncated from
above. Its use is explained after its proof, and it will be instrumental
to show that some probabilities tend to $0$. It is inspired from 
a technique used in Cline and Hsing (1991) as well as Ng, Tang, Yan and 
Yang (2004).

\Lemma{\label{MGF}
  Let $Z$ be a centered random variable with distribution function $H$. 
  For any positive $\lambda$ and $a$, for any positive $\eta$ less than $1$,
  $$\displaylines{\quad
    \log E\exp \bigl({\lambda Z}\II\{\, Z\leq a\,\}\bigr)
    \leq \eta\lambda E|Z| 
    -\lambda EZ\II\{\, \lambda Z\leq\log (1-\eta)\,\}
    \hfill\cr\noalign{\vskip 3pt}\hfill
    {}+ e^{\lambda a} \oH \Bigl({\log (1+\eta)\over \lambda}\Bigr) \, .
    \quad\cr}
  $$
  }

\Proof Let $M(\lambda)$ be the expected value of 
$\exp\bigl(\lambda Z\II\{\,Z\leq a\,\}\bigr)$. The inequality $1+x\leq e^x$ yields
$$\eqalign{
  M(\lambda )
  &{}=1+\int_{-\infty}^a e^{\lambda z} -1\d H(z) \cr
  &{}\leq\exp\Bigl( \int_{-\infty}^a e^{\lambda z}-1\d H(z)\Bigr) \, . \cr
  }
$$
Since $e^x-1$ is nonnegative and at most $(1+\eta)x$ on 
$[\,0,\log (1+\eta)\,]$ and is nonpositive and at most $(1-\eta)x$ 
on $[\,\log (1-\eta),0\,]$, we have
$$\eqalign{
  \log M(\lambda)
  &{}\leq \int_{-\infty}^a e^{\lambda z}-1 \d H(z) \cr
  &{}\leq (1+\eta) \int_0^\infty 
    \lambda z\II\{\, \lambda z\leq \log (1+\eta)\,\}\d H(z) \cr
  &\qquad {}+(1-\eta)\int_{-\infty}^0 
    \lambda z\II\{\, \log (1-\eta)\leq\lambda z\,\}\d H(z) \cr
  &\qquad {}+ e^{\lambda a} \oH\Bigl( {\log (1+\eta)\over\lambda}\Bigr) 
    \, . \cr
  }
$$
Since $Z$ is centered,
$$\displaylines{\qquad
  \int_{\RR}\lambda z 
  \II\{\, \log (1-\eta)\leq \lambda z\leq\log (1+\eta)\,\} \d H(z)
  \hfill\cr\hfill
  \leq -\int_{-\infty}^{\log (1-\eta)/\lambda} \lambda z  \d H(z) \, ,
  \qquad\cr}
$$
and the result follows.\hfill\qed

\bigskip

We will use Lemma \MGF\ in the following situation. Consider a sequence
$(a_j)_{j\geq 1}$ of positive real numbers and a sequence of independent
and equidistributed and centered random variables $(Z_j)_{j\geq 1}$. 
We write $Z$ for a random variable having the same distribution as $Z_1$. 
Substituting $\lambda$ with $\lambda g_j$ and $a$ with $a_{n-j}$
in Lemma \MGF\ we obtain the Chernoff type inequality, valid for all positive
$\eta$ less than $1$, any positive $\lambda$, all sequences of positive 
reals $a_j$, nonnegative $g_j$, and any real numbers $t$ and $s_n$ 
such that $t-s_n$ is positive,
$$\displaylines{
  \log P\Bigl\{\, \sum_{0\leq j< n} g_jZ_{n-j}\II\{\, Z_{n-j}\leq a_{n-j}\,\}
                  >t-s_n 
               \,\Bigr\}
  \hfill\cr\hfill
  {}\leq -\lambda (t-s_n)
  +\eta \lambda \g{n} E|Z|
  -\lambda \sum_{0\leq j<n} g_j 
  EZ\II\{\, \lambda g_j Z\leq \log (1-\eta)\,\}
  \cr\hfill
  {}+ \sum_{0\leq j< n}e^{\lambda g_j a_{n-j}}
  \oH\Bigl( {\log (1+\eta)\over \lambda g_j}\Bigr) \, .
  \qquad\qquad\equadef{LargeDev}
  \cr}
$$
We will take $\lambda$ small in this bound. The a priori strange formulation
of this inequality, using a rather mysterious $t-s_n$ instead of a single 
variable is on purpose and designed to make the remainder of this paper 
easier to read.

\bigskip


\subsection{A Karamata type theorem}%
The following result is an easy extension of the direct half of Karamata's 
theorem (see Bingham, Goldie and Teugels, 1989, Proposition 1.5.10).
Recall that if $b$ is a function with limit infinity at infinity, 
then, if it exists, $b^\leftarrow$ is an asymptotic inverse of $b$, that
is, a function such that $b\circ b^\leftarrow\sim b^\leftarrow \circ b\sim \Id$
at infinity.

\Lemma{\label{Karam}%
  Let $b$ be a regularly varying function of positive index $\beta$. If
  $\alpha>1\vee (1/\beta)$ then for any positive real number $r$
  $$
    \int_r^\infty \oF\bigl( t+b(u)\bigr) \d u
    \sim {1\over\beta} B\Bigl({1\over\beta},\alpha-{1\over\beta}\Bigr)
    (b^\leftarrow\oF)(t) 
  $$
  as $t$ tends to infinity.
}

\bigskip

\Proof Let $\epsilon$ be a positive real number. The change of variable
$u=\lambda b^\leftarrow(t)$ and regular variation of $\oF$ and $b$
show that
$$
  \int_{\epsilon b^\leftarrow(t)}^{b^\leftarrow(t)/\epsilon} 
  \oF\bigl( t+b(u)\bigr)\d u
  \sim (b^\leftarrow\oF)(t) \int_\epsilon^{1/\epsilon} 
  (1+\lambda^\beta)^{-\alpha} \d\lambda \, .
$$
Next, since $b$ is ultimately positive, tends to infinity at infinity, so
does $b^\leftarrow$. Therefore, by the monotonicity of $\oF$,
$$
  \int_r^{\epsilon b^\leftarrow(t)} \oF\bigl( t+b(u)\bigr)\d u
  \leq \epsilon (b^\leftarrow\oF)(t) \bigl(1+o(1)\bigr)
$$
as $t$ tends to infinity.
Furthermore, by Karamata's theorem and monotonicity of $\oF$,
$$\eqalign{
  \int_{b^\leftarrow(t)/\epsilon}^\infty \oF\bigl(t+b(u)\bigr)\d u
  &{}\leq \int_{b^\leftarrow(t)/\epsilon}^\infty \oF\circ b(u)
   \d u \cr
  &{}\sim {\epsilon^{\alpha\beta-1}\over\alpha\beta-1} (b^\leftarrow \oF)(t) 
   \, .
   \cr}
$$
Since $\alpha\beta>1$,
$$
  \int_r^\infty \oF\bigl(t+b(u)\bigr) \d u
  \sim (b^\leftarrow\oF)(t) \int_0^\infty (1+\lambda^\beta)^{-\alpha}
  \d\lambda
$$
as $t$ tends to infinity. The change of variable $x=\lambda^\beta$ then 
yields the result.\hfill\qed

\bigskip


\section{Some asymptotic analysis related to analytic functions}%
The purpose of this section is prove some purely analytical results
related to analytic functions which will be needed to prove our results
on the maximum of $(g,F)$-processes and their trajectories. 

In the first subsection we restate known results in a form suitable for our
purpose. The second subsection introduces a family of functions, $\Psi_n$,
$n\geq 1$, associated to an analytic function. The notation is not fortuitous,
for if $S_n$ is a $(g,F)$-process, then the function $\psi_n(x)=\max_{i\geq 0}
E(S_i|X_n=x)$ involved in our heuristic will be related to the function
$\Psi_n$ in a simple way. In order to apply the methodology presented in
section \fixedref{3}, we need to have some information on $\psi_n^{-1}$,
and for this reason, we will obtain some basic information on 
$\Psi_n^{-1}$. This will be done in the third subsection when the 
sequence of Taylor coefficients $(g_n)_{n\geq 0}$ tends to $0$, in the
fourth subsection when that sequence diverges toward infinity, and in the
fifth subsection when $(g_n)_{n\geq 0}$ has a positive finite limit.

\bigskip

Throughout this subsection, $g(x)=\sum_{i\geq 0} g_ix^i$ is a real analytic 
function on $(-1,1)$, regularly varying 
at $1$ with positive index $\gamma$. In particular,
$\lim_{x\to 1; x<1}g(x)=+\infty$. Recall that for any positive
$n$, we write $\g{n}$ for $\sum_{0\leq j< n} g_j$.

\bigskip

The following notation will save some 
unsightful $\epsilon$, $t_0$, $n_0$ as well as various quantifiers.

\bigskip

\noindent{\bf Notation.} Throughout this section, if $(a_n)$ and $(b_n)$
are two sequences, we say that `$a_n$ is bounded from above by
an equivalent of $b_n$' and write $a_n\lsim b_n$  if $a_n\leq b_n 
\bigl(1+o(1)\bigr)$ as $n$ tends to infinity, or, equivalently, if
$\limsup_{n\to\infty} a_n/b_n\leq 1$. Similarly we define in
an obvious manner what it is to be asymptotically bounded from
below, and write $\gsim$ for this relation. Both relations are transitive.

\bigskip


\subsection{Preliminaries}%
In this section we restate some known result in a form suitable for
our analysis.

Our first lemma essentially restates Karamata's Tauberian theorem for
power series (Bingham, Goldie and Teugels, 1989, Corollary 1.7.3), 
and adds some uniformity to it. Recall that throughout this paper we
assume that $(g_n)_{n\geq 0}$ is asymptotically equivalent to a 
monotone sequence.

\Lemma{\label{BS}%
  The following asymptotic equivalences hold as $n$ tends to infinity, 
  uniformly in $x$ in any compact subset of the positive half-line,

  \smallskip

  \noindent (i) $\ds g_{\lfloor nx\rfloor}
  \sim {\ds x^{\gamma-1}\over\ds \Gamma(\gamma)}\,{\ds g(1-1/n)\over\ds n}$,

  \smallskip

  \noindent (ii) $\ds \g{nx}
  \sim {\ds x^\gamma\over\ds\Gamma(1+\gamma)}\, g(1-1/n)$.
}

\bigskip

\Proof For a fixed $x$, the result is Corollary 1.7.3 in Bingham, Goldie
and Teugels (1989). 
Uniformity follows by the same proof or the following one. We note that for
any $\epsilon$ positive,
$$
  \sup_{\epsilon\leq x\leq 1/\epsilon} 
  \Bigl| { g_{\lfloor nx\rfloor}\lfloor nx\rfloor \Gamma(\gamma)\over
           g(1-1/\lfloor nx\rfloor) } -1 \Bigr|
  \leq
  \sup_{m\geq \epsilon n} \Bigl| {g_mm\Gamma(\gamma)\over g(1-1/m)} -1\Bigr|
  \, . 
$$
The pointwise version implies that this upper bound tends to $0$ as $n$ tends
to infinity. Since $\lfloor nx\rfloor/nx$ tends to $1$ uniformly in $x$ in
any compact set of the positive half-line, and by the
uniform convergence theorem (Bingham, Goldie and Teugels, 1989, Theorem 1.5.2)
$g(1-1/\lfloor nx\rfloor)/g(1-1/n)$ tends to $x^\gamma$ uniformly as well,
it follows that
$$
  \lim_{n\to\infty}\sup_{\epsilon\leq x\leq 1/\epsilon} 
  \Bigl| { ng_{\lfloor nx\rfloor}\over g(1-1/n) } - 
         { x^{\gamma-1}\over\Gamma(\gamma) } \Bigr| = 0 \, .
$$
A similar argument proves the uniformity in the convergence 
of $\g{nx}/g(1-1/n)$.\hfill\qed

\bigskip


\subsection{The functions \poorBold{$\Psi_n$} and their inverses}%
In this subsection we develop some asymptotic estimates for some functions
derived from a real analytic function on $(-1,1)$ with a singularity
at $1$. Some of the results presented may be of independent interest
and fit in the rich corpus of Tauberian theorems in the realm of
analytic functions.

Recall that $U$ is a function which satisfies $g\bigl(1-1/U(t)\bigr)
\sim t$ as $t$ tends to infinity, and that $g_n$, $n\geq 0$, is the sequence
of the Taylor coefficients of $g$ at the origin. In this section, we do not
assume that these coefficients are nonnegative, but that only finitely
many of them may be negative. We also assume that
$$
  g^*=\sup_{n\geq 0} g_n \quad \hbox{is positive.}
$$

We consider the following functions,
$$
  \Psi_n(x)=0\vee \max_{k\geq 0} (g_k x -\g{n+k})\, ,\qquad n\geq 1\,.
$$
Since Lemma \BS\ implies that $\lim_{k\to\infty}\g{n+k}/g_k=\infty$,
these functions are defined for all $n$ and $x$ nonnegative, that is
to say, since the maximum is attained, it is proper to write a maximum
instead of a supremum.  Lemma \BS\ implies that there exists $n_0$
such that for any $n$ at least $n_0$, both $g_n$ and $\g{n}$ are
nonnegative. Clearly, if all the coefficients $g_n$ are nonnegative,
we can take $n_0$ to be $0$.  For $n$ at least $n_0$, our first lemma
gives another expression for $\Psi_n$ and shows that this function is
increasing and convex on the half-line where it is positive. For this
purpose, we define
$$
  K_+=\{\, k\in\NN \, : g_k>0 \,\} \, ,
$$
and
$$
  x_n=\min_{k\in K_+} \g{n+k}/g_k \, .
$$

\Lemma{\label{psiMonotone}%
  For any $n$ at least $n_0$,

  \smallskip

  \noindent (i) $\Psi_n$ vanishes on $[\,0,x_n\,]$;

  \smallskip

  \noindent (ii) $\Psi_n$ is positive and increasing
  on $(x_n,\infty)$. Moreover, on the half-line $(x_n,\infty)$,
  $$
    \Psi_n (x)=\max_{k\in K_+} g_kx-\g{n+k} \, ;
  $$

  \smallskip

  \noindent (iii) $\Psi_n$ is continuous and convex on the nonnegative
  half-line;

  \smallskip

  \noindent (iv) If $n_0\leq n\leq m$, then $\Psi_m\leq \Psi_n$ 
  and $x_n\leq x_m$.
}

\bigskip

\Proof (i) Let $x$ be a nonnegative real number at most equal to $x_n$. The
definition of $x_n$ implies that $g_k x-\g{n+k}$ is nonpositive for any $k$ 
in $K_+$. If $k$ does not belong to $K_+$ then $g_kx-\g{n+k}$ is 
nonpositive, for both $g_k$ and $-\g{n+k}$ are nonpositive. Therefore,
$\Psi_n$ vanishes at $x$. 

\noindent (ii) If $x$ is larger than $x_n$, then $g_k x-\g{n+k}$ is positive
for some $k$ in $K_+$ and so is $\Psi_n(x)$. Since $g_kx-\g{n+k}$ is 
nonpositive for
$k$ not in $K_+$, this proves that $\Psi_n$ has the representation given in 
(ii). This representation shows that $\Psi_n$ is increasing 
on $(x_n,\infty\,]$.

\noindent (iii) The representation obtained in (ii) and the proof of (i)
show that $\Psi_n(x)=0\vee\max_{k\in K_+}g_kx-\g{n+k}$ on the nonnegative
half-line. As the supremum of nondecreasing linear functions, $\Psi_n$
is convex and therefore continuous.

\noindent (iv) If $n_0\leq n\leq m$ then $\g{n+k}\leq \g{m+k}$ for any
$k$ and the result follows.\hfill\qed

\bigskip

For $n$ at least $n_0$, the expression for $\Psi_n$ in Lemma
\psiMonotone.ii shows that $\lim_{x\to\infty} \Psi_n(x)=+\infty$.
Lemma \psiMonotone.ii--iii imply that $\Psi_n$ is invertible as a map
from $(x_n,\infty)$ to the positive half-line. Therefore, for $n$ at
least $n_0$, it is meaningful to define the inverse $\Psi_n^{-1}$ on
the positive half-line. We extend it to $0$ by continuity, defining
$\Psi_n^{-1}(0)=x_n$. The following lemma provides an expression for
that inverse.

\Lemma{\label{psiInv}%
  Let $n$ be at least $n_0$. For any nonnegative $t$,
  $$
    \Psi_n^{-1}(t)=\min_{k\in K_+} {t+\g{n+k}\over g_k} \, .
  $$
}

\Proof Since $\Psi_n$ is invertible, for any positive $t$ and any $k$ in $K_+$,
$$
  t=\Psi_n\circ\Psi_n^{-1} (t)\geq g_k\Psi_n^{-1}(t)-\g{n+k} \, .
$$
Therefore,
$$
  \Psi_n^{-1}(t)\leq \min_{k\in K_+} {t+\g{n+k}\over g_k} \, .
$$
To prove that this upper bound is sharp, assume that it is not, so that
there exists a positive $\epsilon$ such that for any positive $g_k$,
$$
  \Psi_n^{-1}(t)\leq {t+\g{n+k}\over g_k} -\epsilon \, .
$$
Since $\Psi_n$ is onto, there exists $x$ such that $\Psi_n(x)=t$. Then
$$
  x=\Psi_n^{-1}\circ\Psi_n(x) \leq {\Psi_n(x)+\g{n+k}\over g_k}-\epsilon \, ,
$$
and therefore
$$
  \Psi_n (x)\geq x g_k -\g{n+k} +\epsilon g_k \, .
$$
In particular, considering this inequality for a value of $k$ which 
maximizes $xg_k-\g{n+k}$, this last quantity being then equal to $\Psi_n(x)$,
we obtain
$$
  \Psi_n(x) 
  \geq \Psi_n(x)+\epsilon g_k \, .
  \eqno{\equa{psiInvA}}
$$
If $g_k$ were equal to
$0$ then $t=\Psi_n(x)=-\g{n+k}$ would be negative, for $n$ is at least $n_0$.
Therefore, $g_k$ is positive and \psiInvA\ yields $\Psi_n(x)>\Psi_n(x)$
which is a contradiction.\hfill\qed

\bigskip

Given Lemma \psiInv, we write $k_n(t)$ for an integer such that
$$
  \Psi_n^{-1}(t)= {t+\g{n+k_n(t)}\over g_{k_n(t)}} \, .
$$
Such an integer may not be unique, but whatever statement we will make about
it will not depend on its particular choice.

\bigskip


\subsection{Approximation of \poorBold{$\Psi_n^{-1}$} when 
\poorBold{$(g_n)_{n\geq 0}$} tends to 0}%
When the sequence $(g_n)_{n\geq 0}$ tends to $0$ at infinity, the minimization
involved in Lemma \psiInv\ can be made explicit for large argument $t$. 
For this purpose, let $k^*$ be the smallest integer at which the sequence
$(g_n)_{n\geq 0}$ achieves its maximum; thus $k^*$ is the smallest integer
for which $g_k=g^*$.

\Lemma{\label{PsiInvFormulag}%
  There exists a nonnegative $t_0$ such that for any $n$ at least $n_0$
  and any $t$ at least $t_0$,
  $$
    \Psi_n^{-1}(t)={t+\g{n+k^*}\over g^*} \, .
  $$
}

\Proof Given Lemma \psiInv, it is clear that the proposed expression is
an upper bound for $\Psi_n^{-1}$. If $k$ is in $K_+$ and larger than
$k^*$, then $g^*/g_k\geq 1$ and $g_{[n+k^*,n+k)}\geq 0$. Therefore,
$$\eqalign{
  {t+\g{n+k}\over g_k}
  &{}={t+\g{n+k^*}\over g^*} {g^*\over g_k} + {g_{[n+k^*,n+k)}\over g_k} \cr
  &{}\geq {t+\g{n+k^*}\over g^*} \, . \cr
  }
$$
If $k$ is in $K_+$ and less than $k^*$ we write $(t+\g{n+k})/g_k$ as the 
sum of $(t+\g{n+k^*})/g^*$ and
$$
  {t+\g{n+k^*}\over g^*} \Bigl( {g^*\over g_k}-1\Bigr)
  - {g_{[n+k,n+k^*)}\over g_k} \, .
$$
This quantity is positive if
$$
  t> g_{[n+k,n+k^*)} {g^*\over g^*-g_k} - \g{n+k^*} \, .
$$
Since $k^*$ is minimal, the maximum of this lower bound over all $k$ in $K_+$
and less than $k^*$ and all $n\geq n_0$ is finite. Call $t_0$ its 
maximum.\qed

\bigskip

It follows from Lemmas \BS.ii and \PsiInvFormulag\ that
$$
  \Psi_n^{-1}(t)\sim {t+\g{n}\over g^*}
  \eqno{\equa{PsiInvAsymp}}
$$
as $n$ tends to infinity, uniformly in $t$ in $[\, t_0,\infty)$. We then 
deduce the following asymptotic equivalence.

\Lemma{\label{PsiInvUg}%
  For any $x$ and $c$ in any compact subset of the positive half-line, as
  $t$ tends to infinity,
  $$
    \Psi_{\lfloor xU(t)\rfloor}^{-1} (t)
    \sim {t\over g^*} \Bigl( 1+{x^\gamma\over\Gamma(1+\gamma)}\Bigr) \, ,
  $$
  and, for any $n$ at least $n_0$ and $t$ at least $t_0$, the equality
  $k_n(t)=k^*$ holds.
}

\bigskip

\Proof This is immediate from Lemmas \BS\ and \PsiInvFormulag.\hfill\qed

\bigskip

Our next two lemmas provide some bounds for $\Psi_n^{-1}$. 

\Lemma{\label{PsiInvLbgA}%
  For any $n$ at least $n_0$ and any $t$ at least $t_0$,
  $$
    \Psi_n^{-1} (t)\geq t/g^* \, .
  $$
}

\Proof The result follows from the formula in Lemma \PsiInvFormulag\ since
$\g{n+k^*}$ is nonnegative.\hfill\qed

\bigskip

\Lemma{\label{PsiInvLbgB}%
  As $n$ tends to infinity
  $$
    \inf_{t\geq t_0}\Psi_n^{-1} (t)
    \gsim {g(1-1/n)\over g^*\Gamma(1+\gamma)} \, .
  $$
}

\Proof This follows from Lemmas \BS\ and \PsiInvFormulag, since
$t$ is at least $t_0$ and hence positive.\hfill\qed

\bigskip


\subsection{Asymptotic analysis of \poorBold{$\Psi_n^{-1}$} 
when \poorBold{$(g_n)_{n\geq 0}$} tends to infinity}%
The purpose of this subsection is to 
derive an asymptotic equivalent for $\Psi_n^{-1}$ 
and $k_n(\,\cdot\,)$ when
$n$ is of order $U(t)$ and the argument is of order $t$, as well as some 
bounds for this function when the sequence of coefficients, $(g_n)_{n\geq 0}$
tends to infinity.

When $\gamma$ is $1$, we assume that
\setbox1=\vbox{\hsize=3.2in\noindent
  $(g_n)_{n\geq 0}$ is asymptotically equivalent to an increasing 
  sequence with limit $+\infty$}
$$
  \box1
  \eqno{\equa{PrelimA}}
$$
and moreover that
\setbox1=\vbox{\hsize=3.2in\noindent
  the Karamata representation of $g$ satisfies \HypGOneB\ and \HypGOneC.}
$$
  \box1
  \eqno{\equa{PrelimB}}
$$

For any positive real number $a$, we 
introduce the functions
$$
  \xi_a(x)=\min_{y> 0} {a+(x+y)^\gamma\over\gamma y^{\gamma-1}} \, .
$$
Since $\gamma$ is fixed throughout this section, the notation $\xi_a$ does
not keep track of the dependence of this function on $\gamma$. In general
an explicit form for the minimum cannot be found, however, if $\gamma$ is $1$,
then $\xi_a(x)=a+x$.
With respect to Theorem \ThG, note that the function $\rho_\gamma$ is
equal to $\xi_{\Gamma (1+\gamma)}$. 

The following shows that the minimum involved in the definition of $\xi_a(x)$
is achieved at a unique point and that all the functions $\xi_a$ can be 
recovered from a knowledge $\xi_1$.

\Lemma{\label{xiA}%
  (i) The function 
  $$
    y\mapsto {a+(x+y)^\gamma\over\gamma y^{\gamma-1}}
  $$ 
  has a unique minimum on the positive half-line. 

  \smallskip

  \noindent (ii) If $a\leq b$, then $\xi_a\leq\xi_b$.

  \smallskip

  \noindent (iii) The identity $\xi_a(x)=a^{1/\gamma}\xi_1(a^{-1/\gamma}x)$
  holds.
}

\bigskip

\Proof (i) The derivative of the function vanishes at the minimum. Thus,
the minimizer $y$ satisfies
$$
  \Bigl({x+y\over y}\Bigr)^{\gamma-1} - {\gamma-1\over\gamma}
  {a+(x+y)^\gamma\over y^\gamma} = 0 \, .
$$
Setting $s=y/(x+y)$, this equation asserts that
$$
  s= {\gamma-1\over\gamma} \Bigl( a {(1-s)^\gamma\over x^\gamma}+1\Bigr) \, .
$$
The left hand side of this equality is increasing in $s$ in $(0,1\,]$ while
the right hand side is decreasing. Therefore, the equality is achieved
for a unique $s$.

\noindent (ii) This is obvious. 

\noindent (iii) We substitute $y$ with 
$a^{-1/\gamma}y$ in the definition of $\xi_1(a^{-1/\gamma}x)$.\hfill\qed

\bigskip

Knowledge of the behavior of $\xi_a$ at the origin and at infinity will also
be useful and some information is given now.

\Lemma{\label{xiAsymp}%
  The following hold for $\gamma$ greater than $1$.
  \smallskip

  \noindent (i) $\xi_a(0)=a^{1/\gamma}(\gamma-1)^{(1/\gamma)-1}$.

  \smallskip

  \noindent (ii) $\xi_0(1)=\bigl({\ds\gamma\over\ds\gamma-1}\bigr)^{\gamma-1}$.

  \smallskip
  \noindent (iii) $\xi_a(x)\sim x\xi_0(1)$ as $x$ tends to infinity.

  \smallskip
  
  Moreover, if $\gamma=1$, then $\xi_a(x)=a+x$ and (i)--(iii) hold provided 
  they are extended by continuity as $\gamma$ tends to $1$.
  }

\bigskip

\Proof Assume that $\gamma$ is greater than $1$.
Standard calculus shows that the function to minimize to calculate
$\xi_a(0)$ achieves its minimum at $y=a^{1/\gamma}(\gamma-1)^{1/\gamma}$, 
while that to calculate $\xi_0(1)$ achieves its minimum at $y=\gamma-1$. 
Parts (i) and (ii) follow.

To prove (iii), Lemma \xiA.iii yields, with $xa^{1/\gamma}$ in place of $x$ 
and setting $a=x^{-\gamma}$, the identity
$$
  \xi_1(x)=a^{-1/\gamma}\xi_a(xa^{1/\gamma}) = x\xi_{x^{-\gamma}}(1) \, .
$$
It suffices to prove that $\xi_a(1)$ tends to $\xi_0(1)$ as $a$ tends
to $0$ from above. On the one hand, if $a$ is positive, Lemma \xiA.ii shows
that $\xi_a\geq \xi_0$, and on the other hand, with $y=\gamma-1$,
$$
  \xi_a(1)
  \leq {a+\gamma^\gamma\over\gamma (\gamma-1)^{\gamma-1}}
$$
and this upper bound tends to $\xi_0(1)$ as $a$ tends to $0$.\hfill\qed

\bigskip

We can now derive an asymptotic equivalent for $\Psi_n^{-1}$ when $n$ is
of order $U(t)$ and the argument is of order $t$. To proceed, we note that
given Lemma \xiA.i, it is legitimate to define $\kappa_a(x)$ as the unique
positive real number such that
$$
  \xi_a(x)={a+(x+\kappa_a(x))^\gamma\over\gamma\kappa_a(x)^{\gamma-1}} \, .
$$
If $\gamma$ is greater than $1$, then
$$
  \lim_{\epsilon\to 0} {a+(x+\epsilon)^\gamma\over\gamma\epsilon^{\gamma-1}}
  =+\infty \, .
$$
Thus, for $\gamma$ greater than $1$, the function $\kappa_a$ maps compact
subsets of the positive half-line to compact subsets of the positive
half-line. In particular, on any compact subset of the positive half-line, 
$\kappa_a$ is lower bounded by a positive constant.

We remark that if $\gamma=1$, taking $0^0$ to be $1$, we have
$\kappa_a(x)=0$.
In the next result, the function $\varepsilon$ is that involved in
the Karamata representation of $g$. Recall that \PrelimB\ holds, that is
\HypGOneB\ and \HypGOneC\ are assumed to hold when $\gamma$ is $1$.

\Lemma{\label{psiInvU}%
  The following asymptotic equivalents hold uniformly in $x$ and $c$ in 
  any compact subset of the positive half-line as $t$ tends to infinity,
  $$
    \Psi_{\lfloor xU(t)\rfloor}^{-1}(ct)\sim U(t)\xi_{c\Gamma(1+\gamma)}(x)
    \, ,
  $$
  and
  $$
    k_{\lfloor xU(t)\rfloor}(ct)\sim
    \cases{ \kappa_{c\Gamma(1+\gamma)}(x) U(t)& if $\gamma>1$,\cr
            \noalign{\vskip 4pt}
            (c+x)(\Id\,\varepsilon)\circ U(t)   & if $\gamma=1$.\cr}
  $$
}

\Proof We distinguish two cases, according to whether $\gamma$ is greater than
$1$ or not.

\noindent{\it Case $\gamma>1$}.
Write $n=\lfloor xU(t)\rfloor$. In the minimization defining 
$\Psi_n^{-1}(ct)$, we consider three ranges of $k$. First, if $k\sim yU(t)$
for $y$ in some compact subset of the positive half-line, Lemma \BS\ shows
that $(ct+\g{n+k})/g_k$ is equal to
$$\displaylines{\qquad
  { ct+g\Bigl( 1-{\ds 1\over\ds (x+y)U(t)}\Bigr) 
             {\ds 1+o(1)\over\ds\Gamma(1+\gamma)}
           \over
           g\Bigl(1-{\ds 1\over\ds yU(t)}\Bigr) }
         yU(t) \Gamma(\gamma) \bigl( 1+o(1)\bigr) 
  \hfill\cr\hfill
  {}= U(t) { c\Gamma (1+\gamma)+(x+y)^\gamma \bigl(1+o(1)\bigr)
               \over
               \gamma y^{\gamma-1} } \bigl( 1+o(1)\bigr) \, ,\qquad \cr
  }
$$
with the $o(1)$ terms being uniform in $x$ and $y$ in any compact subset of 
the positive half-line. Therefore, by Lemma
\xiA.i, provided $\epsilon$ is small enough, 
the minimum of $(ct+\g{n+k})/g_k$ over $k$ in the range $\epsilon U(t)\leq
k\leq U(t)/\epsilon$ is asymptotically equivalent to 
$\xi_{c\Gamma(1+\gamma)}(x)U(t)$ as $t$ tends to infinity and any $k$ 
minimizing in this range is asymptotically equivalent 
to $U(t)\kappa_{c\Gamma(1+\gamma)}(x)$.

If $0\leq k\leq \epsilon U(t)$ and $g_k$ is positive, then for $t$ large enough
and $\epsilon$ small enough, $(ct+\g{n+k})/g_k$ is at least
$$
  {ct+\g{n}\over 2g_{\lfloor\epsilon U(t)\rfloor}}
  \sim U(t){c\Gamma(1+\gamma)+x^\gamma\over 2\gamma \epsilon^{\gamma-1}}
  \geq 2 U(t)\xi_{c\Gamma(1+\gamma)}(x) \, ,
$$
where the last inequality uses our earlier observation that when $\gamma$ is
greater than $1$, the function $\kappa_a$ is lower bounded by a positive
constant on any compact subset of the positive half-line.

Finally, if $k\geq U(t)/\epsilon$, then $(ct+\g{n+k})/g_k$ is asymptotically
bounded from below by an equivalent of
$$
  { \g{k}\over g_k }\sim {k\over\gamma} \geq {U(t)\over\epsilon\gamma} \, ,
$$
which is asymptotically greater 
than $U(t)\xi_{c\Gamma (1+\gamma)}(x)$, again provided that $\epsilon$ is
small enough.

\noindent{\it Case $\gamma=1$}. This case is more involved, and for clarity
of the argument, we split the proof into several steps. The first one consists
in proving the result for a continuous analogue of the minimization
problem involved in the variational form of $\Psi_n^{-1}$. Recall that 
$g(1-1/t)=t\ell (t)$ for some slowly varying function $\ell$ having
the Karamata representation
$$
  \ell(x)=d(x)\exp\int_1^x {\varepsilon (u)\over u} \d u \,
$$
ultimately and where $\varepsilon$ satisfies \HypGOneB\ and \HypGOneC.

\noindent{\it Step 1.} Assume that the function $d(\,\cdot\,)$ is constant.
For any fixed positive $c$, consider the function
$$
  \phi_t(y)={\ds ct+g\Bigl(1-{1\over xU(t)+y}\Bigr)\over \ell (y)} \, .
$$
Asymptotically in $t$, this function is a continuous analogue of
the function $y\mapsto (ct+\g{xU(t)+y})/g_{\lfloor y\rfloor}$. In view of
Lemma \psiInv,
we are seeking the minimum value of $\phi_t$ as well as its minimizing
argument. First, the minimizer has to tend to infinity with $t$, for 
if $y$ stays bounded then $\phi_t(y)\sim (c+x^\gamma)t/\ell(y)$ and,
since we assume that \PrelimA\ holds so that $\ell$ tends to infinity, 
this asymptotic equivalent can be made 
smaller by increasing $y$. Second, the minimizer has to be of smaller order
than $U(t)$ because if $y\sim \theta U(t)$ for some $\theta$ in a compact
subset of the positive half-line, then
$$
  \phi_t(y) \sim {t\over\ell\circ U(t)} (c+x+\theta)
$$
which attains its minimum for $\theta$ vanishing; and, moreover, if $y$ is 
of order larger than $U(t)$, then the same argument as in the case $\gamma>1$ 
show that $\phi_t(y)$ cannot be minimum.

Next, differentiating $\phi_t$ and after substitution of $\Id\ell$ for
$g(1-1/\Id)$, the minimizer satisfies
$$\displaylines{\qquad
  0 = 
  { \ell\bigl( xU(t)+y\bigr) + \bigl( xU(t)+y\bigr)\ell'\bigl( xU(t)+y\bigr)
    \over
    \ell (y) }
  \hfill\cr\hfill
  {}- {\ell'\over\ell^2} (y) \Bigl(ct+g\Bigl(1-{1\over xU(t)+y}\Bigr)\Bigr) 
  \, .
  \qquad\cr}
$$
Since the Karamata representation of $\ell$ with a constant function $d$
implies $\Id\ell'/\ell=o(1)$, and since the minimizer is $o\bigl(U(t)\bigr)$,
it follows that, after factoring $1/\ell(y)$ and simplifying,
$$
  0=\ell\bigl( xU(t)\bigr)\bigl(1+o(1)\bigr)-
  {\varepsilon(y)\over y} \Bigl(ct+g\Bigl(1-{1\over xU(t)+y}\Bigr)\Bigr) \, .
$$
Consequently, since $\gamma=1$ and 
$$
  g\Bigl( 1-{1\over xU(t)+y}\Bigr) 
  \sim g\Bigl( 1-{1\over xU(t)}\Bigr)
  \sim xt \, ,
$$
the minimizer satisfies
$$
  {\varepsilon (y)\over y} 
  \sim {\ell\circ U(t)\over t} {1\over c+x} \, .
$$
Furthermore, 
$$
  \ell\circ U\sim g(1-1/U)/U \sim \Id/U \, ,
  \eqno{\equa{PsiInvUAA}}
$$
at infinity, and therefore the minimizer satisfies
$$
  { \varepsilon (y)\over y}
  \sim {1\over U(t) (c+x)} \, .
$$
Assumption \HypGOneB\ then implies $y\sim (c+x)U(t)\varepsilon \circ U(t)$
--- compare with the value for $k_{\lfloor xU(t)\rfloor}(ct)$ given
in the statement of the lemma. For such value of $y$, we have
$$
  \phi_t(y)
  \sim { t\over\ell\bigl( U(t)\varepsilon\circ U(t)\bigr)} (c+x) \, .
$$
Since $\log\varepsilon (e^t)$ is self-neglecting and therefore self-controlled,
Theorem 3.12.5 in Bingham, Goldie and Teugels (1989) shows that 
$\ell \bigl( U\varepsilon(U)\bigr)\sim \ell(U)$ as $U$ tends to infinity.
Combined with \PsiInvUAA, this yields
$$
  \phi_t(y)
  \sim {t\over \ell\circ U(t)}(c+x)
  \sim U(t)(c+x)
$$
as $t$ tends to infinity.

\noindent{\it Step 2.} This step consists in showing that the minimum of
$\phi_t(y)$ has some form of 
continuity with respect to the asymptotic behavior of $g$. Assume now that we
have another function, $g_1$ asymptotically equivalent to $g$ at
$1-$.  With obvious notation, this new function gives rise to the
corresponding functions $U_1$ and $\ell_1$.  Then, for any positive $\eta$,
$$\eqalignno{
  {ct+g\Bigl(1-{\ds 1\over\ds xU(t)+y}\Bigr)\over \ell(y)}
  &{}\leq { ct+(1+\eta)g_1\Bigr( 1-{\ds 1\over\ds xU_1(t)+y}\Bigr) 
            \over (1-\eta)\ell_1(y) } \cr
  &{}\leq {1+\eta\over 1-\eta} \,
   { ct+g_1\Bigl(1-{\ds 1\over\ds xU_1(t)+y}\Bigr) \over \ell_1(y) }\qquad\quad
  &\equa{psiInvUA}\cr}
$$
as $t$ tends to infinity, and uniformly in the range $y$ nonnegative. 
Step 1 of this proof shows that even though we do not assume $\ell_1$ to
be smooth, the minimizer of the corresponding $\phi_{1,t}(y)$ function is
$o\bigl( U_1(t)\bigr)$. Therefore, \psiInvUA\ and the analogous lower bound 
obtained by permuting $g$ and $g_1$ show that the minimizer of $\phi_{1,t}$
is asymptotically equivalent to that of $\phi_t$; moreover, $\phi_{1,t}$ and
$\phi_t$ have asymptotically equivalent minimum values. It follows that
the conclusion of step 1 remains valid if we only assume that the 
function $d(\,\cdot\,)$
in the Karamata representation of $\ell$ has a limit and that the function
$\varepsilon(\,\cdot\,)$ in that representation satisfies \HypGOneB\ 
and \HypGOneC.

\noindent{\it Step 3.} Recall that $\gamma$ is $1$ here, so that 
both $\Gamma(\gamma)$ and $\Gamma(1+\gamma)$ are $1$ as well. Going back to 
the problem of 
evaluating $\Psi_{\lfloor xU(t)\rfloor}^{-1}(ct)$, we have $g_k\sim
\ell(k)$ and, by Karamata's theorem, $\g{k} \sim k\ell (k)$, this equivalent
being uniform in the range of $k$ of order $U(t)$.  This allows us to replace 
the discrete minimization
to calculate $\Psi_{\lfloor xU(t)\rfloor}^{-1}(ct)$ by the continuous one
solved in the first step. Using that $\xi_a(x)=a+x$ when $\gamma$ is $1$,
this proves the lemma.\hfill\qed

\bigskip

In the preceding lemma, writing $n$ for $\lfloor xU(t)\rfloor$, we obtain
$$
  \Psi_n^{-1}(t) \sim U(t)\xi_{\Gamma (1+\gamma)} \bigl(n/U(t)\bigr)
  \eqno{\equa{psiInvAsymp}}
$$
as $t$ tends to infinity. Lemma \xiAsymp\ asserts that if $x$ is large,
then $\xi_{\Gamma (1+\gamma)}(x)$ is about $x\xi_0(1)$. Consequently,
we expect that if $n/U(t)$ is large then $\Psi_n^{-1}$ is about $n\xi_0(1)$.
The following bounds show that in some sense this is indeed the case.

\Lemma{\label{psiInvBounds}%
  Let $\epsilon$ be a positive real number less than $1$.

  \noindent (i) For $n$ large enough, for any positive $t$,
  $$
    \Psi_n^{-1}(t)\geq (1-\epsilon) n\xi_0(1) \, .
  $$

  \noindent (ii) For $n$ at least $U(t)/\epsilon$ and $t$ large enough,
  $$
    \Psi_n^{-1}(t) \leq (1+\epsilon) n\xi_0(1) \, .
  $$
}

\Proof (i) Clearly, $\Psi_n^{-1}(t)$ is at least $\min_{k\in K_+}\g{n+k}/g_k$.
Setting $k=yn$ with $y$ in a compact set of the positive half-line, we obtain,
as $n$ tends to infinity,
$$
  {\g{n+k}\over g_k} 
  \sim {1\over \gamma} 
   { g\Bigl(1-{\ds 1\over\ds n(1+y)}\Bigr) 
     \over 
     g\Bigl(1-{\ds 1\over\ds ny}\Bigr)} ny
  \sim {n\over\gamma} {(1+y)^\gamma\over y^{\gamma-1}}
  \geq n\xi_0(1) \, .
$$
Next, let $\delta$ be a positive real number. For any positive $k$ at most
$\delta n$ and for any $g_k$ positive, Lemma \BS\ yields
$$
  {\g{n+k}\over g_k}
  \gsim {\g{n}\over g_{\lfloor\delta n\rfloor}}
  \sim {n\over\gamma} {1\over\delta^{\gamma-1}} \, ,
$$
while for any $k$ at least $n/\delta$, it yields
$$
  {\g{n+k}\over g_k}
  \gsim {\g{k}\over g_k}
  \gsim {n\over\gamma} {1\over\delta} \, . 
$$
The result follows by taking $\delta$ such that 
$\delta^{1-\gamma}\wedge \delta^{-1}\geq\xi_0(1)$. Note that when $\gamma$ is
$1$, the result still holds because then $\xi_0(1)=1$.

\noindent (ii) Recall that $g(1-1/\cdot)$ is regularly varying with 
nonvanishing index. Hence, it is asymptotically equivalent to a monotone
function (Bingham, Goldie and Teugels, 1989, Theorem 1.5.3).
Since $U(t)\leq \epsilon n$, for $t$ large enough,
$$
  t
  \sim  g\Bigl( 1-{1\over U(t)}\Bigr)
  \lsim g\Bigl( 1-{1\over \epsilon n}\Bigr)
  \sim g\Bigl( 1-{1\over n}\Bigr)\epsilon^\gamma \, ,
$$
as $n$ tends to infinity. Therefore, taking $k=yn$ with $y$ fixed, when
$\gamma$ is greater than $1$, we obtain that for $t$ large enough, 
$\Psi_n^{-1}(t)$ is at most
$$\displaylines{
  {   g\Bigl(1-{\ds 1\over\ds n}\Bigr)\epsilon^\gamma +
      g\Bigl( 1-{\ds 1\over\ds n(1+y)}\Bigr)
      { \ds \bigl( 1+o(1)\bigr) \over \ds\Gamma(1+\gamma)}
    \over
      g\Bigl( 1-{\ds 1\over\ds ny}\Bigr) 
  }
  ny\Gamma(\gamma) \bigl( 1+o(1)\bigr) 
  \hfill\cr\hfill
  {}\sim n 
  {\Gamma (1+\gamma)\epsilon^\gamma+(1+y)^\gamma
   \over \gamma y^{\gamma-1}}
   \, . \cr
  }
$$
If $\gamma$ is greater than $1$, the result follows by taking $y$ minimizing 
$(1+y)^\gamma /y^{\gamma-1}$, upon noting that the minimizing value is 
positive. If $\gamma$ is $1$, the result follows by
taking $y=\epsilon$.\hfill\qed

\bigskip

While the previous lemma gives valuable information on $\Psi_n^{-1}$
when $n$ is large, it does not give any estimate for $\Psi_n^{-1}$ when
$n$ is moderate, say, and $t$ is large. The next result fills this gap.
It should be compared to \psiInvAsymp.

\Lemma{\label{psiInvLB}%
  There exists $n_1$ such that 
  $$
    \min_{n\geq n_1}\Psi_n^{-1}(t)/U(t)
    \gsim \xi_{\Gamma (1+\gamma)}(0) \, .
  $$
  as $t$ tends to infinity. Moreover, if all the $g_n$ are nonnegative, we
  can take $n_1$ to be $0$.
}

\bigskip

\Proof We first prove the following claim.

\bigskip

\Claim{
  There exists $n_1$ such that for any $n$ at least $n_1$ and any
  $k$ nonnegative, $\g{n+k}\geq \g{k}$.
}

\bigskip

\noindent Indeed, for any fixed $k$, define
$$
  m_k=\min\{\, n\geq 1\, :\, \forall i\geq n\, ,\, g_{[k,i+k)}>0\,\} \, .
$$
This integer exists since the partial sums $\g{n}$ diverge as $n$ tends to 
infinity. Recall that $n_0$ is the smallest integer such that both $g_n$ and 
$\g{n}$ are positive if $n$ is at least $n_0$. If $k$ is at least $n_0$, 
then $m_k=1$. Therefore the sequence $(m_k)_{k\geq 1}$ is bounded and admits
a maximum element which we call $\widetilde n_1$. We then set $n_1$ to be the 
maximum of $n_0$ and $\widetilde n_1$.

\bigskip

Having proved the claim, for $n$ at least $n_1$,
$$
  \Psi_n^{-1}(t) \geq \min_{k\in K_+} {t+\g{k}\over g_k} \, .
$$
If $k=\lfloor yU(t)\rfloor$ for some $y$ in a compact subset of the positive 
half-line, then
$$
  {t+\g{k}\over g_k}
  \sim U(t) {\Gamma (1+\gamma)+y^\gamma\over \gamma y^{\gamma-1}}
  \geq U(t)\xi_{\Gamma (1+\gamma)}(0)
$$
as $t$ tends to infinity.

If $k$ is at most $\delta U(t)$, then
$$
  {t+\g{k}\over g_k}
  \gsim {t\over g_{\lfloor \delta U(t)\rfloor}}
  \sim U(t)\Gamma (\gamma) \delta^{1-\gamma} \, ,
$$
while if $k$ is at least $U(t)/\delta$, then
$$
  {t+\g{k}\over g_k}
  \geq {\g{k}\over g_k }
  \sim {k\over\gamma}
  \geq {U(t)\over\gamma\delta} \, .
$$
The result follows by choosing $\delta$ small enough so 
that $\Gamma (1+\gamma)\delta^{1-\gamma}
\wedge \delta^{-1} \geq\gamma\xi_{\Gamma (1+\gamma)}(0)$ which as noted
previously holds trivially when $\gamma$ is $1$.\hfill\qed

\bigskip


\subsection{Approximation of \poorBold{$\Psi_n^{-1}$} 
when \poorBold{$(g_n)_{n\geq 0}$} has a positive and finite limit}%
In this subsection we consider the case where the sequence $(g_n)_{n\geq 0}$
has a positive and finite limit. The
minimization involved in Lemma \psiInv\ may or may not be made explicit,
according to whether the supremum of the sequence is achieved or not. To
be more precise, recall that $g^*=\sup_{n\geq 0} g_n$, and, if it exists
let $k^*$ be the smallest integer such that $g_{k^*}=g^*$. Note that
such number does not exists for a sequence such as 
$\bigl(1-(n+1)^{-1}\bigr)_{n\geq 0}$. 

\Lemma{\label{PsiInvFormulaO}%
  If $g^*$ is attained, that is $k^*$ is well defined, then there exists a
  nonnegative real number $t_0$ such that for any $n$ at least $n_0$
  and any $t$ at least $t_0$,
  $$
    \Psi_n^{-1}(t)={t+\g{n+k^*}\over g^*} \, .
  $$
  Otherwise, for any $n$ at least $n_0$, and any positive $t$,
  $$
    \Psi_n^{-1}(t)\geq {t+\g{n}\over g^*} \, ,
  $$
  and for any positive $\epsilon$ there exists $k$ such that for any positive
  $t$,
  $$
    \Psi_n^{-1}(t)\leq (1+\epsilon){t+\g{n}\over g^*} 
    + (1+\epsilon){g_{[n,n+k)}\over g^*} \, .
  $$
}

\Proof If $k^*$ exists, the proof of Lemma \PsiInvFormulag\ is still valid
and yields the result. Hence, we assume that all the $g_n$ are less than
their limit $g^*$. For $n$ at least $n_0$, the inequality $\g{n+k}\geq \g{n}$
holds and, since all the $g_n$ are less than $g^*$, the formula for
$\Psi_n^{-1}$ in Lemma \psiInv\ implies the given lower bound 
for $\Psi_n^{-1}$.

To prove the upper bound, let $k$ be any integer such 
that $g_k\geq g^*/(1+\epsilon)$. Then, the formula for $\Psi_n^{-1}$ in Lemma
\psiInv\ shows that
$$
  \Psi_n^{-1}(t)\leq (1+\epsilon) {t+\g{n+k}\over g^*} \, ,
$$
and the result follows by writing $\g{n+k}$ as $\g{n}+g_{[n,n+k)}$.\hfill\qed

\bigskip

We then obtain the following asymptotic equivalence.

\Lemma{\label{PsiInvUO}%
  Uniformly in $x$ and $c$ in any compact subset of the positive half-line,
  $$
    \Psi_{\lfloor xU(t)\rfloor}^{-1}(t)\sim {t\over g^*} (1+x)
  $$
  as $t$ tends to infinity.
}

\bigskip

\Proof The lemma follows from Lemma \PsiInvFormulaO.\hfill\qed

\bigskip

Our next lemma is stated so that it can be easily referred to. It involves
a real number $t_0$ defined in Lemma \PsiInvFormulaO\ when $k^*$ exists, and
otherwise, one can take $t_0$ to be $1$.

\Lemma{\label{PsiInvLbO}%
  For any $n$ at least $n_0$ and any $t$ at least $t_0$,
  $$
    \Psi_n^{-1}(t)\geq t/ g^*
  $$
  and, as $n$ tends to infinity,
  $$
    \inf_{t\geq t_0}\Psi_n^{-1}(t) \gsim {g(1-1/n)\over g^*} \, . 
  $$
}

\Proof The proof is the same as that of Lemmas \PsiInvLbgA\ and 
\PsiInvLbgB.\hfill\qed

\bigskip


\def\preveq{(\the\sectionnumber .\the\equanumber)}
\def\prevs{\the\sectionnumber .\the\snumber }

\section{Proof of the results of section \fixedref{2}}%
The proof follows by an application of Theorem \UsefulStatement, after
completion of all the steps described in section \fixedref{3.1}. While most
arguments depend on the asymptotic behavior of $(g_n)_{n\geq 0}$ at infinity, 
the calculation of the conditional expectation, encoded in the 
functions $\psi_{i,n}$ and $\psi_n$, can be done once and for all. Indeed,
without any loss of generality, we assume that $EX_i=-1$. Then
$$
  \psi_{i,n}(x)
  =E(S_i|X_n=x)
  =\cases{ -\g{i} & if $i<n$, \cr
           (x+1)g_{i-n}-\g{i} & if $i\geq n$.\cr }
  \eqno{\equa{psiin}}
$$
It follows that
$$\eqalign{
  \psi_n(x)
  &{}=\max_{i\geq 0} \psi_{i,n}(x)\cr
  &{}=\max_{0\leq i<n} (-\g{i})\vee
   \max_{i\geq n}\bigl( (x+1)g_{i-n}-\g{i}\bigr)\, .\cr
  }
$$
Therefore, if $x$ is such that $\psi_n(x)$ is both positive and greater than
$\max_{i\geq 0} -\g{i}$, then 
$$
  \psi_n(x)=0\vee\max_{i\geq n} (x+1)g_{i-n}-\g{i} \, .
$$
Writing $t_1$ for $0\vee\max_{i\geq 0} -\g{i}$, 
this shows that, with the notation of section \fixedref{5.2}, for any
$n$ positive and any $x$ in the preimage under $\psi_n$ 
of $(t_1,\infty)$,
$$
  \psi_n(x)=\Psi_n(x+1) \, .
$$
Therefore, for $t$ greater than $t_1$,
$$
  \psi_n^{-1}(t)=\Psi_n^{-1}(t)-1\, .
$$

Define the sequence of independent and identically distributed centered
random variables $(Z_i)_{i\geq 1}$ by $Z_i=X_i-\mu$, $i\geq 1$.
A useful remark for completing step 6 is that
$$
  S_i-\psi_{i,n}(X_n)=\sum_{\ss 0\leq j<i\hfill\atop\ss i-j\not= n\hfill}
  g_j Z_{i-j} \, .
$$

Finally we will use also the following weak law of large number.

\Lemma{\label{LLN}%
  Let $(S_n)$ be a $(g,F)$-process with negative mean innovations. Then,
  $$
    \lim_{n\to\infty} {S_n\over g(1-1/n)} = {\mu\over \Gamma(\gamma+1)}
  $$ 
  in probability; in other words, $S_n/ES_n$ converges to $1$ in probability.
}

\bigskip

\Proof Lemma \BS.i shows that only a finite number of $g_n$ may be
nonpositive.  Therefore, here, we can assume without any loss of
generality that all the coefficients $g_n$ are positive. Lemma \BS\
and the uniform convergence theorem for regularly varying functions
(Bingham, Goldie and Teugels, 1989, Theorem 1.2.1) imply that
uniformly in $\alpha$ in any compact subset of $(0,1\,]$,
$$
  \lim_{n\to\infty} g_{\lfloor\alpha n\rfloor}/ \g{n}
  = 0 \, .
$$
Therefore, Theorem 3 in Jamison, Orey and Pruitt (1965) implies
$$
  \lim_{n\to\infty} S_n /\g{n} = \mu
$$
in probability. The result then follows from Lemma \BS.ii.\hfill\qed

\bigskip

\noindent{\bf Notation.} In the remainder of this paper we write $s_n$ for 
the expected value of $S_n$, that is, $s_n=\mu\g{n}$.

\bigskip


\def\preveq{(\the\sectionnumber .\the\subsectionnumber .\the\equanumber)}
\def\prevs{\the\sectionnumber .\the\subsectionnumber .\the\snumber }

\subsection{Proof of Theorem \Thg\ -- upper bound}%
We complete all the steps described in section \fixedref{3}. Recall that
without any loss of generality, we assume that $\mu$ is $-1$. Also, replacing
$t$ by $t/(-\mu g^*)$, and $g$ by $g/g^*$, we assume without loss of
generality that $g^*=1$.

\bigskip

\noindent {\it Step 1.} Since $\psi_n^{-1}=\Psi_n^{-1}-1$ on 
$(t_1,\infty)$, 
relation (S1) follows from Lemma \PsiInvUg\ with $\chi$ being $\Id$ and
$\rho(x)=1+x^\gamma/\Gamma(1+\gamma)$.

\bigskip

\noindent {\it Step 2.} Recall the asymptotic equivalence given 
in \rcalculation. Lemma \PsiInvLbgA\ and regular variation of $\oF$ 
and $U$ imply that as $t$ tends to infinity,
$$\eqalign{
  \sum_{n\leq \epsilon U(t)} \oF\circ \psi_n^{-1} (t)
  &{}\lsim \sum_{n\leq \epsilon U(t)} \oF(t) \cr
  &{}\sim \epsilon (U\oF)(t)\, , \cr
  }
$$
while Lemma \PsiInvLbgB\ implies
$$
  \sum_{n\geq U(t)/\epsilon} \oF\circ\psi_n^{-1} (t)
  \lsim \sum_{n\geq U(t)/\epsilon} 
  \oF\Bigl( {g(1-1/n)\over \Gamma(1+\gamma)}\Bigr) \, .
$$
Using the regular variation of $\oF$ and $g$, the approximation of the sum
in this upper bound by a Riemann integral and Karamata's theorem, we obtain
an asymptotic upper bound equivalent to
$$
  \Gamma(1+\gamma)^\alpha
  \int_{U(t)/\epsilon}^\infty \oF\circ g(1-1/u)\d u
  \sim \Gamma(1+\gamma)^\alpha 
  {\epsilon^{\alpha\gamma-1}\over\alpha\gamma-1} (U\oF)(t) \, .
$$
Since $\alpha\gamma$ is larger than $1$, this completes the proof of (S2).

\bigskip

\noindent{\it Step 3.} We prove that the maximum of the process 
is unlikely to occur at a time of smaller order than $U(t)$.

\Lemma{\label{SmallGStepThree}%
  The following limit holds,
  $$
    \limeps\limsupt 
    { P\{\, \exists n \, :\, n\leq\epsilon U(t)\, , \, S_n>t\,\}
      \over
      (U\oF)(t) } = 0 \, .
  $$
}

\Proof Let $n$ be at least $n_0$. Recall that $s_n$ is the expectation 
of $S_n$.Since $t-s_n$ is at least $t$, if $S_n$ exceeds $t$ 
then $\sum_{0\leq j<n} g_j Z_{n-j}$ exceeds $t$ as well. By the standard
estimate for weighted convolution of distribution functions with regularly
varying tails and Bonferroni's inequality, for any fixed $k$ and any positive
$\delta$,
$$
  \limeps\limsupt { P\{\, \exists n\,:\, n\leq\epsilon U(t)\, ,\, 
  \sum_{0\leq j<k} g_jZ_{n-j} >\delta t\,\}\over (U\oF)(t) } = 0 \, .
$$
Therefore, taking $\delta$ to be less than $1/2$, it suffices to prove that
for some $k$,
$$
  \limeps\limsupt { P\{\, \exists n\,:\, k<n\leq\epsilon U(t)\, ,\, 
  \sum_{k\leq j<n} g_jZ_{n-j} >\delta t\,\}\over (U\oF)(t) } = 0 \, .
$$
Since
$$
  P\{\, \exists j \, :\, j\leq\epsilon U(t)\, ,\, Z_j>t\,\}
  \lsim \epsilon (U\oF)(t) \, ,
$$
it suffices to prove that
$$\displaylines{
  \limeps\limsupt \epsilon U(t)\max_{k\leq n\leq\epsilon U(t)} 
  { P\{\, \sum_{k\leq j<n} g_jZ_{n-j}\II\{\,Z_{n-j}\leq t\,\} >\delta t\,\}
    \over (U\oF)(t) } 
  \hfill\cr\hfill
  {}= 0 \, .\qquad
  \equa{SmallGStepThreeA}\cr}
$$
Using that $-\log\oF\sim \alpha\log$ at infinity, and taking $\lambda$ of
the form $ct^{-1}\log t$ in inequality \LargeDev\ with $\delta t$ in place 
of $t-s_n$ there, the logarithm of the ratio
$$
  P\Bigl\{\, \sum_{k\leq j<n} g_j Z_{n-j}\II\{\, Z_{n-j}\leq t\,\} >\delta t\,
  \Bigr\} \bigm/\oF(t)\, ,
$$
for $n$ at most $\epsilon U(t)$, is ultimately at most
$$\displaylines{\quad
  -c\delta \log t +\eta c t^{-1}\log t \,\g{n}E|Z|
  \hfill\equa{SmallGStepThreeB}\cr\noalign{\vskip 6pt}\hfill
  \eqalign{
  & {}- ct^{-1}\log t \sum_{0\leq j<n} g_j EZ
    \II\{\, \lambda g_jZ\leq\log(1-\eta)\,\} \cr
  & {}+\epsilon U(t)\exp \Bigl( c\log t\max_{k\leq j<U(t)} g_j\Bigr)\oH
    \Bigl( {t\log(1+\eta)\over c\log t}\Bigr)
    + 2\alpha \log t \, .\cr
  }
  \quad\cr}
$$
Referring to the second summand in this bound, $\g{n}$ is at most 
$\g{\epsilon U(t)}\sim\epsilon^\gamma t/\Gamma(1+\gamma)$. Thus, 
the second summand is 
ultimately at most $2\eta c\epsilon^\gamma E|Z|\log t$. Since 
$\lambda$ tends to $0$ as $t$ tends to infinity and so
$EZ\II\{\, \lambda Z\leq \log (1-\eta)\,\}=o(1)$, the third summand
is negligible compared to the second one. Therefore, ultimately, 
\SmallGStepThreeB\ is at most
$$\displaylines{
  \bigl( -c\delta+ 2\eta c\epsilon^\gamma E|Z|+2\alpha\bigr)\log t
  \hfill\cr\hfill
  {}+\epsilon \exp\Bigl( \log U(t)+c\max_{j\geq k} g_j\log t
  + \log \oH\Bigl( {t\log (1+\eta)\over c\log t}\Bigr)\Bigr)\, .
  \quad\equa{SmallGStepThreeC}\cr
  }
$$
We take $\eta$ small enough so that $-\delta+2\eta\epsilon^\gamma E|Z|$ is 
negative. Then,
we take $c$ large enough so that whenever $\epsilon$ is less than $1$,
$$
  c(-\delta +2\eta\epsilon^\gamma E|Z|)+2\alpha <-2 \, ,
$$
say. Since $\alpha\gamma>1$ and $(g_n)_{n\geq 0}$ converges to $0$, we 
can then fix $k$ large enough so that
$$
  \gamma^{-1}+c\max_{j\geq k} g_j-\alpha <0 \, .
$$
For such $k$, we obtain that
$$
  \limt \log U(t) + c\max_{j\geq k} g_j \log t 
  + \log\oH\Bigl({t\log (1+\eta)\over c \log t}\Bigr)
 =-\infty \, .
$$
Hence, ultimately, \SmallGStepThreeC\ is at most $-\log t$. It follows that
\SmallGStepThreeA\ holds as well as the conclusion of the lemma.\hfill\qed

\bigskip

\noindent {\it Step 4.} We prove that the maximum of the process is unlikely
to occur at a time of larger order than $U(t)$.

\Lemma{\label{SmallGStepFour}%
  The following limit holds,
  $$
    \limeps\limsupt 
    { P\{\, \exists n\, : \, n\geq U(t)/\epsilon \, , \,S_n>t\,\}
      \over (U\oF)(t) } = 0 \, .
  $$
}

\Proof Referring to step 4 in section \fixedref{3}, it suffices to show that
\StepFourAlt\ holds. In the current context, 
$\chi\circ U^{-1}\sim g(1-1/\Id)$ at infinity. Thus, given two positive 
real number $\theta$ and $\eta$, set
$$
  a_n=\theta \Bigl( g\Bigl( 1-{1\over n}\Bigr) + {t\over\eta}\Bigr) \, .
$$ 
If all $X_n$ are at most $a_n+\mu$, then
$$
  S_n-s_n=\sum_{0\leq i<n} g_iZ_{n-i} \II\{\, Z_{n-i}\leq a_{n-i}\,\} \, .
$$
Hence, the probability involved in \StepFourAlt\ is at most
$$
  \sum_{n\geq U(t)/\epsilon} P\Bigl\{\,\sum_{0\leq i<n} g_iZ_{n-i}
  \II\{\, Z_{n-i}\leq a_{n-i}\,\} > t-s_n\,\Bigr\} \, .
$$
The usual estimate for the tail of weighted convolutions of heavy-tail 
distribution functions shows that for any fixed $k$,
$$\displaylines{\qquad
  \sum_{n\geq U(t)/\epsilon} P\Bigl\{\, \sum_{0\leq j<k} g_jZ_{n-j}
  \II\{\, Z_{n-j}\leq a_{n-j}\,\}>{t-s_n\over 2} \,\Big\}
  \hfill\cr\hfill
  \eqalign{
  {}\lsim{}& k2^\alpha \sum_{n\geq U(t)/\epsilon} \oF(t-s_n) \cr
  {}\sim{} & k2^\alpha\int_{U(t)/\epsilon}^\infty 
             \oF\Bigl(t+{g(1-1/u)\over\Gamma(1+\gamma)}\Bigr) \d u \cr}
  \qquad\cr}
$$
which, by the proof of Lemma \Karam, is negligible compared to $(U\oF)(t)$
as first $t$ tends to infinity and then $\epsilon$ tends to $0$. Therefore, it
suffices to prove that for some fixed $k$,
$$\displaylines{
  \limeps\limsupt \sum_{n\geq U(t)/\epsilon} 
  { P\Bigl\{ \sum_{k\leq i<n} g_iZ_{n-i}\II\{\, Z_{n-i}\leq a_{n-i}\,\}>
    {\ds t-s_n\over\ds 2}
    \,\Bigr\} \over (U\oF)(t) } 
  \hfill\cr\hfill
  {}= 0 \, . \qquad
  \equa{SmallGStepFourB}
  \cr}
$$
Using inequality \LargeDev, the logarithm of each summand is at most 
asymptotically bounded by an equivalent of
$$\displaylines{\quad
    -\lambda {t-s_n\over 2} +\eta \lambda \g{n}E|Z| 
    -\lambda \sum_{0\leq j<n} g_j EZ\II\{\, \lambda g_j Z\leq\log (1-\eta)\,\}
  \hfill\cr\hfill
    {}+\sum_{k\leq j<n} e^{\lambda g_ja_{n-j}}
    \oH\Bigl( {\log(1+\eta)\over\lambda g_j}\Bigr) 
    - \Bigl({1\over\gamma}-\alpha\Bigr) \log t \bigl( 1+o(1)\bigr) \, .
    \qquad \equa{SmallGStepFourB}\cr
  }
$$
We take $\lambda$ of the form $2c(t-s_n)^{-1}\log g(1-1/n)$ for a 
constant $c$ to be determined later.

Referring to the successive terms in \SmallGStepFourB, we have
$$
  \lambda {t-s_n\over 2}=c\log g(1-1/n) \sim c\gamma\log n
$$
as $n$ tends to infinity. Furthermore, since $t-s_n\geq -s_n=\g{n}$,
$$\eqalign{
  \eta\lambda\g{n} E|Z|
  &{}\sim 2\eta c{\log g(1-1/n)\over t-s_n} \g{n} E|Z| \cr
  &{} \lsim 2 \eta c\gamma\log n E|Z| \, .\cr
  }
$$
Next, since $\lambda$ tends to $0$ as $t$ tends to infinity and
uniformly in $n\geq U(t)/\epsilon$,
$$
  \lambda \sum_{0\leq j<n} g_j EZ\II\{\, \lambda g_j Z\leq \log (1-\eta)\,\}
  = \lambda o(\g{n})
  = o(\log n) \, . 
$$
Moreover, since $n\geq U(t)/\epsilon$ and therefore, 
$g(1-1/n)\gsim \epsilon^{-\gamma}t$, we have $t\lsim -\epsilon^\gamma
\Gamma(1+\gamma) s_n$, and since $\gamma$ is at most $1$,
$$\eqalign{
  \lambda g_j a_{n-j}
  &{}\lsim {2c\log g(1-1/n)\over t-s_n} \theta
    \Bigl( g\Bigl( 1-{1\over n}\Bigr)+{t\over\eta}\Bigr) \cr
  &{}\lsim 4c\gamma \theta \Gamma(1+\gamma)
   \Bigl(1+{\epsilon^\gamma\over\eta}\Bigr) \log n \, . \cr
  }
$$
Finally, for the same reason, $\log t\lsim \log g(1-1/n)\sim\gamma\log n$.
In particular, this implies
$$\eqalign{
  \log \oH\Bigl( {\log(1+\eta)\over \lambda g_j}\Bigr)
  \lsim \alpha\log\lambda
  \sim -\alpha\log (t-s_n)
  &{}\lsim -\alpha\log (-s_n)\cr
  &{}\sim -\alpha\gamma\log n \, .\cr
  }
$$
Therefore,
$$
  \sum_{k\leq j<n} e^{\lambda g_ja_{n-j}}
  \oH\Bigl( {\log (1+\eta)\over\lambda g_j}\Bigr) 
  \lsim n^{1-\alpha\gamma+(4c\gamma\theta\Gamma(1+\gamma)) +o(1)} 
  \, .
$$
It follows that \SmallGStepFourB\ is asymptotically bounded by an equivalent
of
$$
  \Bigl( -c\gamma +2\eta c\gamma E|Z| +o(1) + \alpha-{1\over\gamma}\Bigr)\log n
  + n^{1-\alpha\gamma+ 4c\gamma \theta\Gamma(1+\gamma)+o(1)} \, .
  \eqno{\equa{SmallGStepFourC}}
$$
We take $\eta$ less than $1/6E|Z|$. We take $c$ large enough so that
$$
  -{c\gamma\over 3}+\alpha-{1\over \gamma} <-3 \, .
$$
Since $\alpha\gamma$ is larger than $1$, we can take $\theta$ small 
enough so that 
$$
  1-\alpha\gamma + 8c\gamma\theta\Gamma(1+\gamma) < 0 \, .
$$
These choices lead to that \SmallGStepFourC\ as well as \SmallGStepFourB\ are 
asymptotically bounded by  $-3\log n$. Since
$$
  \sum_{n\geq U(t)/\epsilon} n^{-2}\sim {\epsilon\over U(t)} \, ,
$$
This proves the lemma\hfill\qed

\bigskip

\noindent{\it Step 5.} We now prove that for $M$ to exceed $t$ it is 
likely that we must have at least one random variable $X_n$ to be large.
Recall that $B_j$ is the event $\{\, X_j\leq\theta t\,\}$. For any $t$ large
enough, up to increasing $\theta$ slightly, we can replace $X_j$ by $Z_j$
in the definition of $B_j$, and so we set $B_i=\{\, Z_i\leq\theta t\,\}$.

\Lemma{\label{SmallGStepFive}
  For any positive $\epsilon$, there exists  positive $\theta$
  such that
  $$
    P\Bigl( \exists n\in I \,:\, S_n>t\,;\, \bcap_{0\leq i\leq n} B_i\Bigr)
    = o(U\oF)(t)
  $$
  as $t$ tends to infinity.
}

\bigskip

\Proof Recall that $N$ denotes $\lfloor U(t)/\epsilon\rfloor$. On 
$\bcap_{1\leq i\leq n}B_i$ we have
$$
  \sum_{0\leq i<n} g_iZ_{n-i}
  =\sum_{0\leq i<n} g_i Z_{n-i}\II\{\, Z_{n-i}\leq\theta t\,\} \, .
$$
We apply inequality \LargeDev\ and use that $n$ is at most $N$ to obtain that
the logarithm of $P\{\, S_n>t\,;\, \bcap_{1\leq i\leq n}B_i\,\}$ is at most
$$\displaylines{\qquad
  -\lambda t +\eta\lambda\g{N}E|Z|
  -\lambda\sum_{0\leq j<n} g_j EZ\II\{\, \lambda g_jZ\leq\log (1-\eta)\,\}
  \hfill\cr\hfill
  {}+\sum_{0\leq j<n} e^{\lambda g_j\theta t}
  \oH\Bigl( {\log (1+\eta)\over\lambda g_j}\Bigr) \, .
  \qquad
  \equa{SmallGStepSixA}
  }
$$
We choose $\lambda$ of the form $ct^{-1}\log t$ where $c$ will be specified
later. With this choice, we examine all the terms in \SmallGStepSixA. We have
$\lambda t=c\log t$. Lemma \BS\ implies
$$
  \lambda\g{N}
  \sim c{\log t\over t} {g\bigl(1-\epsilon/U(t)\bigr)\over \Gamma(1+\gamma)}
  \sim {c\epsilon^{-\gamma}\over\Gamma(1+\gamma)}\log t \, .
$$
Since $\lambda$ tends to $0$, we also have, referring to the third
summand in \SmallGStepSixA,
$$
  -\lambda \sum_{0\leq j<n} g_jEZ\II\{\, \lambda g_j Z\leq \log (1-\eta)\,\}
  \leq\lambda\g{n}o(1)
$$
as $t$ tends to infinity. Finally, the bound
$$
  \lambda g_j\theta t
  \leq c\theta\log t
$$
shows that
$$
  \sum_{0\leq j<n} e^{\lambda g_j\theta t} 
  \oH\Bigl( {\log (1+\eta)\over\lambda g_j}\Bigr)
  \lsim {U(t)\over\epsilon} t^{c \theta} 
  \oF\Bigl( {t\over\log t}\Bigr)c^\alpha\log^{-\alpha} (1+\eta) \, .
  \eqno{\equa{SmallGStepSixB}}
$$
Since $\alpha\gamma>1$, Potter's bounds imply that if
$$
  {1\over\gamma}+c\theta-\alpha <0
  \eqno{\equa{SmallGStepSixC}}
$$
then \SmallGStepSixB\ tends to $0$ as $t$ tends to infinity. Since this is
the case by choosing $c$ and $\theta$ such that $c\theta$ is sufficiently 
small, \SmallGStepSixA\ is bounded by an asymptotic equivalent of
$$
  \Bigl( -c+{\eta c\epsilon^{-\gamma}E|Z|\over\Gamma (1+\gamma)}\Bigr)
  \log t \, .
$$
Let $p$ be any positive number.
We take $\eta$ such that $\eta\epsilon^{-\gamma}/\Gamma (1+\gamma)<1/2$
say. Then, we take $c$ large enough so that $-c/2 <-p-2/\gamma$. Then,
we choose $\theta$ small enough so that \SmallGStepSixC\ holds. This shows
that 
$$
  \max_{n\in I} P\{\, S_n>t\,;\, \bcap_{1\leq i\leq n}B_i\,\}
  \leq t^{-p-2/\gamma}
$$
ultimately in $t$. Therefore,
$$
  P\Bigl\{\, \exists n\in I\,:\, S_n>t\,;\,\bcap_{1\leq i\leq n}B_i\,\Bigr\}
  \leq \epsilon^{-1}U(t) t^{-p-2/\gamma} \, ,
$$
which, by Potter's bounds is $o(t^{-p})$ as $t$ tends to infinity. Taking
$p$ greater than $\alpha$ proves the lemma.\hfill\qed

\bigskip

\noindent{\it Step 6.} We prove the law of large number 
which allows one to approximate all the $S_i$, $0\leq i\leq U(t)/\epsilon$, 
given that $X_n$ is large. Recall that the sets $B_n$ and $C_n$ involved
in step 6 depend on the parameter $\theta$, while $D_n$ depends on a
parameter $\delta$ and the set $I$ depends on a parameter $\epsilon$.
In the following lemma, $D_n$ refers in fact to the one sided event
$$
  D_n=\bcap_{1\leq i\leq N}\{\, S_i-\psi_{i,n}(X_n)\leq \delta t\,\} \, .
$$

\Lemma{\label{SmallGStepSix}%
  Let $\delta$ be a positive real number. For any $\epsilon$ positive small
  enough, there exists a positive $\theta$ such
  that
  $$
    P\Bigl( \bcup_{n\in I}B_n^\complement \cap C_n\cap D_n^\complement\Bigr) 
    = o\bigl(r(t)\bigr)
  $$
  as $t$ tends to infinity.
}

\bigskip

\Proof The result follows from the same estimate as in the previous lemma, 
taking $p$ to be large enough in that proof upon using that
$S_i-\psi_{i,n}(X_n)=\sum_{\ss 0\leq j<i\hfill\atop\ss i-j\not= n\hfill}
g_j Z_{i-j}$.\hfill\qed

\bigskip

Having completed steps 1--5 and the one-sided version of step 6, 
Theorem \UsefulStatement\ yields the upper bound pertaining to Theorem \Thg.

\bigskip


\subsection{Proof of Theorem \Thg\ -- lower bound}%
To prove a lower bound matching the upper bound, we could use a tail balance
condition and a two-sided version of the events $B_n$ to check 
the two-sided version of (S6) --- which would then
follow from the proof of Lemma \SmallGStepSix --- and then
check that (S7) holds. To verify (S7) is particularly easy
because \psiin\ and Lemma \PsiInvUg\ show that $i(n,x)=n+k^*$. 
However, in order not to impose a
tail balance condition, we give a proof inspired by Zachary's (2004) 
probabilistic proof of Veraverbeke's theorem. Zachary's proof, suitably 
modified, is remarkably robust to the choice of the process.

For this proof, we keep assuming, without any loss of generality
that the mean $\mu$ of $X_i$ is $-1$ and that
$g^*$ is $1$. Recall that  $k^*$ is the smallest integer $k$ such 
that $g_k$ is
maximal, that is equal to $g^*$. Recall we are in the case where 
the sequence $(g_i)_{i\geq 0}$ is nonnegative and tends to $0$ at infinity, 
so that the sequence attains its maximum value, assumed to be positive, and,
by our convention, $g_{k^*}=1$.
Let $\widehat S_n=\sum_{i\geq 0; i\not= k^*} g_iX_{n-i}$.

We write $s_n$ for the expectation of $S_n$ and $\hat s_n$ for that of
$\widehat S_n$. Lemma \BS.ii shows that $s_n\sim \hat s_n$ as $n$
tends to infinity. Moreover, Lemma \LLN\ implies that for any positive
$\epsilon$ and any $n$ larger than some $n_2$,
$$
  P\{\, \widehat S_n>(1+\epsilon)s_n\,\} \geq 1-\epsilon \, .
$$

If the event
$$
  \{\, \widehat S_n>(1+\epsilon) s_n\, ;\, 
       X_{n-k^*}>t-(1+\epsilon)s_n\,\}
$$
occurs, then $S_n$ is greater than $t$, and so is $M$. Consequently, applying
Bonferroni's inequality, the probability that $M$ is greater than $t$ is at 
least
$$\displaylines{\qquad
  \sum_{n\geq n_2}P\{\, \widehat S_n>(1+\epsilon)s_n\,;\, 
  X_{n-k^*}>t-(1+\epsilon)s_n\,\}
  \hfill\cr\noalign{\vskip 3pt}\hfill
  {}-\sum_{\ss m,n\geq n_2\atop\ss m\not= n} 
    P\{\, X_{n-k^*}>t-(1+\epsilon)s_n\,;\, 
                          X_{m-k^*}>t-(1+\epsilon)s_m\,\} \, .
  \quad\equadef{LBSmallGa}\cr
  }
$$
Since $\widehat S_n$ and $X_{n-k^*}$ are independent, the first sum 
in \LBSmallGa\ is at least
$$
  (1-\epsilon) \sum_{n\geq n_2} \oF\bigl( t-(1+\epsilon)s_n\bigr) \, .
$$
Since $s_n=-\g{n}$, Lemma \BS\ and regular variation of both $\oF$ and
$g$ imply
$$
  \sum_{n\geq n_2}\oF\bigl( t-(1+\epsilon)s_n\bigr)
  \sim \int_1^\infty \oF\Bigl( t+(1+\epsilon){g(1-1/u)\over\Gamma(1+\gamma)}
  \Bigr)\d u
$$
as $t$ tends to infinity. Under the claim to be proved that the second sum 
in \LBSmallGa\ is asymptotically negligible with respect to the first sum,
since $\epsilon$ is arbitrary, upon using Lemma \Karam, we obtain
$$
  \liminf_{t\to\infty} {P\{\, M>t\,\}\over 
  \ds\int_1^\infty \oF\Bigl( t+{g(1-1/u)\over\Gamma(1+\gamma)}\Bigr) \d u }
  \geq 1\, .
  \eqno{\equa{LBSmallGb}}
$$
The second sum in \LBSmallGa\ is less than 
$$
  \sum_{m,n\geq n_2} \oF\bigl( t-(1+\epsilon)s_n\bigr)
  \oF\bigl( t-(1+\epsilon)s_m\bigr)
  =\Bigl( \sum_{n\geq n_2} \oF\bigl( t-(1+\epsilon)s_n\bigr)\Bigr)^2 \, .
$$
By the previous arguments, this last quantity is of smaller order than the 
first sum in \LBSmallGa. This proves \LBSmallGb, which, using Lemma \Karam, 
is the lower bound pertaining to Theorem \Thg.

\bigskip


\subsection{Proof of Theorem \ThG\ -- upper bound}%
Again, we complete all the steps described in section \fixedref{3}. 
To prove this
upper bound, we assume without any loss of generality that $\mu$ is $-1$.

\bigskip

\noindent{\it Step 1.} The asymptotic equivalence (S1) follows from the
equality $\psi_n^{-1}=\Psi_n^{-1}-1$ and Lemma \psiInvU. We
now take $\chi=U$ and $\rho=\xi_{\Gamma(1+\gamma)}$. Thus, following \rDef,
here
$$
  r(t)=(\Id\oF)\circ U(t)\int_0^\infty \xi_{\Gamma(1+\gamma)}^{-\alpha}(v)\d v 
  \, .
$$

\bigskip

\noindent{\it Step 2.} Let $\epsilon$ be a positive real number. Lemma
\psiInvLB\ implies that
$$\eqalign{
  \sum_{n_1\leq n\leq \epsilon U(t)}\oF\circ \Psi_n^{-1}(t)
  &{}\leq \sum_{n_1\leq n\leq \epsilon U(t)} 
    \oF\bigl( U(t)\xi_{\Gamma(1+\gamma)}(0)\bigr) \cr
  &{}\lsim\epsilon\xi_{\Gamma(1+\gamma)}^{-\alpha}(0)(\Id\oF)\circ U(t) \, .\cr
  }
$$
Furthermore, regular variation of $\oF$ and Lemma \psiInvBounds\ show that
$$
  \sum_{n\geq U(t)/\epsilon}\oF\circ\Psi_n^{-1}(t)
  \lsim \sum_{n\geq U(t)/\epsilon}\oF\bigl(n\xi_0(1)\bigr) \, .
$$
This last sum can be approximated by an integral, which, by Karamata's theorem
is asymptotically equivalent to
$$
  {U(t)\over (\alpha-1)\epsilon}\oF\bigl( U(t)\xi_0(1)/\epsilon\bigr)
  \sim {\epsilon^{\alpha-1}\over\alpha-1}\xi_0(1)^{-\alpha}(\Id\oF)\circ U(t)
$$
as $t$ tends to infinity. In view of \rcalculation, this completes step 2.

\bigskip

\noindent {\it Step 3.} Our next lemma shows that the maximum is unlikely
to occur at a time of order smaller than $U(t)$. Its proof is inspired by
that of Lemma 2.4 in Mikosch and Samorodnitsky (2000).

\Lemma{\label{BigGSmallTime}%
  The following limit holds
  $$
  \lim_{\epsilon\to 0} \limsup_{t\to\infty} 
  { P\{\, \exists n\, :\, n\leq \epsilon U(t) \, ,\, S_n\geq t\,\}
    \over
    (\Id\oF)\circ U(t) } = 0 \, .
  $$
}

\Proof Let $Z$ be a random variable having the same distribution as $Z_1$ say.
For $n$ large enough,
$$
  S_n-s_n=\sum_{0\leq j<n} g_j Z_{n-j} \leq 2g_n\sum_{0\leq i<n} |Z_{n-i}| \, .
$$
Therefore, with $N=\lfloor \epsilon U(t)\rfloor$ and $t$ large enough,
$$\displaylines{\qquad
  P\bigl\{\, \exists n\, :\, n\leq \epsilon U(t) \, :\, S_n>t\,\bigr\}
  \hfill\equadef{BigGSmallTimeA}\cr\noalign{\vskip 5pt}\hfill
  \eqalign{
    {}\leq{}& P\Bigl\{\, 2g_N\sum_{0\leq i<N} |Z_{N-i}|>t\,\Bigr\} \cr
    {}={}   & P\Bigl\{\, \sum_{0\leq i<N} |Z_{n-i}|-E|Z|>{t\over 2 g_N}-NE|Z|\,
               \Bigr\} \,. \cr
  }
  \qquad\cr}
$$
To apply the large deviation result of Nagaev (1969 a,b) and Cline and
Hsing (1991)
stated as Lemma A.1 in Mikosch and Samorodnitsky (2000), we check that
for some positive $\delta$ and any $t$ large enough
$$
  {t\over 2 g_N} -NE|Z|>\delta N \, .
  \eqno{\equa{BigGSmallTimeB}}
$$
Since
$$
  g_N\sim { g\Bigl( 1-{\ds 1\over\ds \epsilon U(t)}\Bigr)
            \over 
            \epsilon U(t)\Gamma(\gamma)
          }
  \sim {\epsilon^{\gamma-1}\over\Gamma(\gamma)} {t\over U(t)} \, ,
$$
the left hand side of \BigGSmallTimeB\ is asymptotically equivalent to
$$
  \Bigl({1\over 2}\epsilon^{1-\gamma}\Gamma(\gamma)-\epsilon E|Z|\Bigr)U(t) 
  \, ,
$$
while the right hand side is asymptotically equivalent 
to $\delta\epsilon U(t)$. Since $\gamma$ is at least $1$, we see that 
\BigGSmallTimeB\ holds if $\epsilon$ is small enough and $t$ is large enough.
Hence, applying Lemma A.1 in Mikosch and Samorodnitsky (2000), \BigGSmallTimeA\
is at most
$$\displaylines{\qquad
  2NP\bigl\{\, |Z|-E|Z|\geq (t/2g_N)-NE|Z|\,\bigr\}
  \hfill\cr\noalign{\vskip 3pt}\hfill
  \eqalign{
  {}\sim{}& 2N\oF_*\bigl( (t/2g_N)-NE|Z|\bigr) \cr
  {}\sim{}&2 \epsilon 
           \Bigl( {1\over 2}\epsilon^{1-\gamma}\Gamma(\gamma)
           -\epsilon E|Z|\Bigr)^{-\alpha} (\Id\oF_*)\circ U(t) \, . \cr
  }
  \qquad\cr}
$$
The result follows since $\gamma$ is at least $1$ and 
$\oF_*\asymp \oF$.\hfill\qed

\bigskip

\noindent{\it Step 4.} We now prove that the maximum of the process is 
unlikely to occur at a time of order larger than $U(t)$.

\Lemma{\label{BigGLargeTime}%
  The following holds,
  $$
    \lim_{\epsilon\to 0} \limsup_{t\to\infty}
    { P\{\, \exists n\, :\, n\geq U(t)/\epsilon \, , \, S_n>t\,\}
      \over
      (\Id\oF)\circ U(t)
    } = 0 \, .
  $$
}

\Proof It suffices to prove \StepFourAlt. Thus, we need to evaluate
$$\displaylines{\qquad
  P\Bigl( \bigl\{\,\exists n\, :\, n\geq U(t)/\epsilon \, ,\, S_n>t \,\bigr\}
  \hfill\cr\hfill
  \bcap 
          \bigl\{\, 
             \forall n\geq 1 \, , \, X_n\leq \theta 
             \bigl( n+U(t)/\epsilon\bigr)
           \,\,\bigr\}\Bigr) \, .
  \qquad\equa{BigGLargeTimeA}\cr}
$$
Define
$$
  a_n=\theta \bigl(n+U(t)/\epsilon\bigr) -\mu \, ,
$$
so that $X_n\leq \theta \bigl(n+U(t)/\epsilon\bigr)$ is equivalent to
$Z_n\leq a_n$. If $Z_n$ is at most $a_n$ for all $n$, then
$$
  S_n-s_n =\sum_{0\leq i<n} g_i Z_{n-i}\II\{\, Z_{n-i}\leq a_{n-i}\,\} \, .
$$
Hence \BigGLargeTimeA\ is at most
$$
  \sum_{n\geq U(t)/\epsilon} 
  P\Bigl\{\, \sum_{0\leq i<n} g_i Z_{n-i}\II\{\, 
  Z_{n-i}\leq a_{n-i}\,\} >t-s_n\,\Bigr\} \, .
$$
We use \LargeDev\ to bound each probability involved in this
sum, that is to bound
$$
  \log P\Bigl\{\, \sum_{0\leq i<n} g_i Z_{n-i}\II\{\, 
  Z_{n-i}\leq a_{n-i}\,\} >t-s_n\,\Bigr\} \, .
  \eqno{\equa{BigGLargeTimeB}}
$$
In \LargeDev, we take 
$$
  \lambda=c{\log g(1-1/n)\over t-s_n}
$$
for some positive number
$c$ to be determined later. On the range $n\geq U(t)/\epsilon$, our chosen
$\lambda$ tends to $0$ as $t$ tends to infinity, uniformly in $n$. Moreover,
as $t$, and hence $n$, tends to infinity,
$$
  \lambda \leq c{\log g(1-1/n)\over -s_n} 
  \sim c \Gamma (1+\gamma)
  {\log g(1-1/n)\over g(1-1/n)} \, .
$$
Using the Karamata representation (Bingham, Goldie and Teugels, 1989, 
Theorem 1.3.1), $\log g(1-1/n)\sim \gamma\log n$ as $n$ tends to infinity. 
Therefore on the range $n\geq U(t)/\epsilon$, since we assume that 
the sequence $(g_n)_{n\geq 0}$ is asymptotically equivalent to a monotone
sequence,
$$
  \lambda \max_{0\leq i\leq n} g_i 
  \sim \lambda g_n 
  \lsim c\gamma^2 {\log n\over n} \, . 
$$
In particular, uniformly on that range of $n$,
$$
  \max_{0\leq i<n} -EZ\II\{\, \lambda g_iZ\leq\log (1-\epsilon)\,\}
  = o(1)
$$
as $t$ tends to infinity. Furthermore, referring to the last sum involved in
\LargeDev,
$$
  \max_{0\leq i<n}\lambda g_i a_{n-i}
  \leq \lambda \max_{0\leq i<n} g_i a_n
  \lsim 2 c\gamma^2\theta \log n \, .
$$
Therefore, since by Potter's bounds 
$\oH(n/\log n)= O\bigl((\log n) /n\bigr)^{\alpha-\epsilon}$, an application 
of \LargeDev\ show that \BigGLargeTimeB\ is at most
$$\displaylines{\qquad
  -c\log g(1-1/n) +\epsilon c E|Z|\log g(1-1/n)
  \hfill\cr\hfill
  + \log g(1-1/n)o(1) + n^{3\theta c\gamma^2+1-\alpha}
  \log^{\alpha-\epsilon} n \, .
  \qquad\cr}
$$
Taking $\epsilon$ small enough, we first choose $c$ so that, say,
$$
  -c + \epsilon c E|Z| \leq -\alpha-3 \, ,
$$
and then $\theta$ small enough so that
$$
  3\theta c\gamma^2+1-\alpha < 0 \, .
$$
Then, as $t$ tends to infinity and uniformly in $n\geq U(t)/\epsilon$,
$$
   P\Bigl\{\, \sum_{0\leq i<n} g_i Z_{n-i}\II\{\, 
  Z_{n-i}\leq a_{n-i}\,\} >t-s_n\,\Bigr\} 
  = O(n^{-\alpha-2}) \, .
$$
Since $\Id^{-\alpha-1}=o(\Id\oF)$, it follows that \BigGLargeTimeA\
is $o(\Id\oF)\circ U(t)$ as $t$ tends to infinity, and this
concludes the proof of the lemma.\hfill\qed

\bigskip

\noindent{\it Step 5.} We can now prove that if the process exceeds $t$ at 
a time between $\epsilon U(t)$ and $U(t)/\epsilon$, then at least one of 
the $X_i$ has to exceed $U(t)$.

\Lemma{\label{BigGAllZSmall}
  For any positive $\epsilon$ and $p$, there exists $\theta$ such that
  $$
    P\Bigl( \bcup_{n\in I} 
    \Bigl(\{\, S_n>t\,\} \cap\bcap_{0<j\leq n} 
    \{\, Z_j\leq \theta U(t)\,\} \Bigr)\Bigr) = o(t^{-p}) \, .
  $$
}

\Proof Recall that $n_0$, defined before Lemma \psiMonotone, is an integer
such that whenever $n$ is at least $n_0$, both $g_n$ and $\g{n}$ are 
nonnegative.
Let $t$ be sufficiently large so that $\epsilon U(t)$ is at least $n_0$.
If all $Z_j$, $0<j\leq n$, are at most $\theta U(t)$, then
the event $S_n>t$ occurs if and only if
$$
  \sum_{0\leq j<n} g_jZ_{n-j}\II\{\, Z_{n-j}\leq\theta U(t)\,\} > t-s_n \, .
$$
Since $t-s_n$ exceeds $t$, the logarithm of the probability of that
event is bounded as in the next lemma by (\fixedref{6.3.5}) hereafter, with
$\delta=1$ say. We conclude as in the proof of the next lemma by 
taking $\lambda=ct^{-1}\log t$
with $c$ large enough and using Bonferroni's inequality.~ \hfill\qed

\bigskip

\noindent{\it Step 6.} As for the proof of Theorem \Thg, we complete only
the one-sided version of step 6, namely the version where $D_n$ is defined as
$$
  D_n=\bcap_{1\leq i\leq N}\{\, S_i-\psi_{i,n}(X_n)\leq\delta t\,\} \,.
$$
As we argued when proving Lemma \LLN, Theorem 3 in 
Jamison, Orey and Pruitt
(1965) yields a weak law of large numbers on the weighted 
sum $\sum_{0\leq j <n}g_j Z_{n-j}/g_{[0,n)}$ as $n$ tends to infinity. 
The next lemma shows that this weak law of large numbers holds with some
uniformity with respect to the weights.

\Lemma{\label{BigGLLNUniform}%
  Let $\delta$ be a positive real number. For any $\epsilon$ positive small
  enough there exists a positive $\theta$ such that
  $$
    P\Bigl\{\, \bcup_{n\in I} B_n^\complement\cap C_n\cap D_n^\complement
     \,\Bigr\}
    = o\bigl(r(t)\bigr)
  $$
  as $t$ tends to infinity.
}

\bigskip

\Proof Let $N$ be $\lfloor U(t)/\epsilon\rfloor$. On $C_n$,
$$
  \sum_{\ss 0\leq j<i\atop\ss i-j\not=n} g_jZ_{i-j}
  = \sum_{\ss 0\leq j<i\atop\ss i-j\not=n} g_jZ_{i-j}\II\{\, Z_{i-j}\leq
  \theta U(t)\,\} \, .
$$
In the following, we use the one-sided form of $D_n$, 
namely 
$$
  D_n=\bcap_{i\leq N}\{\, S_i-\psi_{i,n}(X_n)\leq\delta t\,\} \, .
$$
Applying inequality \LargeDev\ and 
using that $n$ is at most $N$, the 
logarithm of the probability that $B_n^\complement\cap C_n\cap D_n^\complement$
occurs is at most
$$\displaylines{\qquad
  -\lambda\delta t +\eta\lambda g_{[0,N)}E|Z|
  -\lambda \sum_{0\leq j<N} g_j 
  EZ\II\{\, \lambda g_jZ\leq \log (1-\eta)\,\}
  \hfill\cr\hfill
  +\sum_{0\leq j <N} e^{\lambda g_j\theta U(t)} 
  \oH\Bigl({\log (1+\eta)\over \lambda g_j}\Bigr) \, . \qquad
  \equadef{BigGLLNUniformA}\cr
  }
$$
Note that because we used $N$ in this bound, it holds for all $n$ such
that $\epsilon U(t)\leq n\leq U(t)/\epsilon$.

We choose $\lambda =ct^{-1}\log t$ where $c$ will be specified later. With
this choice, we examine all terms in \BigGLLNUniformA. We have
$\lambda\delta t= c\delta\log t$. Moreover, Lemma \BS\ implies
$$
  \lambda g_{[0,N)}
  \sim c {\log t\over t}\, {g\bigl(1-\epsilon/U(t)\bigr)\over\Gamma (1+\gamma)} 
  \sim {c\epsilon^{-\gamma}\over\Gamma (1+\gamma)}\log t \, .
$$
It also implies
$$
  \lambda g_N
  \sim c {\log t\over t}\, {g\bigl( 1-\epsilon/U(t)\bigr)\over
    \Gamma(\gamma) U(t)/\epsilon}
  \sim {c\epsilon^{-\gamma+1}\over\Gamma (\gamma)} {\log t\over U(t)} \, .
$$
Therefore, since $(\log t)/U(t)$ tends to $0$ at infinity, the bound
$$
  -Z\II\{\, \lambda g_j Z\leq \log (1-\eta)\,\}
  \leq -Z\II\{\, 2\lambda g_NZ\leq \log (1-\eta)\,\}
$$
yields
$$
  -\lambda\sum_{0\leq j<N}g_j EZ\II\{\, \lambda g_jZ\leq\log (1-\eta)\,\}
  = o(\lambda g_{[0,N)})\, .
$$
Finally, we have
$$
  \lambda g_N\theta U(t)
  \sim {c\epsilon^{-\gamma+1}\over\Gamma (\gamma)}\theta \log t \, ,
$$
so that for $t$ large enough,
$$\eqalign{
  \sum_{0\leq j<N} e^{\lambda g_j\theta U(t)} 
  \oH\Bigl({\log(1+\eta)\over\lambda g_j}\Bigr)
  &{}\lsim Nt^{2c\epsilon^{-\gamma+1}\theta/\Gamma(\gamma)} 
   \oF\bigl({\log (1+\eta)\over\lambda g_N}\Bigr) \cr
  &{}\sim U(t) t^{2c\epsilon^{-\gamma+1}\theta/\Gamma (\gamma)} 
   \oF\Bigl({U(t)\over\log t}\Bigr) O(1) \, .\cr
  }
$$
We obtain that \BigGLLNUniformA\ is at most equivalent to
$$
  c\log t \Bigl( -\delta +{\eta \epsilon^{-\gamma}\over\Gamma(1+\gamma)}E|Z|
                 +o(1) \Bigr)
  + U(t) t^{2c\epsilon^{-\gamma+1}\theta/\Gamma(\gamma)} 
  \oF\Bigl( {U(t)\over\log t}\Bigr) O(1) \, .
$$
Let $p$ be an arbitrary positive real number. We take $c$ large 
enough and $\eta$ small enough so that
$$
  c\Bigl( -\delta +{\eta\epsilon^{-\gamma}\over\Gamma(1+\gamma)}E|Z| \Bigr)
  \leq -p-2\gamma-3 \, .
$$
We take $\theta$ small enough so that
$$
  U(t) t^{c\epsilon^{-\gamma+1}\theta/2\Gamma(\gamma)} 
  \oF\Bigl( {U(t)\over\log t}\Bigr) = o(1) \, .
$$
Such $\theta$ exists because Potter's bound applied to both $\oF$ and
$U$ ensure that the function $U(t)\oF\bigl( U(t)/\log t\bigr)$ tends to $0$ at
infinity at a rate at least some positive power of $1/t$. Therefore, we obtain,
as $t$ tends to infinity,
$$
  \sup_{n\in I}P(B_n^\complement\cap C_n\cap D_n^\complement) 
  = o(t^{-p-2\gamma-2}) \, .
$$
The result follows by an application of Bonferroni's inequality, upon using
Potter's bound to bound $U(t)$ and taking $p$ to be greater 
than the negative of the index of regular variation of $r$, 
that is $(1-\alpha)/\gamma$ here.\hfill\qed

\bigskip

Having completed steps 1 through 5 indicated in section \fixedref{3.1} as 
well as the one-sided
version of step 6, the upper bound result follows by an application of Theorem 
\UsefulStatement.\hfill\qed

\bigskip


\subsection{Proof of Theorem \ThG\ -- lower bound}%
As with the proof of Theorem \Thg, proving the lower bound by
an application of Theorem
\UsefulStatement\ upon completing step 7 requires a tail balance condition
on the distribution function $F$ in order to prove the two-sided version of
Step 6. This extra assumption is not needed
with the following proof, again adapted from Zachary's (2004) work.

We assume without loss of generality, that $\mu$ is $-1$.
Let $\epsilon$ be a positive real number less than $1$ and consider the 
corresponding set $I$. Let $\delta$ be a positive real number, and let $p$ be
an integer depending on $n$, $t$ and $\delta$, such 
that $g_p\psi_n^{-1}\bigl((1+\delta)t\bigr)+s_{n+p}$ is maximum and 
therefore asymptotically equivalent to $(1+\delta)t$. For $n$ 
in $I$, Lemma \psiInvU\ shows that $\psi_n^{-1}\bigl((1+\delta)t\bigr)$ 
is of order 
$U(t)$. Arguments very similar to that of the proof of Lemma \psiInvU\ 
show that $p/U(t)$ remains in a compact subset of the nonnegative half-line
when $n$ stays in $I$. Then Lemma \BS\ implies that
$s_{n+p}/t$ stays bounded over $n$ in $I$ and as $t$ tends to infinity.
Therefore, we can find $\eta$ small enough so that $\min_{n\in I}\delta t+\eta
s_{n+p}$ is positive for any $t$ large enough.

We consider the events
$$
  A_{n,t}=\{\, S_{n+p}-g_p X_n \geq (1+\eta) s_{n+p} \hbox{ and }
  X_n\geq\psi_n^{-1}\bigl((1+\delta)t\bigr)\,\} \, .
$$
If $A_{n,t}$ occurs, then
$$\eqalign{
  S_{n+p}
  &{}\geq g_p X_n+(1+\eta)s_{n+p} \cr
  &{}> g_p\psi_n^{-1}\bigl((1+\delta)t\bigr)+(1+\eta)s_{n+p}  \cr
  &{}\geq (1+\delta)t +\eta s_{n+p}\cr
  &{}\geq t \, , \cr
  }
$$
and therefore $M\geq t$. Consequently, provided $t$ is large enough, the 
inclusion $A_{n,t}\subset \{\, M>t\,\}$ holds for every $n$. It follows 
that for $t$ large enough,
$$\eqalignno{
  P\{\, M>t\,\} 
  &{}\geq P\bigl(\bigcup_{n\in I}A_{n,t}\bigr) \cr
  &{}\geq \sum_{n\in I} P(A_{n,t}) 
   -\sum_{\ss n,m\in I\atop\ss m\not= n} P(A_{n,t}\cap A_{m,t}) \, .
  &\equadef{LBBigGa}\cr
  }
$$
Since $S_{n+p}-g_pX_n$ and $X_n$ are independent, for $n$ larger than $n_0$,
$$
  P(A_{n,t}) 
  = P\{\, S_{n+p}-g_p X_n > (1+\eta) s_{n+p}\,\}
  P\bigl\{\, X_n\geq \psi_n^{-1}\bigl((1+\delta)t\bigr)\,\bigr\} \, .
$$
Therefore, Lemma \LLN\ implies that provided $n$ is large enough,
$P(A_{n,t})$ is at least 
$(1-\epsilon)\oF\circ\psi_n^{-1}\bigl((1+\delta)t\bigr)$. Hence, 
if we can prove 
that the second sum in \LBBigGa\ is negligible compare to the first one, then
$$
  P\{\, M> t\,\} 
  \geq (1-\epsilon)\sum_{n\in I}\oF\circ\psi_n^{-1}
  \bigl((1+\delta)t\bigr) \, .
$$
Then, since $\epsilon$ and $\delta$ are arbitrary, the arguments used to 
derive \rDef\ shows that, in view of $\chi=U$ and 
$\rho=\xi_{\Gamma(1+\gamma)}$,
$$
  \liminf_{t\to\infty} {P\{\, M> t\,\}\over (\Id\oF)\circ U(t)}
  \geq \int_0^\infty\xi_{\Gamma (1+\gamma)}^{-\alpha}(v) \d v \, .
$$
But the double sum in \LBBigGa\ is at most the square of the first one and 
hence is of order $o(\Id\oF)\circ U(t)$.\hfill\qed

\bigskip


\subsection{Proof of Theorem \ThOne}%
Part of the proof is analogous to that of Theorem \Thg. 
Convergence of the sequence $(g_n)_{n\geq 0}$ to $g_\infty$ implies 
asymptotic equivalence $g(1-x)\sim g_\infty/x$ as $x$ tends to $0$, and
therefore $U\sim \Id/g_\infty$ at infinity.

\bigskip

\noindent {\it Step 1.} Argue as in step 1 of Theorem \Thg, using Lemma
\PsiInvUO\ instead of Lemma \PsiInvUg, to show that we may take $\chi\sim \Id$ 
at infinity, which yields $\psi^{-1}_{\lfloor xU(t)\rfloor}
(t)\sim t (1-\mu x)/g^*$ as $t$ tends to infinity.

\bigskip

\noindent {\it Step 2.} The arguments used in step 2 of the proof of
Theorem \Thg\ carry over, substituting Lemma \PsiInvLbO\ for Lemmas 
\PsiInvLbgA\ and \PsiInvLbgB.

\bigskip

\noindent {\it Step 3.} We prove that the process is very unlikely to 
reach the level $t$ at a time of smaller order than $t$.

\Lemma{\label{OStepThree}%
  The following limit hold,
  $$
    \lim_{\epsilon\to 0}\limsupt 
    {P\{\, \exists n \, :\, 0\leq n\leq \epsilon U(t) 
           \, ,\, S_n>t\,\}\over t\oF(t) } = 0 \, .
  $$
}

\Proof Arguing as in the beginning of Lemma \SmallGStepThree 
--- see \SmallGStepThreeA\ --- it suffices to
prove that for any positive $\delta$, there exists some 
positive $\theta$ such that
$$
  \lim_{\epsilon\to 0}\limsupt \epsilon t\max_{1\leq n\leq \epsilon t}
  { P\{\, \sum_{0\leq j<n} g_jZ_{n-j}\II\{\, Z_{n-j}\leq \theta t\,\}>\delta t
     \,\}
    \over t\oF(t) } = 0 \, .
  \eqno{\equa{OStepThreeA}}
$$
Let $c$ be a positive real number to be determined later.
Using inequality \LargeDev\ with $\lambda=c(\delta t-s_n)^{-1}\log t$ 
and $a_i=\theta t$, the logarithm of the ratio
$$
  P\Bigl\{\, \sum_{0\leq j<n} g_jZ_{n-j}\II\{\, Z_{n-j}\leq \theta t\,\}
             >\delta t \,\Bigr\}
  \bigm/ \oF(t)
$$
is ultimately at most
$$\displaylines{\quad
  -c(1+\epsilon/\delta)^{-1}\log t + \eta cE|Z|\log t 
  \hfill\cr\noalign{\vskip 3pt}\qquad\qquad
  {}- {c\log t\over\delta t-s_n}\sum_{0\leq j<n}
  g_jEZ\II\Bigl\{\, {c\log t\over \delta t-s_n}g_j Z\leq \log (1-\eta)\,\Bigr\}
  \hfill\cr\qquad\qquad
  {}+\sum_{0\leq j<n} e^{\lambda g^*\theta t}\oH\Bigl( {\log (1+\eta)\over g^*}
  {\delta t\over c\log t}\Bigr) + (\alpha+\epsilon)\log t \, .
  \hfill\equa{OStepThreeB}\cr}
$$
Referring to the third summand in this bound, since $n\leq \epsilon t$, it is
at most
$$
  -{c\log t\over\delta} g^*\epsilon E|Z|
  \II\Bigl\{\, Z\leq {\delta \log (1-\eta)\over g^* c}{t\over\log t}\,\Bigr\}
  = o(\log t) \, .
$$
Let $\nu$ be a positive real number such that $1-\alpha+2\nu$ is negative.
The fourth summand is ultimately at most
$$\displaylines{\qquad
  n e^{(cg^*\theta/\delta)\log t} \oH\Bigl( {\log (1+\eta)\over g^*} 
  {\delta t\over c\log t}\Bigr)
  \hfill\cr\hfill
  \eqalign{
  {}\leq{}& \epsilon t^{1+cg^*\theta/\delta} 
            \oH\Bigl( {\log (1+\eta)\over g^*}{\delta t\over c\log t}\Bigr) \cr
  {}\leq{}& \epsilon t^{1+(cg^*\theta/\delta)-\alpha+\nu} \, ,\cr
  }
  \qquad\cr}
$$
where we used Potter's bound to obtain the last inequality. Therefore,
\OStepThreeB\ is ultimately at most
$$
  \bigl( -c(1+\epsilon/\delta)^{-1}+\eta cE|Z|+\alpha+\epsilon+o(1)\bigr) 
  \log t 
  + \epsilon t^{1+(cg^*\theta/\delta)-\alpha+\nu} \, .
$$
We take 
$$
  \eta=1/2E|Z| \, ,\qquad 
  \epsilon=\delta/3 \, , \qquad
  c=8\alpha
$$ 
and $\theta$ small enough so that $cg^*\theta/\delta\leq \nu$, which, 
given how $\nu$ was defined, guarantees that
$1+(cg^*\theta/\delta)-\alpha+\nu$ is negative. This proves \OStepThreeA\ as 
well as the lemma.\hfill\qed

\bigskip

\noindent {\it Step 4.} We need to prove that the process is very unlikely
to reach the level $t$ at a time of larger order than $t$. This follows
from Lemma \SmallGStepFour\ whose proof, and hence, conclusion, remains
valid in the present context.

\bigskip

\noindent {\it Step 5.} Similarly to the previous step, Lemma \SmallGStepFive\
remains valid in the present context.

\bigskip

\noindent {\it Step 6.} Similarly to the previous step, Lemma \SmallGStepSix\
remains valid.

\bigskip

An application of Theorem \UsefulStatement\ yields the upper bound. The proof
of the lower bound of Theorem \Thg\ carries over in the present setting, and
this concludes the proof of Theorem \ThOne.

\bigskip


\subsection{Proof of Theorem \VgSigned}%
We only sketch the proof. Assume 
without loss of generality that $\mu=-1$ and define as before
$$
  \psi_{i,n}(x)
  =E(S_i\mid X_n=x)
  =\cases{ s_i & if $i<n$, \cr
           (x+1)g_{i-n}+s_i & if $i\geq n$.\cr}
$$
For $x$ positive, we define
$$
  \psi_{+,i}(x)=\max_{\ss k\geq 0\hfill\atop\ss g_k>0\hfill} xg_k+s_{i+k}
$$
and for $x$ negative define
$$
  \psi_{-,i}(x)=\max_{\ss k\geq 0\hfill\atop\ss g_k<0\hfill} xg_k+s_{i+k} \, .
$$
By the same arguments as in our heuristic, we expect to prove that to 
reach the level $t$, either $\psi_{+,n}(X_n)>t$ or $\psi_{-,n}(X_n)>t$
for some $n$.
The actual proof can be done by redefining $B_n$ as the two-sided event
$\{\, |X_n|\leq\theta\chi(t)\,\}$ and using the tail balance condition.
Thus we have
$$
  P\{\, M>t\,\}
  \sim \sum_{n\geq n_1} \oF\circ\psi_{+,n}^{-1}(t)
  +\sum_{n\geq n_1}F\circ\psi_{-,n}^{-1}(t) \, .
$$
Similarly to what we proved previously, one has
$$
  \psi_{+,n}(x)\sim g^*x+s_n
$$
as $x$ tends to infinity and
$$
  \psi_{-,n}(x)\sim g_*x+s_n
$$
as $x$ tends to minus infinity. The result follows as in the proof of Theorem
\Thg.\hfill\qed

\bigskip


\subsection{Proof of Theorem \TrajG}%
The theorem is proved by applying Theorem \LimitingTraj. The tail balance
condition \TailBalance\ guarantees that (S6) holds with the two-sided event
$D_n$; this can be seen by exactly the same arguments we used to prove
Lemma \BigGLLNUniform, using two-sided versions of the events $B_n$. 
Thus, it remains to show that \PathA\ holds.

We first asume without loss of generality that $\mu=-1$. Equality \psiin\
shows that for any positive real number $\lambda$, $\tau$ and $y$,
\finetune{\hfuzz=8pt}
$$\displaylines{\quad
  h_t(\lambda,\tau,y)
  \hfill\cr\qquad
  {}=t^{-1}\psi_{\lfloor\lambda U(t)\rfloor,\lfloor \tau U(t)\rfloor}
   \bigl(\chi(t)y\bigr)
  \hfill\cr\qquad
  {}=-t^{-1}\g{\lambda U(t)}
  \hfill\cr\qquad\qquad
    {}+\II\{\, \lfloor\lambda U(t)\rfloor\geq\lfloor\tau U(t)\rfloor\,\}
   t^{-1}g_{\lfloor\lambda U(t)\rfloor-\lfloor\tau U(t)\rfloor}
   \bigl( \chi(t)y+1\bigr) \, . \hfill\cr
  }
$$
\finetune{\hfuzz=0pt}
Therefore, setting
$$
  h(\lambda,\tau,y)
  = -{\lambda^\gamma\over\Gamma(1+\gamma)} +\II\{\, \lambda\geq\tau\,\}
  {(\lambda-\tau)^{\gamma-1}\over\Gamma(\gamma)} y  \, ,
$$
and using that we take $\chi$ equal to $U$ in this case, we have, using Lemma
\BS, the pointwise convergence
$$
  \limt h_t(\lambda,\tau,y)
  = h(\lambda,\tau,y) \, .
$$
Let $\epsilon$ be a positive real number. We prove that $h_t(\,\cdot\,,\tau,y)$
tends to $h(\,\cdot\,,\tau,y)$ in $\D[\,0,1/\epsilon\,]$. For this, note that
$$
  \II\{\, \lfloor\lambda U(t)\rfloor\geq\lfloor\tau U(t)\rfloor\,\}
  = \II\Bigl\{\, {\tau U(t)\over \lfloor\tau U(t)\rfloor}\lambda \geq\tau
  \,\Bigr\} \,  .
$$
Set $v_t(\lambda)=\lambda \lfloor\tau U(t)\rfloor/\bigl(\tau U(t)\bigr)$.
Then Lemma \BS\ implies that, as $t$ tends to infinity,
$$\displaylines{
  h_t\bigl(v_t(\lambda),\tau,y\bigr)-h(\lambda,\tau,y)
  =-t^{-1}\g{v_t(\lambda) U(t)}+{\lambda^\gamma\over\Gamma(1+\gamma)} 
  \hfill\cr\hfill
  {}+\II\{\,\lambda\geq\tau\,\} 
    \Bigl( t^{-1}g_{ \left\lfloor{\lfloor\tau U(t)\rfloor\over\tau}\lambda
                     \right\rfloor
                     -\lfloor\tau U(t)\rfloor
                    } 
           \bigl(\chi(t)y+1\bigr)
           -{(\lambda-\tau)^{\gamma-1}\over\Gamma(\gamma)}y
    \Bigr)
  \cr
  }
$$
tends to $0$ uniformly in $\lambda$ in any fixed compact subset of the positive
half-line. Moreover, since $\gamma$ is positive and the $g_i$ are
ultimately positive, Lemma \BS\ also shows 
that $t^{-1}\g{\lambda U(t)}$ tends to $\lambda^\gamma/\Gamma(1+\gamma)$
uniformly on any interval of the form $[\,0,1/\epsilon\,]$. Therefore, taking
$\epsilon$ to be less than $\tau$, this shows that
$h_t\bigl(v_t(\lambda),\tau,y\bigr)-h(\lambda,\tau,y)$ converges uniformly
to $0$ on $[\,0,1/\epsilon\,]$.
Since $v_t$
tends to the identity uniformly in $[\,0,1/\epsilon\,]$, it follows from the
definition of the Skorohod topology (see Billingsley, 1968, definition of
the distance $d$ in section 14) that for every $\tau$ and $y$ the
function $h_t(\,\cdot\,,\tau,y)$ converges to $h(\,\cdot\,,\tau,y)$ 
in $\D[\,0,1/\epsilon\,]$.

We then apply Theorem \LimitingTraj\ to obtain Theorem \TrajG\ when $\mu=-1$.
For a general negative mean $\mu$, let $(X_i)_{i\geq 0}$ be as before a 
sequence of independent and identically distributed random variables 
with mean $-1$ and let 
$\widetilde X_i=(-\mu)X_i$, $i\geq 0$. We agree to cover by a tilde whatever
quantity is calculated on the $\widetilde X_i$ and to leave 
uncovered quantities calculated on the
$X_i$. Then, with the notation of section \fixedref{2},
\setbox1=\vbox{\baselineskip=18pt\par
\halign{$#$\hfill\cr
\widetilde S_n=(-\mu)S_n \, ,\cr
\widetilde M>t \hbox{ if and only if } M>t/(-\mu) \, , \cr
\widetilde N_t=N_{t/(-\mu)} \, , \cr
\widetilde J_t=J_{t/(-\mu)} \, ,\cr
\widetilde \tau_t={\ds\widetilde J_t\over\ds U(t)}
                 ={\ds U(t/(-\mu))\over\ds U(t)} \tau_{t/(-\mu)}\, , \cr
\widetilde Y_t={\ds\widetilde X_{\widetilde J_t}\over\ds U(t)}
              =(-\mu) {\ds U\bigl(t/(-\mu)\bigr)\over\ds U(t)} Y_{t/(-\mu)} 
              \, , \cr
\widetilde\fS_t(\lambda)=-\mu {\ds S_{\lfloor\lambda U(t)\rfloor}\over\ds t}
                        = \fS_{t/(-\mu)}
                        \Bigl(\lambda {\ds U(t)\over \ds U\bigl(t/(-\mu)\bigr)}
                        \Bigr) \, .\cr
}
}

$$\box1$$

\noindent It follows that the limiting random 
variables $(\widetilde\fS,\widetilde\tau,\widetilde Y)$ satisfy
$$
  \widetilde\fS=\fS\bigl((-\mu)^{1/\gamma} \,\cdot\,) \, ,\quad
  \widetilde\tau=(-\mu)^{-1/\gamma}\tau\, \quad \hbox{and}\quad
  \widetilde Y=(-\mu)^{1-1/\gamma}Y \, ,
$$
and this completes the proof of Theorem \TrajG.\hfill\qed

\bigskip

\noindent{\bf Acknowldegements.} Ph.B.\ thanks Florin Avram for showing him
some works on combinatorics which was particularly enlightening and which,
a year later, was very inspiring when working on this paper.

\bigskip


\noindent{\bf References}
\medskip

{\leftskip=\parindent \parindent=-\parindent
 \par

\hfuzz=2pt
J.\ Akonom, Chr.\ Gouri\'eroux (1987). A functional central limit theorem
for fractional processes, discussion paper 8801, CEPREMAP, Paris.

S.\ Asmussen (1987). {\sl Applied Probability and Queues}, Wiley.

S.\ Asmussen (2000). {\sl Ruin Probabilities}, World Scientific.

S.\ Asmussen, C.\ Kl\"uppelberg (1996). Large deviation results for 
subexponential tails, with applications to insurance risk, {\sl Stoch.\ Proc.\
Appl.}, 64, 103--125.

S.\ Asmussen, H.\ Schmidli, V.\ Schmidt (1999). Tail probabilities for 
non-standard risk and queueing process with subexponential jumps, {\sl Adv.\
Appl.\ Probab.}, 31, 422--447.

F.\ Baccelli, S.\ Foss (2004). Moments and tails in monotone-separable 
stochastic networks, {\sl Ann.\ Probab.}, 14, 612--650.

F.\ Baccelli, S.\ Foss, M.\ Lelarge (2005). Tails in generalized Jackson 
networks with subexponential service-time distributions, {\sl J.\ Appl.\ 
Probab.}, 42, 513--530.

\hfuzz=0pt
Ph.\ Barbe, M.\ Broniatowski (1998). Note on functional large deviation 
principle for fractional ARIMA processes, {\it Statistical Inference
for Stochastic Processes}, 1, 17--27.

Ph.\ Barbe, A.-L.\ Foug\`eres, Chr.\ Genest (2006). On the tail behavior of 
sums of dependent risks, {\sl Astin Bulletin}, 36, 361--373.

J.\ Beran (1994). {\sl Statistics for Long-Memory Processes}, Chapman and Hall.

P.\ Billingsley (1968). {\sl Convergence of Probability Measures}, Wiley.

N.H.\ Bingham, C.M.\ Goldie, J.L.\ Teugels (1989). {\sl Regular Variation}, 
2nd ed. Cambridge University Press.

E.\ Bolthausen (1993). Stochastic processes with long range interactions 
of the paths, {\sl Contemporary Mathematics}, 149, 297-319.

A.A.\ Borovkov (1971). {\it Stochastic Processes in Queueing Theory}, Springer.

A.A.\ Borovkov (2003). Large deviation probabilities for random walks in 
the absence of finite expectation jumps, {\sl Probab.\ Theor.\ Relat.\ Fields},
125, 421--446.

B.L.J.\ Braaksma, D.\ Stark (1997). A Darboux-type theorem for slowly
varying functions, {\sl J.\ Comb.\ Theor., A}, 77, 51--66.

S.S.\ Choi, Th.M.\ Cover, I.\ Csiz\'ar (1987). Conditional limit theorem under
Markov conditioning, {\sl IEEE Trans.\ Inform.\ Theory}, IT-33, 788--801.

D.B.H.\ Cline, T.\ Hsing (1991). Large deviation probability for
sums of random variables with heavy or subexponential tails, preprint,
Texas A\&M Univ.

A.\ Dembo, O.\ Zeitouni (1992). {\sl Large Deviation Techniques and
Applications}, Bartlett.

P.\ Embrechts, D.D.\ Lambrigger, M.V.\ W\"utrich (2008). Multivariate extremes
and the aggregation of dependent risks: example and counter-examples, preprint.

P.\ Embrechts, J.\ Ne\v slehov\'a, M.V.\ W\"utrich (2008). Additivity 
properties for Value-at-Risk under Archimedean dependence and heavy-tailedness,
{\sl Insurance, Mathematics and Economics}, to appear.

P.\ Embrechts, N.\ Veraverbeke (1982). Estimates for the probability of ruin
with special emphasis on the possibility of large claims, {\it Insurance:
Math.\ Economics}, 1, 55--77.

S.\ Foss, Z.\ Palmowski, S.\ Zachary (2005). The probability of exceeding
a high boundary on a random time for a heavy-tailed random walk, {\sl Ann.\
Appl.\ Probab.}, 15, 1936--1957.

C.W.J.\ Granger (1980). Long memory relationships and the aggregation of
dynamic models, {\sl J.\ Econometrics}, 14, 227--238.

B.\ Jamison, S.\ Orey, W.\ Pruitt (1965). Convergence of weighted averages
of independent random variables, {\sl Zeit.\ Wahrsch.\ Theor.\ Verw.\ Geb.},
4, 40--44.

P.S.\ Kokoszka, M.S.\ Taqqu (1995). Fractional ARIMA with stable innovations,
{\sl Stoch.\ Proc.\ Appl.}, 60, 19--47.

D.\ Konstantinides, T.\ Mikosch (2005). Large deviations and ruin probabilities
for solutions of stochastic recurrence equations with heavy-tailed innovations,
{\sl Ann.\ Probab.}, 33, 1992--2035.

D.\ Korshunov (1997). On distribution tail of the maximum of a random walk,
{\sl Stoch.\ Proc.\ Appl.}, 72, 97--103.

A.J.\ McNeail, R.\ Frey, P.\ Embrechts (2005). {\it Quantitative Risk 
Management: Concept, Techniques, Tools}, Princeton University Press.

T.\ Mikosch, G.\ Samorodnitsky (2000). The supremum of a negative drift random
walk with dependent heavy-tail steps, {\sl Ann.\ Appl.\ Probab.}, 10, 
1025--1064.

A.\ Montanari, R.\ Rosso, M.S.\ Taqqu (1997). Frationally differenced ARIMA
models applied to hydrologic time series: identification, estimation, and 
simulation, {\sl Water Ressources Research}, 33, 1035--1044.

A.\ Nagaev (1969a). Integral limit theorem for large deviations when Cram\'er's
condition is not fulfilled, I, II, {\it Theory Probab.\ Appl.}, 14, 51--64, 
193--208.

A.\ Nagaev (1969b). Limit theorems for large deviations when Cram\'er's 
conditions are violated, {\sl Izv.\ Akad.\ Nauk UzSSR Ser.\ Fiz.-Mat. Nauk},
6, 17--22 (in Russian).

K.W.\ Ng, Q.\ Tang, J.-A.\ Yan, H.\ Yang (2004). Precise large deviations
for sum of random variables with consistently varying tails, {\sl J.\ Appl.\
Prob.}, 41, 93--107.

H.\ Nyrhinen (2005). Power estimates for ruin probabilities, {\sl Adv.\ Appl.\
Probab.}, 37, 726--742.

K.B.\ Oldham, J.\ Spanier (2006). {\sl The Fractional Calculus}, Dover.

A.G.\ Pakes (1975). On the tails of waiting-time distributions, {\sl J.\
Appl.\ Probab.}, 12, 555--564.

P.C.B.\ Phillips (1987). Time series regression with a unit root, 
{\sl Econometrica}, 55, 277--301.

D.\ Pollard (1984). {\sl Convergence of Stochastic Processes}, Springer.

S.I.\ Resnick (2007). {\sl Heavy-Tail Phenomena, Probabilistic and Statistical
Modeling}, Springer.

G.\ Samorodnitsky, M.S.\ Taqqu (1992). Linear models with long-range dependence
and finite or infinite variance, in {\sl New Directions in Time Series 
Analysis, Part II}, D.\ Brillinger, P.\ Caines, J.\ Geweke, E.\ Parzen, M.\
Rosenblatt eds., Springer.

K.\ Tanaka (1999). The nonstationary fractional unit root, {\sl Econometric 
Theory}, 15, 549--582.

S.R.S.\ Varadhan (1966). Asymptotic probabilities and differential equations,
{\sl Comm.\ Pure Appl.\ Math.}, 19, 261--286.

N.\ Veraverbeke (1977). Asymptotic behavior of
Wiener-Hopf factors of a random walk. {\sl Stoch. Proc. Appl.}, 5, 27--37.


Q.\ Wang, Y.-X.\ Lin, Ch.\ Gulati (2002). Asymptotics for general 
nonstationary fractionally integrated process without prehistoric
influence, {\sl J.\ Appl.\ Math.\ Decision Sci.}, 6, 255-269.

M.\ Woodroofe (1982). {\sl Nonlinear Renewal Theory in Sequential Analysis},
SIAM.

S.\ Zachary (2004). A note on Veraverbeke's theorem, {\sl Queuing Systems},
46, 9--14.

}


\bigskip

\setbox1=\vbox{\halign{#\hfil&\hskip 40pt #\hfill\cr
  Ph.\ Barbe            & W.P.\ McCormick\cr
  90 rue de Vaugirard   & Dept.\ of Statistics \cr
  75006 PARIS           & University of Georgia \cr
  FRANCE                & Athens, GA 30602 \cr
  philippe.barbe@math.cnrs.fr                        & USA \cr
                        & bill@stat.uga.edu \cr}}
\box1

\bye

\vfill\eject

\eightpoints

{\bf LEFT OVER}

\def\Vsum{{\bf Vsum}}
The reason we present the three forms, is that our extension will show that
\Vsum is perhaps the most intuitive one, as well as a necessary step to
obtain the other two. Our analysis will also show that while analogue 
of \Vint can be 
found, its extensions are not quite natural. On the other hand, we will
find analogues of \VK in our more general context, derived from the analogue
of \Vsum through the proper Karamata type theorem.

$$
  P\{\, M>t\,\}\sim \sum_{n\geq 1} \oF(t-n\mu)
  \eqno{\equa{Vsum}}
$$
as $t$ tends to infinity.  Hence, since $\mu$ is negative and $t$
tends to infinity in this asymptotic equivalence, it relates the tail
of the maximum of the random walk to the integrated tail of the
increments. An additional argument, comparing the sum to an integral
and making a change of variable, yields the more elegant statement,
$$
  P\{\, M>t\,\} \sim {1\over -\mu} \int_t^\infty \oF(u)\d u \, .
  \eqno{\equa{Vint}}
$$
as $t$ tends to infinity.
To rewrite the integral equivalent \Vint\ into a more explicit form,
one uses an extra argument, namely Karamata's theorem (Bingham, Goldie and 
Teugels, 1989, Proposition 1.5.10), which implies that for $\alpha$ greater
than $1$,
$$
  \int_t^\infty \oF(u)\d u\sim {t\oF(t)\over \alpha-1}
  \eqno{\equa{Karamata}}
$$
as $t$ tends to infinity. Thus, under the slightly stronger assumption
$\alpha$ is greater than $1$, the abv equivalence allows us to rewrite \Vint\ 
as 

\bigskip


\def\preveq{(\the\sectionnumber .\the\equanumber)}
\def\prevs{\the\sectionnumber .\the\snumber }

\section{Moving boundary problems for some centered random walks}%
Veraverbeke's theorem can also be interepreted in term of so-called
`moving boundary' problems, which we will now explain; and this suggests
possible extensions in different directions.

Consider in this section a sequence $(X_n)_{n\geq 1}$ of independent and
idendically distributed random variables, centered, with common distribution
$F$ having regulary varying tail of index $-\alpha$. The corresponding
random walk is defined by $S_n=X_1+\cdots + X_n$ with the convention that
$S_0=0$. Let $\mu$ be a negative real number, and consider the sequence
$(n\mu)_{n\geq 0}$. Veraverbeke's theorem provides an evaluation of the 
probability that the random walk  crosses the boundary $(n\mu)_{n\geq 0}$ 
shifted by $t$, that is, the probability of the event
$$
  \exists n\geq 0 \, : \, S_n>t- n\mu \, .
$$
This suggests extensions where the linear boundary $-n\mu$ is replaced by a 
more general one. We will consider a regularly varying function $b$ and
a boundary defined by $b_n=b(n)$ and shifted by $t$. Thus, we are
seeking to find an asymptotic equivalent of
$$
  P\{\, \exists n\geq 0 \, :\, S_n>t+b_n \,\}
  \eqno{\equa{BdyA}}
$$
as $t$ tends to infinity. 

Before stating our result, a remark is in order. If the variables $X_i$
have a finite variance, then Strassen's law of the iterated logarithm
shows that the probability in \BdyA\ is $1$ whenever $\lim_{n\to\infty}
b_n/\sqrt{2n\log\log n}=0$. Thus it is natural to consider only functions
$b$ whose index of regular variation is at least $1/2$.
The following results assume the that the index is in fact at least $1$.
Recall that a regularly varying function $b$ of positive index tends 
to infinity at infinity. The function
$$
  b^\leftarrow(t)=\inf\{\, x\,:\, b(x)\geq t\,\} \, .
$$
satisfies $b\circ b^\leftarrow\sim \Id$ and $b^\leftarrow\circ b\sim \Id$,
and is called an asymptotic inverse of $b$ (Bingham, Goldie and Teugel, 1989,
Theorem 1.5.12).

\Theorem{\label{BdyTh}%
  Let $(S_n)_{n\geq 1}$ be a random walk associated to a sequence of
  independent random variables with common distribution $F$ with regularly
  varying tail of index $\alpha$. Let $b$ be a regularly varying function
  of index $\beta$. If $\alpha\beta>1$ and $\liminft b(t)/t>0$, then
  $$
    P\{\, \exists n\geq 0\, :\, S_n>t+b(n)\,\}
    \sim \beta B\Bigl({1\over\beta},\alpha-{1\over\beta}\Bigr) 
    (b^\leftarrow\oF)(t)
  $$
  as $t$ tends to infinity.
}

\bigskip

In case of a linear boundary, $b(x)=ax$ with $a$ positive, we recover
Veraverbeke's result, for in this case $\beta=1$ and
$$
  b^\leftarrow (t)= t/a
$$
while
$$
  \beta B\Bigl({1\over\beta},\alpha-{1\over\beta}\Bigr) 
  = {1\over \alpha-1} \, .
$$

\bigskip


\def\preveq{(\the\sectionnumber .\the\equanumber)}
\def\prevs{\the\sectionnumber .\the\snumber }

\section{Proof of Theorem \BdyTh}%
The proof of Theorem \BdyTh\ follows from the general scheme described
in section \fixedref{5} and in particular from Theorem \UsefulStatement.
Thus, we need to complete all the steps mentioned in section \fixedref{5}.

We first calculate
$$
  \psi_{i,n}(x)=E(S_i|X_n=x)
  =\cases{ -b_i & if $i<n$,\cr
           -b_i+x & if $i\geq n$.\cr}
$$
Therefore,
$$
  \psi_n(x)=\max_{1\leq i<n} (-b_i)\vee \max_{i\geq n} (-b_i+x) \, .
$$
It follows that if $\psi_n(x)>\max_{i\geq 1} -b_i$, then 
$$
  \psi_n(x)=\max_{i\geq n} -b_i+x \, .
$$

\bigskip

\noindent{\it Step 1.}
Since $b_n$ is regularly varying with positive index, $\min_{i\geq n}b_i
\sim b_n$ as $n$ tends to infinity, and, consequently, 
$$
  \psi_n(x)\sim -b_n +x
$$
as $n$ tends to infinity and uniformly in $x$. This implies that uniformly
in $t$,
$$
  \psi_n^{-1}(t)\sim t+b_n \, .
$$
Therefore,
$$
  \psi_{\lfloor xb^\leftarrow (t)\rfloor} (t)\sim t(1+x^\beta) \, ,
$$
completing step 1 of section \fixedref{5}.

We see that, with the notation of section \fixedref{5},
$$
  U=b^\leftarrow\, \qquad \chi=\Id \, , \qquad \rho=1+\Id^\beta \, .
$$
Since $\alpha\beta>1$, we have
$$
  r(t)
  =(b^\leftarrow \oF)(t) \int_0^\infty (1+v^\beta)^{-\alpha} \d v
  = \beta B\Bigl({1\over\beta},\alpha-{1\over\beta}\Bigr) (b^\leftarrow\oF)(t) 
  \, .
$$

\bigskip

\noindent{\it step 2.} Note that as $t$ tends to infinity,
$$
  \sum_{n\leq \epsilon b^\leftarrow (t)} \oF\circ\psi_n^{-1}(t)
  \sim \sum_{n\leq \epsilon b^\leftarrow (t)}\oF(t+b_n)
  \leq \epsilon (b^\leftarrow\oF)(t) \, .
  \eqno{\equa{ProofBdyThA}}
$$
Morover, as $t$ tends to infinity,
$$\eqalign{
  \sum_{n\geq b^\leftarrow (t)/\epsilon} \oF\circ \psi_n^{-1}(t)
  &{}\sim\sum_{n\geq b^\leftarrow (t)/\epsilon} \oF (t+b_n) \cr
  &{}\leq \sum_{n\geq b^\leftarrow (t)/\epsilon} \oF(b_n) \cr
  &\sim \int_{b^\leftarrow(t)/\epsilon}^\infty \oF\circ b(u) \d u \, . \cr
  }
$$
Using Karamata's theorem, it follows
$$
  \sum_{n\geq b^\leftarrow (t)/\epsilon}\oF\circ \psi_n^{-1}(t)
  \lsim {\epsilon^{\alpha\beta}\over\alpha\beta-1} (b^\leftarrow\oF)(t) \, .
  \eqno{\equa{ProofBdyThB}}
$$
Estimates \ProofBdyThA\ and \ProofBdyThB\ complete step 2.

\bigskip

\noindent{\it Step 3.} We prove \SThreeAlt, that is,
\finetune{\hfuzz=3pt}
$$\displaylines{\hfill
  \limeps\limsupt {1\over (b^\leftarrow\oF)(t)} 
  P\Bigl\{\, \exists n\leq \epsilon U(t) \, :\, \sum_{1\leq i\leq n}
  X_i\II\{\, X_i\leq \theta t\,\} >t+b_n\,\Bigr\} 
  \cr\hfill
  {}= 0 \, .\qquad
  \equa{ProofBdyThC}
  \cr}
$$
\finetune{\hfuzz=0pt}
Inequality \LargeDev\ shows that the logarithm of
$$
  P\Bigl\{\, \sum_{1\leq i\leq n} X_i\II\{\, X_i\leq\theta t\,\}>t+b_n\,\Bigr\}
  \bigm/ (b^\leftarrow\oF)(t)
$$
is at most
$$\displaylines{\quad
  -\lambda (t+b_n)+\eta\lambda n E|X|
  -\lambda n EX\II\{\, \lambda X\leq \log (1-\eta)\,\}
  \hfill\cr\hfill
  + n e^{\lambda\theta t}\oF\Bigl( {\log (1-\eta)\over \lambda}\Bigr)
  + \Bigl( \alpha-{1\over\beta}\Bigr) \log t \bigl(1+o(1)\bigr) \, .
  \quad \equa{ProofBdyThD}\cr
  }
$$
We take $\lambda$ of the form $ct^{-1}\log t$. If $n\leq \epsilon U(t)$,
then
$$
  b_n\lsim b\bigl( \epsilon U(t)\bigr) \sim \epsilon^\beta t \, .
$$
Thus, for $t$ large enough, \ProofBdyThD\ is at most
$$\displaylines{\qquad
  -c(1+\epsilon^\beta)\log t 
  + \eta ct^{-1}(\log t) \epsilon b^\leftarrow(t) E|X|\bigl( 1+o(1)\bigr)
  \hfill\cr\hfill
  {}+ ne^{c\theta\log t} \oF\Bigl( {\log (1+\eta)\over c\log t} t\Bigr)
  + 2\Bigl(\alpha-{1\over\beta}\Bigr) \, .
  \qquad\cr}
$$
Note also that 
$$
  t^{-1}b^\leftarrow (t)\log t\lsim t^{-1+1/\beta} \log t \lsim \log t
$$
for $\beta$ greater than $1$. Morover, since $n\leq \epsilon U(t)
=\epsilon b^\leftarrow(t)$, Potter's bounds yield
$$
  n t^{c\theta} \oF\Bigl( {\log (1+\eta)\over c\log t}t\bigr)
  = O(t^{\theta c+(1/\beta) -\alpha+(\epsilon/2)})
$$
as $t$ tends to infinity. Therefore, \ProofBdyThD\ is ultimately at most
$$
  \log t\bigl( -c(1+\epsilon^\beta)+2\eta c\epsilon E|X|+2(\alpha-1/\beta)
  \bigr) + t^{\theta c+1/\beta-\alpha+\epsilon} \, .
$$
We take $c$ large enough such that for any $\epsilon$ at most $1/2$, 
$$
  -c(1+\epsilon^\beta)+2\eta c\epsilon E|X| +2(\alpha-1/\beta) <-3/\beta
$$
say. Then, since $1/\beta-\alpha$ is negative, provided $\epsilon$ is small
enough, we can find $\theta$ such that
$$
  \theta c+{1\over\beta}-\alpha+\epsilon < 0 \, .
$$
It follows that for any $\epsilon$ small enough, \ProofBdyThD\ is 
ultimately at most $-2\beta^{-1}\log t$, and this implies \ProofBdyThC.

\bigskip

\noindent{\it Step 4.} We prove \StepFourAlt, that is,
$$\displaylines{\quad
  \limeps\limsupt {1\over (b^\leftarrow\oF)(t)}
  P\Bigl\{\, \exists n\geq b^\leftarrow (t)/\epsilon \, : \,
  S_n>t \,;\, 
  \hfill\cr\hfill
  \bcap_{i\geq 0}\{\, X_i\leq\theta ( b_i+t/\epsilon )\,\}\,\Bigr\}
  = 0 \, .
  \quad
  \equa{ProofBdyThE}
  \cr}
$$
On $\bcap_{i\geq 0}\{\, X_i\leq\theta ( b_i+t/\epsilon )\,\}$
$$
  S_n=\sum_{1\leq i\leq n} X_i\II\{\, X_i\leq \theta (b_i+t/\epsilon)\,\} \, .
$$
Using inequality \LargeDev\ with $\lambda$ of the form $c(t+b_n)^{-1}\log t$,
we bound the logarithm of
$$
  P\Bigl\{\, \sum_{1\leq i\leq n} X_i
             \II\{\, X_i\leq \theta (b_i+t/\epsilon)\,\} >t+b_n\,\Bigr\}
$$
by
$$\displaylines{\quad
  -c\log t +c(t+b_n)^{-1}(\log t)n\eta E|X|\bigl(1+o(1)\bigr)
  \hfill\cr\hfill
  {}+ \sum_{1\leq j\leq n} e^{c(t+b_n)^{-1}(\log t) \theta (b_j+t/\epsilon)}
  \oF\Bigl( {\log (1+\eta)\over c\log t} (t+b_n)\Bigr) \, .
  \quad
  \equa{ProofBdyThF}
  \cr}
$$
Since $n\geq b^\leftarrow (t)/\epsilon$, we have $t\lsim b(\epsilon n)$.
Furthemore, since $\limsup_{n\to\infty} n/b_n<\infty$,
$$
  \lim_{t\to\infty} \sup_{n\geq b^\leftarrow(t)/\epsilon} {n\over t+b_n}
  <\infty \, .
$$
Consequently, there exists a constant $A$ such that, refering to the second 
summand in \ProofBdyThF,
$$
  c\log t{ n\over t+b_n} \eta E|X|\bigl( 1+o(1)\bigr)
  \leq c\eta A \log t
$$
for any $n$ at least $b^\leftarrow(t)/\epsilon$ and any $t$ large enough.
Next, refering to the third summand in \ProofBdyThF, we have, using that
$t\lsim b(\epsilon n)\sim \epsilon^\beta b_n$ and $j\leq n$,
$$
  {b_j+t/\epsilon\over t+b_n}
  \lsim {b_n(1+\epsilon^{\beta-1})\over t+b_n}
  \leq 2
$$
for any $\epsilon$ small enough and $t$ large enough. Therefore, using again
Potter's bounds, the third summand in in \ProofBdyThF\ is ultimately in $t$ 
at most
$$
  n t^{2\theta c} (t+b_n)^{-\alpha+\epsilon}
$$
and, provided $\theta c$ and $\epsilon$ are small enough, 
tends to $0$ as $t$ tends to infinity.
It follows that provided $\theta c$ and $\epsilon$ are small enough, 
\ProofBdyThF\ is at most
$$
  \bigl(-c+c\eta A+o(1)\bigr) \log t
$$
ultimately. By taking $\eta$ smalle enough and $c$ large enough, 
\ProofBdyThE\ follows by Bonferroni's inequality.

\bigskip

\noindent {\it Step 5.} (S5) is proved with very much the same arguments 
as those used in step 4, and the details are omitted.

\bigskip

\noindent {\it Step 6.} Since \psiin\ holds, the same arguments as in step 4 
prove (S6).

\bigskip

\noindent{\it Step 7.} It is easy since the calculations done in step 1 
show that
$i(n,x)\sim n$ uniformly in the range $n\in I_{\delta,t}$ and $\theta\chi (t)
\leq x\leq \chi(t)/\theta$.

\bigskip

Therefore, applying Theorem \UsefulStatement, we proved 
Theorem \BdyTh.\hfill\qed

\bye